\documentclass[fleqn,12pt]{article}
\usepackage{amssymb,amsmath,amsfonts}
\usepackage{color}
\usepackage{graphicx}
\usepackage{latexsym}
\usepackage{amsmath}
\usepackage{amssymb}
\usepackage{graphics}
\usepackage[dvips]{epsfig}
\usepackage[margin=1in]{geometry}

\numberwithin{equation}{section}
\newtheorem{thm}{Theorem}[section]
\newtheorem{coro}[thm]{Corollary}
\newtheorem{lem}[thm]{Lemma}

\newtheorem{rmk}[thm]{Remark}
\newtheorem{defi}[thm]{Definition}

\newcommand{\ea}{\epsilon}
\newcommand{\ta}{\theta}

\newcommand{\da}{\delta}

\renewcommand{\aa}{\alpha}
\newcommand{\za}{\zeta}

\newcommand{\pl}{\partial}
\newcommand{\sa}{\sigma}

\newcommand{\ga}{\gamma}
\newcommand{\ba}{\beta}

\newcommand{\iy}{\infty}

\newcommand{\lt}{\left}
\newcommand{\rt}{\right}
\newcommand{\be}{\begin{equation}}
\newcommand{\bs}{\begin{split}}
\newcommand{\es}{\end{split}}
\newcommand{\ee}{\end{equation}}
\newcommand{\bee}{\begin{equation*}}
\newcommand{\eee}{\end{equation*}}

\newcommand{\ef}{\eqref}

\begin{document}
\begin{center}
\large{ \bf Well-Posedness  for the Motion  of  Physical Vacuum of  the Three-dimensional Compressible Euler Equations with or without Self-Gravitation}
\end{center}
\begin{center}
{Tao Luo, Zhouping Xin, Huihui Zeng}
\end{center}

\begin{abstract} This paper concerns the well-posedness theory of the motion of physical vacuum for the compressible Euler equations with or without self-gravitation.  First, a general uniqueness theorem of classical solutions is proved for the three dimensional general motion. Second, for the spherically symmetric motions, without imposing the compatibility condition of the first derivative being zero at the center of symmetry,  a new local-in-time existence theory is established in a functional space involving less derivatives than those constructed for three-dimensional motions in \cite{10',7,16'} by constructing suitable weights and cutoff functions featuring the behavior of solutions near both the center of the symmetry and the moving vacuum boundary.  \end{abstract}
\tableofcontents

\section{Introduction}
Due to its great physical importance and mathematical challenges, the motion of physical vacuum in compressible fluids has received much attention recently (cf. \cite{29,39,23, 24, 25}),  and significant progress has been made in particularly on the Euler equations (cf. \cite{7, 10, 10', 16, 16', 15}). Physical vacuum problems arise in many physical situations naturally, for example, in the study of the evolution and structure of gaseous stars (cf. \cite{6',cox}) for which vacuum boundaries are natural boundaries.  This paper is concerned with the evolving boundary of stars (the interface of fluids and vacuum states) in a compressible self-gravitating gas flow, which is modeled by
\be\label{2.1} \begin{split}
& \pl_t \rho  + {\rm div}(\rho {\bf u}) = 0 &  {\rm in}& \ \ \Omega(t), \\
 &  \pl_t (\rho {\bf u})  + {\rm div}(\rho {\bf u}\otimes {\bf u})+\nabla_{\bf x} p(\rho) = -\kappa \rho \nabla_{\bf x} \Psi  & {\rm in}& \ \ \Omega(t),\\
 &\rho>0 &{\rm in }  & \ \ \Omega(t),\\
 & \rho=0    &    {\rm on}& \  \ \Gamma(t):=\pl \Omega(t),\\
 &    \mathcal{V}(\Gamma(t))={\bf u}\cdot {\bf n}, & &\\
&(\rho,{\bf u})=(\rho_0, {\bf u}_0) & {\rm on} & \ \  \Omega:= \Omega(0).
 \end{split} \ee
Here $({\bf x},t)\in \mathbb{R}^3\times [0,\iy)$,  $\rho $, ${\bf u} $, $p$ and $\Psi$ denote, respectively, the space and time variable, density, velocity, pressure and gravitational potential given by
 \be\label{potential}\Psi({\bf x}, t)=-G\int_{\Omega(t)} \frac{\rho({\bf y}, t)}{|{\bf x}-{\bf y}|}d{\bf y} \ \ {\rm satisfying} \  \  \Delta \Psi=4\pi G \rho \  \ {\rm in} \ \  \Omega(t)
\ee
with  the gravitational constant $G$ taken to be unity; $\Omega(t)\subset \mathbb{R}^3$, $\Gamma(t)$, $\mathcal{V}(\Gamma(t))$ and ${\bf n}$ represent, respectively, the changing volume occupied by a fluid at time $t$, moving interface of fluids and vacuum states, normal velocity of $\Gamma(t)$ and exterior unit normal vector to $\Gamma(t)$; {$\kappa=1$ or $0$ corresponds to the Euler equations with or without self-gravitation}. We consider a polytropic star: the equation of state is then given by
\be\label{polytropic} p(\rho)=A \rho^{\gamma} \  \ {\rm for} \ \  \gamma>1 \ee
with  the adiabatic constant $A >0$  set to be unity. Let $c=\sqrt {p'(\rho)}$ be the sound speed,  the following condition:
\be\label{pv} -\infty<\nabla_{\bf n}(c^2)<0 \ \  {\rm on} \ \  \Gamma(t), \ee
defines a physical boundary (cf. \cite{7,10',16',23,24,25}).  Equations $\ef{2.1}_{1,2}$ describe the balance laws of mass and momentum, respectively; conditions $\ef{2.1}_{3,4}$ state that $\Gamma(t)$ is the interface to be investigated; $\ef{2.1}_5$ indicates that the interface moves with the normal component of the fluid velocity; and $\ef{2.1}_6$ is the initial conditions for the density, velocity and domain.

The physical vacuum appears in the stationary solutions of system \eqref{2.1} naturally.  Since for a stationary solution, one has
$$\nabla_{\bf x} p(\rho)=-\rho \nabla_{\bf x} \Psi,$$
which yields that in many
physical situations,
$$\nabla_{\bf N}(c^2)=-\frac {(\ga-1)}{\gamma}\nabla_{\bf N} \Psi \in (-\iy, 0)$$
on the interface,
where ${\bf N}$ is the exterior unit normal to the interface pointing from fluids to vacuum.  The physical vacuum makes the study of  free boundary problems of compressible fluids challenging and very interesting, because  standard methods of symmetric hyperbolic systems (cf. \cite{17}) do not apply directly. Recently, important progress has been made in the local-in-time well-posedness theory for the one- and three-dimensional compressible
Euler equations (cf. \cite{16, 10, 7, 10', 16'}). But
for the three-dimensional  compressible Euler-Poisson equations, the gravitational potential $\Psi$ defined by \eqref{potential} is non-local and depends on the unknown domain $\Omega(t)$. This will cause certain difficulties in the analysis.
Moreover, the self-gravitation leads to very rich and interesting physical phenomena for compressible fluids with vacuum (cf. \cite{rein, {17'}, luosmoller1, luosmoller2, liebyau, 6'}).

First, we address the uniqueness of classical solutions for the above free boundary problem. The uniqueness problem of free boundary problems for the  equations of compressible fluids
is subtle. This is particularly so in the presence of vacuum states.  For the physical vacuum free boundary problem of the 3-dimensional  compressible Euler equations, a uniqueness theorem is proved in \cite{10'}  in functional spaces which are smoother of one more degree than the spaces in which the existence theorems are established. This functional space in \cite{10} involves the weighted Sobolev norms of solutions.  In the present paper, we prove a general uniqueness theorem of classical solutions  for $1<\gamma\le 2$ (the most physically relevant regime) only requiring the derivatives appearing in the equations are continuous (indeed, we  can only require that the solutions are in $W^{1, \infty}$ in the whole domain and
$C^1$ in a neighborhood of the boundary).  {The strategy is to  extend the solutions of \ef{2.1} to those of Cauchy problems, for which the physical vacuum \eqref{pv} is crucial.}
Due to vacuum, the uniqueness of the extended solutions to the Cauchy problem is nontrivial because the standard symmetrization method of hyperbolic system does not apply in the presence of physical vacuum.  We use the relative entropy argument (cf. \cite{dafermos}) and potential estimates (cf. \cite{AB}). The advantage of the relative entropy argument is   making the full use of the nonlinear structure of the equations and requiring less regularity as realized by R. DiPerna (cf. \cite{diperna}).  The proof of the uniqueness theorem is  valid for both the compressible Euler-Poisson equations and the compressible Euler equations without self-gravitation. The above approach works for the case when $1<\gamma\le 2$. For the general case of $\gamma>1$, {we study the vacuum dynamics of free boundary problems of the compressible Euler equations without self-gravitation for  spherically symmetric motions, and prove the uniqueness theorem in the class of $C^1\cap W^{1, \infty}(\{{\bf x}\in\mathbb{R}^3: 0<|{\bf x}|\le R(t)\})$ without requiring that the solutions are differentiable  at the center of symmetry. Here the ball $B_{R(t)}$ is the moving domain. It should be noted that we do not require the vacuum boundary being physical in this case.}

We now turn to the existence problem. For a gaseous star,  it is important to consider  spherically symmetric motions since the stable equilibrium configurations are spherically symmetric which minimize the energy among all possible configurations (cf. \cite{liebyau}). As aforementioned, there have been some existence theories available for the  free boundary problems of the three-dimensional compressible Euler equations with physical vacuum (cf. \cite{10', 16'}). However, for  spherically symmetric motions, if the compatibility condition of the first derivatives of solutions being zero at the center of the symmetry is not imposed for the initial data, the initial data are $C^1$ only in the region excluding the origin as  3-spatial dimensional functions, but may not be differentiable at the origin. In this case,  the general existence theories in the three spatial-dimensions in  \cite{10', 16'} do not apply. Moreover, without imposing this compatibility condition at the center of symmetry, the coordinate singularity is very strong and the equation becomes very degenerate near the center of the symmetry.  Indeed, the initial density, $\rho_0$,
 appears as the coefficients in equation \eqref{e1-2'} in the Lagrangian coordinates.  This gives tremendous  difficulties when we use the elliptic estimates to estimate the derivatives in the direction normal to the boundary.  In those estimates, whether the first-order derivatives of the initial density is zero at the origin or not makes a
  distinct difference since we differentiate the system in the direction normal to the boundary in the elliptic estimates. We will choose deliberately a  cut-off function whose effective width depending on the initial density to capture more singular behavior of the solutions near the origin for the case of the non-zero first derivatives of the initial density.
The spherically symmetric solution we construct in this paper without imposing the above mentioned compatibility condition at origin is $C^1$ smooth only in the region excluding the origin, but $W^{1, \infty}$ in the region including the
origin  as the functions of 3-spatial dimensional functions. Therefore, the solution constructed in this paper is
different from those in \cite{16',10'} and exhibits some specially interesting features. For instance,  in the currently available  well-posedness theory for the  free boundary problems of the three-dimensional compressible Euler equations with physical vacuum, it requires by \cite{16'} or \cite{10'}  a weighted norm involving  $ \nabla_{\bf x}^{2\left[ {1}/({\gamma-1})\right]+9} {\bf u} |_{t=0}$  or $ \pl_t^{7+\lceil (2-\ga)/(\ga-1) \rceil}\nabla_{\bf x} {\bf u} |_{t=0} $   to be finite. However, for the three-dimensional spherically symmetric motion without imposing the compatibility condition of the first order derivatives of solutions being zero at the center of the symmetry,  we find in the present work a new  well-posedness theory with the initial data  less regular than those required in \cite{16',10'}.

As mentioned above, one of interesting features and  challenges in the exploration of spherically symmetric motions is to deal with the difficulty caused by the coordinates singularity at the origin (the center of the symmetry), besides the one caused by physical vacuum on the boundary. This is particularly so without imposing the compatibility condition of the first order derivatives of the solution being zero at the center of symmetry.  Indeed, in the well-posedness theory for spherically symmetric motions of viscous gaseous stars modeled by the compressible Navier-Stokes-Poisson equations with  vacuum boundary was shown in \cite{15},  a higher-order energy functional was constructed  which consists of two parts, called the Eulerian energy  near the origin expressed in Eulerian coordinates and the Lagrangian energy described in Lagrangian coordinates away from the origin. This indicates the subtlety of the behavior of solutions near the origin and vacuum boundary. In this paper, we find a uniform way to construct a higher-order energy functional only in Lagrangian coordinates by choosing suitable weights and cutoff functions which work for both the origin and the physical vacuum boundary of which the construction is elaborative.   It is noted such a strategy  works also for the compressible Navier-Stokes-Poisson equations .

 It should be noted here that the detailed proofs of the existence theorems in \cite{16',10'} are given for a initial flat domain of the form $\mathbb{T}^2 \times (0, 1)$, where $\mathbb{T}^2$ is a two-dimensional period box in $x_1$ and  $x_2$. Initially, the reference vacuum boundary is the top boundary
$\Gamma(0) =\{x_3 = 1\}$  while the bottom boundary $\{x_3 = 0\}$ is fixed. The moving vacuum boundary is given by $\Gamma(t) = \eta(t)(\Gamma(0))$ with the flow map $\eta(t)$. In principle, it would be feasible to extend flat domains to  general non-flat ones, for example, utilizing local coordinate charts and changes of coordinates to straighten out the boundary for each chart. However, it seems quite complicated and technically involved.  In this article,  we give a direct proof for non-flat initial domains (balls) of the existence theorem for the free boundary problem with physical vacuum. {It should be noted that the general approach we use here is motivated by \cite{16'}, in particular on the choice of the weights near the vacuum boundary.}

Before closing this introduction,  we would like to review some prior results on the free boundary problems besides the ones aforementioned. There has been a recent explosion of interests in the analysis of inviscid flows, one may refer to \cite{zhenlei,24,25,22,20,26, 29} for compressible motions and to \cite{1,9,19,21,34,40} for incompressible motions. Among these works, it should be mentioned that in \cite{24} a smooth existence theory ({for the sound speed $c$,} $c^\alpha$ is smooth across the interface with  $0 < \alpha \le 1$) was developed for the one-dimensional Euler equations with damping, based on the adaptation of the theory of symmetric hyperbolic systems  which is not applicable to physical vacuum boundary problems for which only $c^2$, the square of sound speed in stead of $c^{\alpha}$ ( $0 < \alpha \le 1$) , is required to be smooth across the interface); in \cite{zhenlei} the well-posedness of the physical vacuum free boundary problem is investigated for the one-dimensional Euler-Poisson equations,
using the methods motivated  by those in \cite{10} for the one-dimensional Euler equations; {existence and uniqueness for the three-dimensional compressible Euler equations modeling a liquid rather than a gas were established in \cite{22} by using Lagrangian variables combined with Nash-Moser iteration to construct solutions. For a compressible liquid, the density  is assumed to be a strictly positive constant on the moving boundary.  As such, the compressible liquid provides a uniformly hyperbolic, but characteristic, system. An alternative proof for the existence of a compressible liquid was given in \cite{35},  employing a solution strategy based on symmetric hyperbolic systems combined with Nash- Moser iteration.} As for viscous flows, there have been many results on the free-boundary Navier-Stokes equations  which cause quite different difficulties in analyses from that for  inviscid flows, so we do not discuss the works in that regime here.

The rest of this paper is organized as follows. In the next section,  we present and prove the uniqueness  of classical solutions to the three-dimensional physical vacuum problem \ef{2.1} when $1<\gamma\le 2$. The rest is devoted to the study of spherically symmetric motions. In Section 3, we formulate the three-dimensional spherically symmetric problem and state the main existence result.
Sections 4-8 are devoted to the case of $\gamma=2$. In Section 4, we describe a degenerate parabolic approximation to the original degenerate hyperbolic system. The uniform estimates for the higher-order energy functional are given in Sections 5-7: some preliminaries are presented in Section 5, the energy estimates in the tangential directions are given in Section 6, and the elliptic estimates in the normal direction for interior and boundary regions are presented respectively in Section 7.  With those
estimates, the existence  can be shown in Section 8. In sections 9 and 10, we will outline, but with enough details, the existence theory for the cases of $1<\gamma<2$ and $\gamma>2$, respectively. Section 11 is devoted to  the uniqueness theorem of classical solutions for the vacuum free-boundary problem of the compressible Euler equations without the self-gravitation in the spherical symmetry setting for all the values of $\gamma>1$, without assuming that the vacuum boundary is physical in the sense of \ef{pv}.

\section{Uniqueness for  three-dimensional Euler-Poisson equations with physical vacuum when $1<\gamma\le 2$.  }
For the three-dimensional free-boundary problem \ef{2.1} with a physical vacuum, we prove the following quite general uniqueness theorem for $1<\gamma\le 2$ in a natural functional space. It should be remarked  that the uniqueness theorems proved in \cite{10,10'} are in the functional spaces which are one more derivative smoother than the spaces  in which the existence theorems are established. Before stating the uniqueness theorem, we give a definition of classical solutions to problem \ef{2.1}.

\begin{defi}
A triple $(\rho, {\bf u}, \Omega(t))$ is called a classical solution to the physical vacuum free boundary
problem \ef{2.1} on $[0, T]$ for $T>0$  if the following conditions hold:

{\rm 1)}  $\Omega(t)=\cup_{k=1}^{m} \Omega^k(t)\subseteq \mathbb{R}^3$ is an open bounded set and $\partial\Omega(0)\in C^2$, where $\Omega^k(t)$ $(k=1,\cdots, m)$ are the connected component of  $\Omega(t)$  satisfying
\be\label{omega}
  \overline{\Omega^j (t) } \ \cap \ \overline{\Omega^k(t)}=\emptyset, \  \  1\le j\ne k \le m, \ \  t\in [0,T];
\ee

{\rm 2)} $(\rho, {\bf u})\in C^1 (\bar{ D})$ satisfies system \ef{2.1} and the physical vacuum condition:
\be\label{physical}
-\infty<\nabla_{\bf n}\lt(\rho^{\gamma-1}\rt)<0 \ \  {\rm on} \ \  \Gamma(t)=\partial\Omega(t), \ee
where ${\bf n}$ is the spatial unit outer norm to $\Gamma(t)$ and
$$D=\{({\bf x}, t): \ \  {\bf x}\in \Omega(t),\  \ t\in [0, T]\}, \ \ \bar{D}=D\cup \pl D.$$
\end{defi}
{Due to the regularities of the solution ${\bf u}\in C^1 (\bar{ D})$ and  $\partial\Omega(0)\in C^2$ in the definition above, we can see easily that
\be\label{oac2}
\bigcup_{0\le t\le T} \Gamma(t)=\bigcup_{0\le t\le T} \pl \Omega(t)=:\widetilde{\partial D} \in C^2.
\ee
Indeed, the interface $\Gamma(t)$ is moving with the fluids given by $ \mathcal{V}(\Gamma(t))={\bf u}\cdot {\bf n}$ on $\partial \Omega(t)$, where  $\mathcal{V}(\Gamma(t))$ is the normal velocity of $\Gamma(t)$; which is equivalent
to saying that $\widetilde{\partial D}$ is foliated by the integral curves of the vector fields $\partial_t +{\bf u}\cdot \nabla_x$.
}

The uniqueness theorem is as follows:

\begin{thm}\label{uniqueness1} {\rm (uniqueness  for the 3-d problem)} Suppose $1<\gamma\le 2$. Let $(\rho_1, {\bf u}_1, \Omega_1(t))$ and $(\rho_2, {\bf u}_2, \Omega_2(t))$ be two classical solutions to
problem \ef{2.1} on $[0, T]$ for $T>0$ in the sense of Definition 2.1, then for $t\in[0,T]$,
\be\label{uniq1}\begin{split}
\Omega_1(t)=\Omega_2(t)   \ \ {\rm and} \ \
 (\rho_1, {\bf u}_1)({\bf x}, t)=(\rho_2, {\bf u}_2)({\bf x}, t),  \ \  {\bf x}\in \Omega_1(t)=\Omega_2(t),
\end{split}\ee
provided that \ef{uniq1} holds for $t=0$.
\end{thm}

{\begin{rmk}\label{rmk2.2} It follows easily from the proof that the uniqueness result stated in Theorem \ref{uniqueness1} holds true for the solutions to \ef{2.1} as stated in definition 2.1 but with the regularity condition $\lt(\rho, {\bf u} \rt)\in C^1(\bar D)$ replaced by a less regular one:
 \be\label{regular'}
 \ (\rho, {\bf u})\in W^{1, \infty} (\bar{ D}) {~\rm and~} \ (\rho, {\bf u})\in C^1( {D}_{\delta} \cup \pl {D}_{\delta}  ),
\ee
where $D_{\delta}\subset D$ is a neighborhood of $\widetilde{\partial D}$.
\end{rmk}}

{\bf Proof of Theorem \ref{uniqueness1}}. The proof is divided into two steps. In step 1, we extend the solutions of \ef{2.1} to those of Cauchy problems. After that, we use the relative entropy argument and potential estimates to prove the uniqueness .

{\bf Step 1} (extension). Suppose that the triple $(\rho, {\bf u}, \Omega(t))$  is a classical solution   to
problem \ef{2.1} on $[0, T]$ in the sense of Definition 2.1.
 We will first extend the solution $(\rho,  {\bf u})$ from the domains $\Omega(t)$ to the whole domain $\mathbb{R}^3$ for
$t\in [0, T]$ such that the extended functions $(\widetilde{\rho},  \widetilde{{\bf u}})$  satisfy
\be\label{extend'}(\widetilde{\rho}, \widetilde{{\bf u}})({\bf x}, t)\in W^{1, \infty}(\mathbb{R}^3\times [0, T]), \ee
and solve the Euler-Poisson equations.

{\bf Step 1.1}. The extension of $\rho$ is clearly given by
\be\label{vacuumextend}
\widetilde{\rho}({\bf x}, t)=  {\rho}({\bf x}, t) \ \ {\rm in} \ \ D   \ \  {\rm and} \ \  \widetilde{\rho}({\bf x}, t)\equiv 0 \ \   {\rm in} \ \ \mathbb{R}^3\times[0,T]\setminus D .\ee
The extension of the vector field ${\bf u}$ is more complicated. In what follows, we extend it from $\Omega(t)$ to a neighborhood of $\Omega(t)$, and then to the rest region.

It follows from the condition \ef{omega} that there exists a small positive constant $\ea$ such that
\bee\label{x1}
\Omega_{\epsilon}^j(t)\cap \Omega_{\epsilon}^k(t)=\emptyset, \ \   1\le j\ne k \le m, \  \ t\in[0,T] ,\eee
where
\bee\label{4211}\Omega_{\epsilon}^j(t) = \Omega^j(t) \cup \{\bar {\bf x}+s {\bf n} (\bar {\bf x}, t): \ \ \bar {\bf x}\in \partial \Omega^j(t), \ \  0\le s\le \epsilon\}, \  \   1\le j \le m.\eee
Moreover, $\ea>0$ is chosen so small that the exponential map:
\be\label{mapin}\begin{split}
 \partial \Omega^j(t)\times [0, \epsilon] &\to \mathbb{R}^3 :
\ \ (\bar{\bf x}, s )& \mapsto  \bar {\bf x}+s {\bf n} (\bar {\bf x}, t)
\end{split}\ee
is injective for $1\le j\le m$ (that is, $\epsilon$ is less than the injectivity radius of $\partial \Omega^j(t)$). It should be noted that the number $\epsilon>0$ can be chosen uniformly for $t\in [0, T]$, because $ \widetilde{\partial  D}\in C^2$ (see \ef{oac2} for details).
Indeed, denote the second fundamental form of $\partial \Omega(t)$ by $\theta (\bar {\bf x}, t)$, then $\|\theta (\bar {\bf x}, t)\|_{C({\widetilde{\partial D}})}\le K_T$ for some positive constant $K_T$ which
may depends on $T$. Therefore, the injectivity radius of $\partial\Omega(t)$ has a positive lower bound for $t\in [0, T]$ (cf. \cite{CL00}).

Let $\eta\in C^{\infty}([0, \ \epsilon])$ be  a cut-off function satisfying
$$0\le \eta \le 1, \ \  \eta (s)=1 \  \ {\rm for} \ \ 0\le s\le \frac{\epsilon}{3}, \ \  \eta (s)=0 \  \ {\rm for}  \ \ \frac{2\epsilon}{3}\le s\le{\epsilon}.$$
For any ${\bf x} \in  \Omega_{\ea}^j(t) \setminus \Omega^j(t)$, define the extension of ${\bf u}$ as
\be\label{extensionforu}\begin{split}
  \widetilde{{\bf u}}\lt({\bf x},t\rt)
  =  \widetilde{{\bf u}}(\bar{\bf x}+s {\bf n} (\bar {\bf x}, t), t)
  =&\eta(s)\lt[{\bf u}(\bar {\bf x}, t)+s\nabla_{\bf x} {\bf u}(\bar{\bf x}, t)\cdot {\bf n}(\bar{\bf x}, t)\rt] \\
 =&\eta(s)\lt[{\bf u}(\bar {\bf x}, t)+\nabla_{\bf x} {\bf u}(\bar{\bf x}, t)\cdot
({\bf x}-\bar{\bf x})\rt], \ \ 0\le s \le \ea.
\end{split}\ee
So, we have extended the vector field ${\bf u}$ from $\Omega(t)$ to $\cup_{j=1}^m \Omega_\ea(t)=: \Omega_\ea (t)$, a neighborhood of $\Omega(t)$.
For the rest region, we simply define
\be\label{uu}
\widetilde{{\bf u}}({\bf x}, t)= {{\bf u}}({\bf x}, t) \  \ {\rm in }\ \ D,  \ \
\widetilde{{\bf u}}({\bf x}, t)={\bf 0} \ \  {\rm for} \ \  {\bf x}\in\mathbb{R}^3 \setminus \ \Omega_\ea(t) \ \  {\rm and}  \ \  t\in [0,T].
\ee

{\bf Step 1.2}. Next, we verify that the extended functions $(\widetilde{\rho}, \widetilde{{\bf u}})({\bf x}, t)$
defined on $\mathbb{R}^3\times [0, T]$ satisfy \eqref{extend'}. The key is the differentiability across the boundary $\widetilde{\partial D}:=\cup_{0\le t\le T} \pl \Omega(t)$.

Before doing so, some notations are needed. For any point $(\bar {\bf x}, \bar t)\in \widetilde{\partial D}$, let $(\tau_0, \tau_1, \tau_2)$ be a basis of the space-time tangent space of $\widetilde{\partial D}$ at $(\bar {\bf x}, \bar t)$ and ${\bf N}={\bf n}(\bar {\bf x}, \bar t)$ be the spatial unit outer normal to $\pl \Omega(\bar{t})$ at $\bar{\bf x}$. Then $(\tau_0, \tau_1, \tau_2, {\bf N})$ forms a basis of $\mathbb{R}^4$. So
$\nabla_{\tau_j}$ ($j=0, 1, 2$) and  $\nabla_{{\bf N}}$ determine all the derivatives $\partial_t$ and $\nabla_{\bf x}$ at the point $(\bar{\bf x}, \bar t)$.
For  $t\in [0, T]$, denote  the interior and exterior sides   of $\widetilde{\partial D}$ (or $\partial \Omega(t)$) by $\widetilde{\partial D}-$ (or $\partial \Omega(t)-$) and $\widetilde{\partial D}+$ (or $\partial \Omega(t)+$), respectively.

For $\widetilde{\rho}$, it follows from
$$\widetilde{\rho}\in C^1(\bar D) \ \ {\rm and} \ \   \widetilde{ \rho}=0 \ \ {\rm on} \  \  \mathbb{R}^3\times [0, T]\setminus D$$
that $\nabla_{\tau_i} \widetilde{\rho}=0$ on both $\widetilde{\partial D}-$ and $\widetilde{\partial D}+$ for  $i=0,1,2$; which implies that the tangential derivatives of $\widetilde{\rho}$ is continuous
across $\widetilde{\partial D}$. For the spatial normal derivative,  it follows from the physical vacuum
condition:
$$-\infty<\nabla_{\bf N}(\widetilde{\rho}^{\gamma-1})<0 \ \ {\rm  on }  \ \  \partial \Omega(t)- , $$
that
$$\nabla_{\bf N}(\widetilde{\rho})=0 \  \ {\rm if} \ \ 1<\ga<2   \ \ {\rm and}    \ \  -\infty<\nabla_{\bf N}(\widetilde{\rho})<0
\  \ {\rm if} \ \  \ga=2  \ \ {\rm on} \ \   {\partial \Omega(t)}-;$$
because of  $\widetilde{\rho}=0$ on $\widetilde{\partial D}$ and the fact
$$\nabla_{\bf N}(\widetilde{\rho})=\frac{1}{\gamma-1}\widetilde{\rho}^{2-\gamma}
\nabla_{\bf N}(\widetilde{\rho}^{\gamma-1}) $$
As on $\widetilde{\partial D}+$, it is easy to see that both the tangential and normal derivatives of $\widetilde{\rho}$ are zero due to $\widetilde{\rho}=0$  in $\mathbb{R}^3\times[0,T]\setminus D$. Thus, we have the following regularity of  $\widetilde{\rho}$:
\begin{equation}\label{extentionrho1} \begin{cases}
\widetilde{\rho}\in C^1\lt(\mathbb{R}^3\times [0, T]\rt)\cap W^{1, \infty}\lt(\mathbb{R}^3\times[0, T]\rt), \qquad {\rm if~}  \  1<\gamma<2, \\
\widetilde{\rho}\in C^1\lt(\overline {D}\rt)\cap C^1\lt(\overline {\mathbb{R}^3\times[0,T] \setminus D}\rt)\cap W^{1, \infty}\lt(\mathbb{R}^3\times[0, T]\rt), \qquad {\rm if~}  \ \gamma=2.
\end{cases}\end{equation}

For $\widetilde{{\bf u}}$, it follows from $\widetilde{{\bf u}}\in C^1(\bar {D})$ and \eqref{extensionforu} that
${\widetilde{\bf u}}$ is continuous across the interface $\widetilde{\pl D}$ which implies that the tangential derivatives of $\widetilde{{\bf u}}$ are  continuous across $\widetilde{\partial D}$; and that $\nabla_{{\bf N}}\widetilde{{\bf u}}$  is continuous across $\widetilde{\partial D}$. Therefore,  it holds that
\be\label{2.12} \widetilde{{\bf u}}\in C^1(\mathbb{R}^3\times [0, T]) \cap W^{1, \infty}(\mathbb{R}^3\times[0, T]). \ee

{\bf Step 1.3} We now verify that $(\widetilde{\rho}, \widetilde{{\bf u}})({\bf x}, t)$ solves the isentropic Euler-Poisson equations point-wisely. Note that
$$\widetilde{\rho} (\cdot, t) \in C^1\lt(\overline {\Omega(t)}\rt)\cap C(\mathbb{R}^3) \ \ {\rm and} \ \ \widetilde{\rho}\equiv 0 \  \ {\rm in} \ \   \mathbb{R}^3\setminus \Omega(t), \  \ t\in [0, T], $$
then we have, by the potential theory (cf. \cite{AB}), that for each $t\in [0, T]$,
\begin{equation}\label{psies}
\widetilde{\psi}(x, t)=-\int_{\Omega(t)}\frac{\widetilde{\rho}(y, t)}{|x-y|}dy=-\int_{\mathbb{R}^3}\frac{\widetilde{\rho}(y, t)}{|x-y|}dy\in C^1(\mathbb{R}^3)\cap W^{1, \infty}(\mathbb{R}^3). \end{equation}
In view of \ef{extentionrho1}, \ef{2.12} and \ef{psies}, we see that the extended functions $(\widetilde{\rho}, \widetilde{{\bf u}})$ solves the Euler-Poisson equations in $\mathbb{R}^3\times[0, T]\setminus \overline {D}$,
since $\widetilde{\rho}\equiv 0$ in this region.
As in $D$, by Definition 2.1,  $(\widetilde{\rho}, \widetilde{{\bf u}})$ of course solves the Euler-Poisson equations.

The remaining task is to verify this on $\widetilde{\partial D}$. Since
the vector field $\partial_t+\widetilde{{\bf u}}\cdot \nabla_{{\bf x}} $ is tangential to
$\widetilde{\partial D}$ and $\widetilde{\rho}=0$ on $\widetilde{\partial D}$, then
$$\lt(\partial_t+\widetilde{{\bf u}}\cdot \nabla_{{\bf x}} \rt) \widetilde{\rho} =0 \ \  {\rm on} \ \
\widetilde{\partial D}+ \ \ {\rm and} \ \ \widetilde{\partial D}-.  $$
It follows from \ef{2.12} and  $\widetilde{\rho}=0$ on $\widetilde{\partial D}$ that
\bee\label{xu}\widetilde{\rho} {\rm div}\widetilde{u}=0 \ \  {\rm on} \ \
\widetilde{\partial D}+ \ \ {\rm and} \ \ \widetilde{\partial D}-. \eee
Therefore, the equation of conservation of mass is verified. Similarly, we have
\begin{equation}\label{2.15}
\widetilde{\rho}(\pl_t+\widetilde{{\bf u}}\cdot \nabla_{{\bf x}} ) \widetilde{{\bf u}}=0 \  \  {\rm on~} \ \  {\widetilde{\partial D}}.
\end{equation}
Moreover, for any tangent vector $\tau$ to $\widetilde{\partial D}$, we have
\be\label{tangentiala}
\nabla_{\tau} P(\widetilde{\rho})\equiv 0 \  \ {\rm on~}  \ \widetilde{\partial D}, \ee
because of $\widetilde{\rho}\equiv 0$ on $\widetilde{\partial D}$.  For any spatial normal ${\bf N}$ to $\partial\Omega(t)$, it holds that
$$\nabla_{\bf N} P(\widetilde{\rho})=\nabla_{\bf N} (\widetilde{\rho}^{\gamma})=\frac{\gamma}{\gamma-1}\widetilde{\rho} \nabla_{\bf N} (\widetilde{\rho}^{\gamma-1})$$
and
$$-\infty<\nabla_{\bf N} (\widetilde{\rho}^{\gamma-1})<0 \  \ {\rm on} \ \  \widetilde{\partial D}- ,  \ \   \widetilde{\rho}=0 \  \ {\rm on} \ \ \widetilde{\partial D}.$$
Thus, we have
$$\nabla_{\bf N} P(\widetilde{\rho})=0 \ \  {\rm on}  \ \  \widetilde{\partial D}-.$$
This, together
with \ef{tangentiala}, verifies that
\be\label{tango}
\nabla_{\bf x} P(\widetilde{\rho})\equiv 0 \  \ {\rm on~} \ \ \widetilde{\partial D}-.  \ee
Since
$$\widetilde{\rho}\equiv 0 \  \ {\rm in} \ \  \mathbb{R}^3\times [0, T]\setminus D,$$
then
$$\nabla_{\bf N} P(\widetilde{\rho})=0 \  \ {\rm on} \ \   \widetilde{\partial D}+,$$
which together with \ef{tangentiala} implies that
\be\label{tango1}
\nabla_{\bf x} P(\widetilde{\rho})\equiv 0 \ {\rm on~} \widetilde{\partial D}+.  \ee
Therefore, it follows from \ef{2.15}, \ef{tango} and \ef{tango1} that the left-hand side of the equation of the balance law of the momentum is zero on $\widetilde{\partial D}$. On the other hand,
in view of  \ef{psies} and the fact that $\widetilde{\rho}\equiv 0$ on $\widetilde{\partial D}$,   the right-hand side  is also zero on $\widetilde{\partial D}$.

{\bf Step 2} (uniqueness).
Now, let $(\rho_1, {\bf u}_1, \Omega_1(t))$ and $(\rho_2, {\bf u}_2, \Omega_2(t))$ be two classical solutions of
problem \ef{2.1} on $[0, T]$ for $T>0$ in the sense of Definition 2.1. We extend those solutions as above by replacing $(\rho, {\bf u}, \Omega(t))$ by $(\rho_i, {\bf u}_i, \Omega_i(t))$ ($i=1, 2$), and denote these extended functions still by $(\rho_i, {\bf u}_i)$ ($i=1, 2$). It is easy to see that, for $i=1, 2,$
\be\label{gep1}\begin{split}
 & \pl_t \rho_i  + {\rm div}(\rho_i {\bf u}_i) = 0 &  {\rm in}& \ \ \mathbb{R}^3\times (0, T], \\
 &  \pl_t (\rho_i {\bf u}_i)  + {\rm div} (\rho_i {\bf u}_i\otimes {\bf u}_i)+\nabla_{\bf x} p(\rho_i) = -\kappa\rho_i \nabla_{\bf x} \Psi_i  & {\rm in}& \ \ \mathbb{R}^3\times (0, T],\\
 &\rho_i>0 &{\rm in }  & \ \ \Omega_i(t),\\
 & \rho_i=0    &    {\rm in}& \  \  \mathbb{R}^3\setminus \Omega_i(t),\\
 \end{split} \ee
 where
\be\label{potentiali}
\Psi_i ({\bf x}, t)=-\int_{\mathbb{R}^3}\frac{\rho_i({\bf y}, t)}{|{\bf x}-{\bf y}|}d{\bf y}, \quad
{\bf x}\in \mathbb{R}^3, \ \  t\in [0, T],\ee
$\kappa=0$ or $1$,
and \ef{extend'}, \ef{extentionrho1},  \ef{2.12}, \ef{psies} hold for $(\widetilde{\rho}, \widetilde{{\bf u}}, \Omega(t), \widetilde{\psi})$=$(\rho_i, {\bf u}_i, \Omega_i(t), \psi_i )$, $i=1, 2$. In what follows, we define the relative entropy-entropy flux pairs and derive some potential estimates.

{\bf Step 2.1}. For $i=1,2$, set
$${\bf u}_i=\lt(u_i^1, u_i^2, u_i^3\rt)^{\rm T}, \ \ {\bf m}_i=\lt(m_i^1, m_i^2, m_i^3\rt)^{\rm T} \ \ {\rm and } \ \ {\bf U}_i=\lt(U_i^0, U_i^1, U_i^2, U_i^3\rt)^{\rm T}, $$
where
$${ m}_i^j  =\rho_i {  u}_i^j , \  \ {  U}_i^0  =\rho_i ,\ \   {  U}_i^j  ={  m}_i^j, \ \  j=1,2,3.$$
Here and thereafter $(\cdot)^{\rm T}$ denotes the transpose. Equations  $\ef{gep1}_{1,2}$ can be written  as
\be\label{U}
 \pl_t {\bf U}_i + \sum_{j=1}^3 \pl_{{\bf x}_j} {\bf F}_j({\bf U}_i) = {\bf R}_i, \ \  i=1,2,
\ee
where ${\bf R}_i=\kappa(0, -\rho_i (\nabla_{\bf x} \Psi_i)^{\rm T})^{\rm T}$ and the flux functions ${\bf F}_j=({  F}_j^0,{  F}_j^1,{  F}_j^2, {  F}_j^3)^{\rm T}$ are given by
\bee
 {\bf F}_1({\bf U}_i ) =\left( {  m}_i^1, \  \frac{{  m}_i^1 {  m}_i^1}{\rho_i} + p(\rho_i) ,\
\frac{{  m}_i^1 { m}_i^2}{\rho_i},\ \frac{{  m}_i^1 {  m}_i^3}{\rho_i}\right)^{\mathrm T},
\eee
\bee
 {\bf F}_2({\bf U}_i ) =\left( {  m}_i^2, \  \frac{{  m}_i^1 {  m}_i^2}{\rho_i} ,\
\frac{{  m}_i^2 {  m}_i^2}{\rho_i}+ p(\rho_i) ,\ \frac{{  m}_i^2 {  m}_i^3}{\rho_i}\right)^{\mathrm T},
\eee
\bee
 {\bf F}_3({\bf U}_i ) =\left( {  m}_i^3, \  \frac{{  m}_i^1 {  m}_i^3}{\rho_i}  ,\
\frac{{  m}_i^2 {  m}_i^3}{\rho_i},\ \frac{{  m}_i^3 {  m}_i^3}{\rho_i}+ p(\rho_i)\right)^{\mathrm T}.
\eee
Denote the entropy $\eta$  and entropy flux function ${\bf q}=({ q}^1, {  q}^2,{  q}^3)^{\rm T}$  by
\be\label{eeflux} \begin{split}
 \eta( {\bf U}_i) =\frac{|{\bf m}_i|^2}{2\rho_i}+\frac{1}{\gamma-1}\rho_i^\gamma \ \ {\rm and} \ \  {\bf q} ({\bf U}_i) =\lt(\frac{|{\bf m}_i|^2}{2\rho_i}+\frac{\gamma}{\gamma-1}\rho_i^\gamma\rt)\frac{{\bf m}_i}{\rho_i}, \ \ i=1,2.
 \end{split} \ee
(For ${\bf x}\in \mathbb{R}^3\setminus \Omega_i(t)$ and $ t\in [0, T]$ where $\rho_i=0$, we set
$ \lt({\bf m}_i/ \rho_i\rt)({\bf x}, t)={\bf u}_i({\bf x}, t) $, $i=1,2$.) Then, we have
\bee
D \eta ({\bf U}_i)  D {\bf F}_j({\bf U}_i ) = D {  q}^j ({\bf U}_i), \ \ j=1,2,3, \ \  i=1, 2,
\eee
where
\bee
D \eta ({\bf U}_i)= \lt(\frac{\pl \eta( {\bf U}_i)}{\pl {  U}_i^0},\frac{\pl \eta( {\bf U}_i)}{\pl {  U}_i^1},\frac{\pl \eta( {\bf U}_i)}{\pl { U}_i^2},\frac{\pl \eta( {\bf U}_i)}{\pl {  U}_i^3}\rt),
\eee
\bee
D {  q}^j ({\bf U}_i) = \lt(\frac{\pl {  q}^j ( {\bf U}_i)}{\pl {  U}_i^0},\frac{\pl {  q}^j ( {\bf U}_i)}{\pl {  U}_i^1},\frac{\pl {  q}^j ( {\bf U}_i)}{\pl {  U}_i^2},\frac{\pl {  q}^j ( {\bf U}_i)}{\pl {  U}_i^3}\rt),
\eee
and  $ D {\bf F}_j({\bf U}_i ) $ represents the Jacobian matrix whose $(k,l)$ element is ${\pl { F}_j^k( {\bf U}_i)}/{\pl {  U}_i^l}$.
Easily, one can derive the equation for the entropy $\eta$  when ${\bf U}_i\in W^{1, \infty}$:
\be\label{etaeta}
\pl_t \eta({\bf U}_i) + \sum_{j=1}^3 \pl_{{\bf x}_j} { q}^j({\bf U}_i) +\kappa {\bf m}_i \cdot \nabla_{\bf x} \Psi_i =0, \ \ i=1,2.
\ee
We can therefore define the relative entropy-entropy flux pairs by
\bee\label{re}\begin{split}
 &\eta^*({\bf U}_1, {\bf U}_2)=\eta({\bf U}_2)-\eta({\bf U}_1)-D \eta({\bf U}_1) ({\bf U}_2-{\bf U}_1),  \\
& {  q}^{* j}({\bf U}_1, {\bf U}_2)={  q}^j({\bf U}_2)-{  q}^j({\bf U}_1)-D \eta({\bf U}_1) ({\bf F}_j({\bf U}_2)-{\bf F}_j({\bf U}_1)),  \ \  j=1, 2, 3. \end{split}\eee
where $\eta$ and ${\bf q}$ are defined by \ef{eeflux}. It follows from \ef{gep1}, \ef{U} and \ef{etaeta} that
\be\label{etastar}\begin{split}
 &\partial_t \eta^*+\sum_{j=1}^3 \partial_{{\bf x}_j}{  q}^{*j} \\
 =  & \lt[D\eta({\bf U}_2)- D \eta({\bf U}_1)\rt]{\bf R}_2
 -D^2 \eta({\bf U}_1) \lt( {\bf R}_1, {\bf U}_2- {\bf U}_1\rt)\\
 &-  \sum_{j=1}^3 D^2 \eta({\bf U}_1) \lt( \pl_{{\bf x}_j}{\bf U}_1 , {\bf F}_j({\bf U}_2)- {\bf F}_j({\bf U}_1) - D {\bf F}_j({\bf U}_1) ({\bf U}_2-{\bf U}_1)  \rt)\\
 = & \kappa \rho_2 ({\bf u}_1 -{\bf u}_2) \cdot \nabla_{{\bf x}}(\Psi_2-\Psi_1)\\
& -  \sum_{j=1}^3D^2 \eta({\bf U}_1) \lt( \pl_{{\bf x}_j}{\bf U}_1 , {\bf F}_j({\bf U}_2)- {\bf F}_j({\bf U}_1) - D {\bf F}_j({\bf U}_1) ({\bf U}_2-{\bf U}_1)  \rt),
\end{split}\ee
where
\begin{equation*}\label{}
 D^2 \eta({\bf U}_1) =\begin{pmatrix}
 &{|{\bf m}_1|^2}/(\rho_1)^3  + \ga (\rho_1)^{\ga-2}  &  -{m_1^1}/(\rho_1)^2 &  -{m_1^2}/(\rho_1)^2 &   -{m_1^3}/(\rho_1)^2\\
&-{m_1^1}/(\rho_1)^2 & {1}/{\rho_1} & 0 & 0\\
&-{m_1^2}/(\rho_1)^2 & 0 & {1}/{\rho_1} &
0\\
&-{m_1^3}/(\rho_1)^2 & 0 & 0 & {1}/{\rho_1}
\end{pmatrix}.\end{equation*}

{\bf Step 2.2}. Next, we will estimate the terms on the right-hand side of \ef{etastar}. Note that
\be\label{eta}\begin{split}
\eta^*=&\frac{1}{\ga-1}\lt[\rho_2^\ga - \rho_1^\ga -\ga \rho_1^{\ga-1}(\rho_2-
\rho_1)\rt]+\frac{1}{2}\rho_2|{\bf u}_2-{\bf u}_1|^2  \end{split}\ee
and
\begin{equation*}\label{122}\begin{split}
& \sum_{j=1}^3D^2 \eta({\bf U}_1) \lt( \pl_{{\bf x}_j}{\bf U}_1 , {\bf F}_j({\bf U}_2)- {\bf F}_j({\bf U}_1) - D {\bf F}_j({\bf U}_1) ({\bf U}_2-{\bf U}_1)  \rt)\\
 = &[p(\rho_2)-p(\rho_1)-p'(\rho_1)(\rho_2-\rho_1)]\sum_{j=1}^3
 \partial_{{\bf x}_j} u_1^j\notag
 +\frac{1}{2}\sum_{i, j=1}^3\rho_2(u_2^i-u_1^i)(u_2^j-u_1
^j)(\partial_{{\bf x}_j} u_1^i+\partial_{{\bf x}_i} u_1^j).
\end{split}\end{equation*}
Then, we have
\bee\label{}\begin{split}
 \lt|\sum_{j=1}^3D^2 \eta({\bf U}_1) \lt( \pl_{{\bf x}_j}{\bf U}_1 , {\bf F}_j({\bf U}_2)- {\bf F}_j({\bf U}_1) - D {\bf F}_j({\bf U}_1) ({\bf U}_2-{\bf U}_1)  \rt)\rt|
\le   C \lt\|  \nabla_{\bf x}  {\bf u}_1 (\cdot, t)   \rt\|_{L^\iy}  \eta^*,
\end{split}\eee
for  some constant $C>0$.
Therefore, we can integrate \ef{etastar} to get
\be\begin{split}\label{test5}
 \int_{\mathbb{R}^3} \eta^*({\bf x}, t)d{\bf x}
 \le & \int_{\mathbb{R}^3} \eta^*({\bf x}, 0)d{\bf x} +\kappa \int_0^t\int_{\mathbb{R}^3}\rho_2 ({\bf u}_1 -{\bf u}_2)\cdot\nabla_{{\bf x}}(\Psi_2-\Psi_1)d{\bf x} d\tau\\
&+C\sup_{0\le \tau \le t}||\nabla_{\bf x} {\bf u}_1(\cdot, \tau)||_{L^\infty}\int_0^t\int_{\mathbb{R}^3}\eta^*({\bf x}, \tau)d{\bf x}d\tau.
\end{split}\ee

To bound the second term on the right-hand side of \ef{test5}, we need a lemma presented in \cite{AB}: suppose $h\in L^{\infty}(\mathbb{R}^3)$ is a function having a compact support, then
$$
 \lt\|\nabla_{{\bf x}} \int_{\mathbb{R}^3}\frac{ {h}({\bf y})}{|{\bf x}-{\bf y}|}d{\bf y}\rt\|_{L^2(\mathbb{R}^3)}^2\le C \lt(\int_{\mathbb{R}^3}|h({\bf x})|^{4/3}d{\bf x}\rt)\lt(\int_{\mathbb{R}^3}|h({\bf x})|d{\bf x}\rt)^{2/3}<\iy, $$
where $C$ is a universal constant. By applying this fact and noting \eqref{potentiali}, we obtain
 \be\begin{split}\label{x10}
 & \int_{\mathbb{R}^3}\lt|\nabla_{{\bf x}}(\Psi_2-\Psi_1)({\bf x}, \tau)\rt|^2 d{\bf x} \\
 \le  & C \lt(\int_{\mathbb{R}^3} |\rho_2-\rho_1|^{4/3}({\bf x},\tau) d{\bf x}\rt)\lt(\int_{\mathbb{R}^3}|\rho_2-\rho_1| ({\bf x},\tau)d{\bf x}\rt)^{2/3} \\
 \le &C \lt(\int_{S(\tau)} |\rho_2-\rho_1|^{4/3}({\bf x},\tau) d{\bf x}\rt)\lt(\int_{S(\tau)}|\rho_2-\rho_1| ({\bf x},\tau)d{\bf x}\rt)^{2/3}
\end{split}\ee
 for any $\tau\in [0, T]$, where
$$S(\tau):= \{{\bf x}: |\rho_1-\rho_2|({\bf x}, \tau)>0\}, \ \  \tau\in [0, T].$$
By virtue of H${\rm \ddot{o}}$lder's inequality, one gets
$$\int_{S(\tau)}|\rho_2-\rho_1|^{4/3}({\bf x}, \tau)d{\bf x}\le \left(\int_{S(\tau)}|\rho_2-\rho_1|^2({\bf x}, \tau)d{\bf x}\right)^{2/3}\lt({\rm Vol} S(\tau)\rt)^{1/3}$$
and
$$\left(\int_{S(\tau)}|\rho_2-\rho_1|({\bf x}, \tau)d{\bf x}\right)^{2/3}\le \left(\int_{S(\tau)}|\rho_2-\rho_1|^2({\bf x}, \tau)d{\bf x}\right)^{1/3}\lt({\rm Vol} S(\tau)\rt)^{1/3}.$$
We thus achieve, using \ef{x10}, that
$$
 \int_{\mathbb{R}^3}\lt|\nabla_{{\bf x}}(\Psi_2-\Psi_1)({\bf x}, \tau)\rt|^2 d{\bf x}
 \le    C  \left(\int_{S(\tau)}|\rho_2-\rho_1|^2({\bf x}, \tau)d{\bf x}\right) \lt({\rm Vol} S(\tau)\rt)^{2/3}.
$$
Note from \ef{eta} that  for $1<\gamma\le 2$,
\be\label{hheta}\begin{split}
\eta^*({\bf x}, \tau)
\ge & C(\ga)\lt(\lt||\rho_2(\cdot, \tau)\rt||_{L^\infty}+\lt||\rho_1(\cdot, \tau)\rt||_{L^\infty}\rt)^{\gamma-2}(\rho_2-\rho_1)^2 +\frac{1}{2}\rho_2|{\bf u}_2-{\bf  u}_1|^2\ge 0. \end{split}\ee
Then, it yields that
 \bee\begin{split}\label{}
 & \int_{\mathbb{R}^3}\lt|\nabla_{{\bf x}}(\Psi_2-\Psi_1)({\bf x}, \tau)\rt|^2 d{\bf x}\\
 \le &   C  \lt(\lt||\rho_2(\cdot, \tau)\rt||_{L^\infty}+\lt||\rho_1(\cdot, \tau)\rt||_{L^\infty}\rt)^{2-\ga} \lt({\rm Vol} S(\tau)\rt)^{2/3}\int_{\mathbb{R}^3} \eta^*({\bf x} ,\tau) d{\bf x}.
\end{split}\eee
Using this and the Cauchy inequality, we have
 \be\label{finaa}\begin{split}
 & \lt|\int_{\mathbb{R}^3}\rho_2 ({\bf u}_1 -{\bf u}_2)\cdot\nabla_{{\bf x}}(\Psi_2-\Psi_1)d{\bf x}\rt|\\
 \le & \int_{\mathbb{R}^3}\rho_2 \lt|{\bf u}_1 -{\bf u}_2\rt|^2  d{\bf x} + \int_{\mathbb{R}^3}\rho_2 \lt|\nabla_{{\bf x}}(\Psi_2-\Psi_1)({\bf x}, \tau)\rt|^2 d{\bf x}\\
 \le &   C  (1+ { Z}(\tau))\int_{\mathbb{R}^3} \eta^*({\bf x} ,\tau) d{\bf x},
\end{split}\ee
where
 $$Z(\tau)=\lt||\rho_2(\cdot, \tau)\rt||_{L^\infty} \lt(\lt||\rho_2(\cdot, \tau)\rt||_{L^\infty}+\lt||\rho_1(\cdot, \tau)\rt||_{L^\infty}\rt)^{2-\ga} \lt({\rm Vol} S(\tau)\rt)^{2/3}.$$

Now, it follows from \ef{test5} and \ef{finaa} that for $t\in[0,T]$,
\bee\label{}\begin{split}
 \int_{\mathbb{R}^3} \eta^*({\bf x}, t)d{\bf x}
 \le &
 C\sup_{0\le \tau \le T}\lt(||\nabla_{\bf x} {\bf u}_1(\cdot, \tau)||_{L^\infty}+ Z(\tau)\rt)\int_0^t\int_{\mathbb{R}^3}\eta^*({\bf x}, \tau)d{\bf x}d\tau,
\end{split}\eee
when
$$\Omega_1(0)=\Omega_2(0)  \ \ {\rm and} \ \ (\rho_1, {\bf u}_1)(x,0)=(\rho_2, {\bf u}_2)(x,0).$$
So, one concludes from \ef{extend'}, \ef{hheta} and Grownwall's inequality  that
$$
 \int_{\mathbb{R}^3} \eta^*({\bf x}, t)d{\bf x}=0, \ \ ({\bf x},t)\in\mathbb{R}^3 \times [0,T],
$$
 and
  $$
 \rho_1({\bf x}, t)= \rho_2 ({\bf x}, t), \ \ ({\bf x},t)\in\mathbb{R}^3 \times [0,T].
$$
In particular,
$${\rm Spt} \rho_1(\cdot, t)= {\rm Spt}\rho_2(\cdot, t), \ \ t\in[0,T],$$
where
$${\rm Spt} \rho_i(\cdot, t)=\{{\bf x}\in \mathbb{R}^3:\ \rho_i({\bf x}, t)>0\}.$$
This implies that
$$\Omega_1(t)=\Omega_2(t), \ \  t\in [0,T].$$
In view of \ef{hheta} and $\ef{gep1}_{3,4}$, we then see that
$${\bf u}_1({\bf x}, t)={\bf u}_2({\bf x}, t),  \ \ ({\bf x},t)\in  \Omega_1(t) \times [0,T]. $$

This finishes the proof of Theorem \ref{uniqueness1}.

\section{Formulation and main existence results for spherically symmetric motions}
Starting from this section, we will focus on  spherically symmetric motions.
For a three-dimensional spherically symmetric motion, that is,
\begin{equation}\label{3.1}\rho({\bf x}, t) = \rho(r, t), \  \ {\bf u}({\bf x}, t) =  u(r, t) {\bf x} /r , \ \  {\rm where}  \ \ u\in  \mathbb{R} \  \ {\rm and}  \ \ r=|{\bf x}|,\end{equation}
system  \eqref{2.1}   can be written as follows: for $0\le t\le T$,
\begin{equation}\label{103}
\begin{split}
&\pl_t (r^2\rho)+\pl_r(r^2\rho u)=0  & {\rm in } & \  \  \lt(0, \  R(t)\rt) , \\
&\rho(\pl_t u+u\pl_r u)+\pl_r p+ {4\pi\rho}r^{-2}\int_0^r\rho(s,t) s^2ds=0 & {\rm in } & \  \ \lt(0, \  R(t)\rt), \\
& \rho>0  & {\rm in } & \  \  \lt[0, \  R(t)\rt)\\
& \rho(R(t),t)=0, \ \ u(0,t)=0, & & \\
& \dot R(t)=u(R(t), t) \ \ {\rm with} \ \ R(0)=1 , &    &\\
&  (\rho, u)= (\rho_0, u_0)   & {\rm on } & \ \ I:= (0, 1).
\end{split}
\end{equation}
Here $\ef{103}_{3,4}$ state that $r=R(t)$ is the free boundary and the center of the symmetry does not move; $\ef{103}_{5}$ describes that the free boundary issues from $r=1$ and moves with the fluid velocity; the initial conditions are prescribed in $\ef{103}_{6}$. The initial domain is taken to be a unit ball  $\{0\le r\le 1\}$. And the initial density of interest is supposed to satisfy
\be\label{156}
\rho_0(r)>0 \ \ {\rm for} \ \  0\le r<1  \ \ {\rm and} \ \ \rho_0(1)=0 ;\ee
and the  physical vacuum condition:
  \be\label{142}
  -\iy< \pl_r \lt(\rho_0^{\ga-1}\rt) <0 \  \ {\rm at} \  \ r=1.
   \ee

To fix the boundary, we transform system \ef{103} into Lagrangian variables.
Without abusing notations and for convenience, we use $x$ ($0\le x\le 1$) as the initial reference variable, and define the Lagrangian variable $r(x, t)$ by
\be\label{rxt}
\partial_t r(x, t)= u(r(x, t), t) \ \ {\rm for} \  \ t>0 \  \ {\rm and} \ \  r(x, 0)=x.
\ee
Thus $\eqref{103}_1$ implies that
$$\int_0^{r(x, t)}\rho(s, t)s^2ds=\int_0^x \rho_0(y) y^2dy.$$
Define the Lagrangian density and velocity by
$$f(x, t)=\rho(r(x, t), t) \ \ {\rm and} \ \  v(x, t)=u(r(x,t), t).$$
 Then the Lagrangian version of system \ef{103} can be written on the reference domain $I$ as
\be\label{419} \begin{split}
&\pl_t(r^2f) +r^2f(\pl_x v)/(\pl_x r)=0    & {\rm in}& \  \    I\times (0, T],\\
& f \pl_t v+ {\pl_x (f^{\gamma})}/(\pl_x r)+ {4\pi f} {r^{-2}}\int_0^x \rho_0(y) y^2dy=0 \ \  &{\rm in}& \  \   I\times (0, T] \\
&  f(1, t)=0, \ \ v(0, t)=0      &  {\rm on}& \  \ (0,T],\\
& (f, v) =(\rho_0, u_0)  &  {\rm on}& \  \   I \times \{t=0\}.
\end{split}
\ee
It follows from solving $\eqref{419}_1$ that
$$
f(x, t)=\left(\frac{x}{r}\right)^2\frac{\rho_0(x)}{\pl_x r(x, t)}.$$
So that system \eqref{419} can be rewritten as
\be\label{419'}\begin{split}
& \rho_0 \lt( \frac{x}{r}\rt)^2\pl_t v   + \pl_x \left[    \lt(\frac{x^2}{r^2}\frac{\rho_0}{\pl_x r}\rt)^\ga    \right] +4\pi \rho_0 \frac{x^2}{r^4}    \int_0^x \rho_0  y^2 dy =0 & {\rm in} & \ \ I \times (0,T],\\
& v(0, t)=0 & {\rm on} & \ \ \{x=0\}\times (0,T],\\
& v(x, 0)= u_0(x)  & {\rm on} & \ \ I \times \{t=0\},
\end{split}
\ee
where the initial density $\rho_0$ satisfying \ef{156} and \ef{142} has been viewed as a parameter.

With the  notations
$$ \sa (x) :=\rho_0^{\ga-1} x  \ \ {\rm and} \ \   \phi(x) :=  4\pi x^{-3}\int_0^x \rho_0  y^2 dy,$$
and the fact $r, \rho_0>0$ in $I \times (0,T]$, equation $\eqref{419'}_1$ can be rewritten as
\be\label{e1-2'}\begin{split}
 x \sa \pl_t v  + \pl_x \left[ \sa^2 \lt( \frac{x}{r}\rt)^{2\ga-2} \lt(\frac{1}{\pl_x r}\rt)^\ga    \right] - 2 \frac{\sa^2}{x}\lt( \frac{x}{r}\rt)^{2\ga-1} \lt(\frac{1}{\pl_x r} \rt)^{\ga -1}+ \phi \sa x^2 \lt(\frac{x}{r}\rt)^2 \\
 + \frac{  2-\ga }{\ga-1}
\sa x \pl_x \lt(\frac{\sa}{x}\rt)  \lt( \frac{x}{r}\rt)^{2\ga-2} \lt(\frac{1}{\pl_x r}\rt)^\ga   =0, \ \ {\rm in} \  \ I \times (0,T].
\end{split}\ee
As $\ga=2$, equation \ef{e1-2'} becomes relatively simple. However, it should be noted that the essential parts for $\ga=2$ and $\ga \neq 2$ are the same
(see equations \ef{eie} and \ef{gi1} later), so that the analysis for $\ga=2$ is applicable for general $\ga$. Therefore, we first present the main results for $\ga=2$ in the rest of this section, following the proof of the results we will then  discuss the case for general $\ga$ in Sections 9 and 10.

For $\ga=2$, we will consider a higher-order energy functional. To this end, we choose a  cut-off function $\za$ satisfying
$$
 \za=1 \  \ {\rm on} \ \ [0,\da], \ \ \za=0 \ \ {\rm on} \ \  [2\da, 1], \ \  |\za'|\le s_0/\da,
$$
for some constant $s_0$, where $\da=\da(\rho_0) $ is a small positive constant depending only on the initial density $\rho_0$ to be determined in Section 7.1.1.
The higher-order energy functional is   defined  to be
\be\label{norm}\begin{split}
  {E}(v, t):= & \lt\| \sa \pl_t^{4} v(\cdot,t) \rt\|_{1}^2 +\lt\|\pl_t^{4 } v (\cdot,t)\rt\|_{0}
  + \sum_{j=1}^2 \lt\{ \lt\|  \sa \pl_t^{4-2j  } v (\cdot,t)\rt\|_{j+1}^2
   +\lt\|   \pl_t^{4-2j  } v (\cdot,t)\rt\|_{j }^2 \rt.\\
   & \lt. +
    \lt\|  \frac{ \pl_t^{4-2j  } v  }{x} (\cdot,t)\rt\|_{j-1}^2
   +\lt\| \sa^{3/2} \pl_t^{5-2j  } \pl_x^{ j +1 }  v (\cdot,t)\rt\|_{0}^2
      +\lt\|  \sa^{1/2} \pl_t^{5-2j } \pl_x^{ j }  v (\cdot,t)\rt\|_{0}^2 \rt. \\
    &\lt.
    +   \lt\|  \pl_t^{5-2j  } v (\cdot,t)\rt\|_{j-1/2}^2 +\lt\| \frac{\pl_t^{5-2j } v }{x}(\cdot,t)\rt\|_{j -1}^2 \rt \}\\
   &  +\sum_{j=1}^2 \lt\{ \lt\|\za \sa \pl_t^{5-2j  } v (\cdot,t)\rt\|_{j+1}^2
   +\lt\|\za  \pl_t^{5-2j  } v (\cdot,t)\rt\|_{j }^2  \rt\} .
\end{split}\ee
Here and thereafter, we use  $\| \cdot\|_s$ to denote the norm of the  standard Sobolev space $\|\cdot\|_{H^s(I)}$ for $s\ge 0$; and define the polynomial function $M_0$ by
\be\label{559} M_0= P(E(v,0)),\ee
where $P$ denotes a generic polynomial function of its argument. Now, we are ready to state the main result.

\begin{thm}\label{existence1} {\rm (existence  for  $\gamma=2$)} Given initial data
$(\rho_0, u_0)$ such that $M_0 < \infty$, conditions \ef{156} and \ef{142} hold and $\rho_0\in C^3([0,1])$, there exists a solution $v(x,t)$ to  problem \eqref{419'} on $[0, T]$ for
$T >0$ taken sufficiently small such that
\be\label{612}\sup_{0\le t\le T}E(v, t)\le 2M_0.\ee
\end{thm}

{This section will be closed by several comments in order. First,  the time  derivatives of $v(x, t)$ at time $t=0$ involved in the definition of $M_0$ can be given in terms of the corresponding spatial derivatives of the initial data $\rho_0$ and $u_0$ due to the compatibility conditions of equation $\ef{419'}_1$. Second, the solution to the spherically symmetric problem \eqref{103} in Eulerian coordinates can be obtained  from the solution constructed in Theorem \ref{existence1}, since the Lagrangian variable $r\in H^2$ and $\pl_x r$ has a positive lower-bound. Finally, we can transform the solution of   problem \ef{103} back to solve the three-dimensional  problem \eqref{2.1} in $W^{1, \infty}(D_T)$, where
$$D_T=\{({\bf x}, t): \ {\bf x} \in \Omega(t), \ \ t\in [0, T]\} \  \ {\rm and} \ \ \Omega(t)=\{{\bf x}\in \mathbb{R}^3: \  |{\bf x}| < R(t)\}.$$
In fact, one can obtain a function
$\rho({\bf x}, t)$ and a vector field ${\bf u}({\bf x}, t)$ via \eqref{3.1} for $({\bf x}, t)\in D_T$ since $u({\bf 0}, t)=0$, and verify that $(\rho, {\bf u})\in C^1(D_T^0)\cap W^{1, \infty}(D_T) $ and \eqref{2.1} holds in $D_T^0$, where
$$\  D_T^0=D_T\setminus   \{{\bf 0}\}\times [0, T]. $$
However,  $(\rho, {\bf u})$ may not be in $C^1(D_T)$ if the compatibility condition of the first derivative being zero at the origin is not required.
}

\section{Parabolic approximations}
Let $\ga=2$ from this section to Section 8. Equation $\ef{419'}_1$ reads
\be\label{e1-3}\begin{split}
x \sa \pl_t v  + \left[ x^2 \sa^2 /(r^2 r'^2) \right]' -2 x^2 \sa^2 /(r^3 r') +  x^4 \phi \sa  /r^2  =0, \ \ {\rm in} \ \  I\times (0, T],
\end{split}\ee
where and in what follows, the notation $'$ denotes the $\pl_x$.
For $\mu>0$, we use the following degenerate parabolic problem to approximate \eqref{419'}:
\be\label{pe1-3}\ \begin{split}
&x \sa \pl_t v + \left[\sa^2 \frac{x^2 } {r^2   r'^2}  \right]' -2   \frac{\sa^2}{x} \frac{x^3}{r^3 r'} + x^2 \phi \sa  \frac{x^2} {r^2}  =\frac{2\mu}{x}\left[(x\sa)^2\left(\frac{v}{x}\right)'\right]'  &{\rm in} &\ \  I\times (0, T],\\
& v(0, t)=0   & {\rm on} &\  \ (0, T],\\
&v(x, 0)= u_0(x) & {\rm on} & \ \  I.
\end{split} \ee
As in \cite{10,10'}, one can show easily the existence and uniqueness of the solution $v_{\mu}$ to the above degenerate parabolic problem in a time interval $[0, T_{\mu}]$ with sufficient smoothness for which our later arguments are legitimate by smoothing the initial data and using the fixed point argument. Next, we will give the uniform estimates independent of $\mu$ to obtain the compactness of the sequence $\{v_{\mu}\}$ and a common time interval  $[0,  T]$ in which the problem \eqref{pe1-3} is solvable for any $\mu>0$, that is,
\begin{lem}\label{mainlem}
For any fixed $\mu>0$, let  $v_{\mu}$ be the smooth solution of \ef{pe1-3} in $[0, T_{\mu}]$. Then there exist constants $C>0$ and  $T\in (0,  T_{\mu}]$ independent of $\mu$  such that for the higher-order energy functional
\be
E(t):= E(v_{\mu}, t)
\ee
 defined in \ef{norm} satisfies the inequality
\be\label{mainestimate}
\sup_{t\in [0, T]}E(t)\le M_0+CT P\lt(\sup_{t\in [0, T]}E(t)\rt),
\ee
where $P(\cdot)$ denotes a generic polynomial function of its argument, and $M_0$ is defined in \ef{559}.
\end{lem}

We will establish the energy estimates in the tangential directions of the boundaries and the elliptic
estimates in the normal direction to prove this lemma. In what follows, for the sake of notational convenience, we omit $\mu$ in $v_{\mu}$, i.e., we denote $v_{\mu}$ by $v$ without ambiguity. Before  performing the detailed estimate, we list  some preliminaries which will be often used later.

\section{Some preliminaries}
In this section, we will present some embedding estimates for weighted Sobolev spaces, and derive some bounds which
follows directly from the definition of the  high order energy functional \ef{norm} and the a priori assumption.

{\bf Embedding of weighted Sobolev spaces}. Set
\be\label{distant}
d(x)=dist(x, \partial I)=\min\{x, 1-x\} \  \ {\rm for}   \ \  x\in I.
\ee
For any  $a>0$ and nonnegative integer $b$, the  weighted Sobolev space  $H^{a, b}(I)$ is given by
$$ H^{a, b}(I) := \lt\{d^{a/2}F\in L^2(I): \ \  \int  d^a|D^k F|^2dx<\infty, \ \  0\le k\le b\rt\}$$
  with the norm
$$ \|F\|^2_{H^{a, b}} := \sum_{k=0}^b \int d^a|D^k F|^2dx.$$
Here and thereafter, we use $\int dx := \int_I dx$  to  denote the spatial integral over the interval $I$.
Then  for $b\ge  {a}/{2}$, it holds the following embedding (cf. \cite{18'}):
 $$ H^{a, b}(I)\hookrightarrow H^{b- {a}/{2}}(I)$$
    with the estimate
  \be\label{wsv} \|F\|_{b- {a}/{2}} \le C \|F\|_{H^{a, b}} .\ee
In particular, we have
\be\label{weightedsobolev}
\|F\|^2_{1- a/2}\le C\int d(x)^a \lt(|F(x)|^2+|DF(x)|^2\rt)dx, \ \ a=1 \ \ {\rm or} \  \ 2.
\ee

{\bf Some consequences of \ef{norm}}. It follows from conditions \ef{156} and \ef{142} that $\sa(x)$ is equivalent to the distance function $d(x)$ defined in \ef{distant}.
Hence, the definition of the energy norm \ef{norm} and the embedding \ef{wsv} yield that
\be\label{odd}\begin{split}
\lt\| \lt(\sa \pl_t^3 v', \ \sa \pl_t v'' \rt )  (\cdot,t)\rt\|_{1/2}^2 + \lt\| \lt(\sa \pl_t^3 v\rt)(\cdot,t)\rt\|^2_{3/2} +  \lt\| \lt(\sa \pl_t v\rt)(\cdot,t) \rt \|_{5/2}^2    \le C {E(t)},
\end{split}\ee
Therefore, it holds that for any $p\in(1,\iy)$,
\be\label{egn}\begin{split}
  \lt\| \lt(x^{-1}v, \ v', \ \sa v'',\ x^{-1} \pl_t v, \ \sa\pl_t v',\  \ \pl_t^2 v, \ \sa\pl_t^2 v',   \sa \pl_t^3 v,   \sa \pl_t^4v
\rt)(\cdot,t)\rt\|_{L^\iy}
\\
  + \lt\| \lt(  \pl_t v',  \ \sa \pl_t v'', \ \pl_t^3 v, \ \sa \pl_t^3 v'
\rt)(\cdot,t)\rt\|_{L^p}\le C\sqrt{E(t)},
\end{split}\ee
where one has used the fact that in one space dimension, $\|\cdot\|_{L^\iy}\le C\|\cdot\|_1$ and $\|\cdot\|_{L^p}\le C\|\cdot\|_{1/2}$ ($1<p<\iy$).
Besides,
 another type of estimates  are also needed.
Noting from  \ef{norm}, \ef{odd}, and the simple fact that for any norm,
\bee\label{}\begin{split}
 \lt\|\pl_t^j v(\cdot, t)\rt\| = & \lt\|\pl_t^j v(\cdot, 0) + \int_0^t \pl_t^{j+1} v(\cdot,s)ds\rt\|   \\
 \le  &   \lt\|\pl_t^j v(\cdot, 0)\rt\| +     \int_0^t\lt\| \pl_t^{j+1} v(\cdot,s)\rt\| ds  \\
 \le & \lt\|\pl_t^j v(\cdot, 0)\rt\| +     t  \sup_{s\in [0,t]}\lt\| \pl_t^{j+1} v(\cdot,s)\rt\|, \ \ j=0,1,2,3;
\end{split}\eee
one can get
\be\label{tnorm}\begin{split}
&\sum_{j=0}^3 \lt\{\lt\|\lt(\sa \pl_t^j v \rt)(\cdot,t)\rt\|_{(5-j)/2}^2+\lt\| \pl_t^j v(\cdot,t)\rt\|_{(3-j)/2}^2\rt\}
 \\
&\ \ \ \ \ + \lt\|\lt(\frac{\pl_t^2 v}{x}\rt)(\cdot,t)\rt\|_0^2 +\lt\|\lt(\frac{\pl_t v}{x}\rt)(\cdot,t)\rt\|_0^2+\lt\|\lt(\frac{ v}{x}\rt)(\cdot,t)\rt\|_1^2
 \le  M_0 + CtP\lt(\sup_{[0,t]}E\rt);
\end{split}\ee
which implies in the same way as in the  derivation of \ef{egn} that for $p\in(1,\iy)$,
\be\label{tegn}\begin{split}
 &\lt\| \lt( \ x^{-1}   v, \ \sa  v',  \ \pl_t  v, \ \sa\pl_t  v',\ \sa \pl_t^2 v,
 \ \sa \pl_t^3 v
\rt)(\cdot,t)\rt\|_{L^\iy}^2 \\
& \ \  \ \ \ \ \ \ \ \ \
 \ \  \  \ \ \ \ \  \  \ \ \ \
  + \lt\| \lt(    v',  \ \sa   v'', \ \pl_t^2 v, \ \sa \pl_t^2 v'
\rt)(\cdot,t)\rt\|_{L^p}^2 \le M_0 + CtP\lt(\sup_{[0,t]}E\rt).
\end{split}\ee
{It should be noted that this paper concerns the local existence, so we always assume the time variable $t\le 1$.}

{\bf The a priori assumptions}.  Let $\mathcal{M}>0$ be a large constant (for instance, $\mathcal{M}=2{M}_0 +1$). Suppose that for $T\in (0,\mathcal{M}/2] $,
$$
\|v(\cdot,t)\|_2\le \mathcal{M}, \ \ t\in[0,T].
$$
Then it holds that for $(x,t)\in(0,1)\times[0,T]$,
\be\label{xr}
\frac{1}{2}\le \frac{r(x,t)}{x}\le \frac{3}{2}, \ \  \frac{1}{2}\le {r'(x,t)} \le \frac{3}{2}.
\ee
This can be achieved by noticing that $r(x,0)=x$ and
for any $(x,t)\in(0,1)\times(0,T]$,
\bee\begin{split}
\lt|\frac{r}{x}-1 \rt| =\lt| \int_0^t \frac{v(x,s)}{x}ds  \rt|=  \lt|\int_0^t \int_0^1 v'(\ta x,s)d\ta ds \rt|  \le t \sup_{s\in[0,t]}\|v'(s)\|_1\le \mathcal{M}T \le \frac{1}{2},
\end{split}\eee
\bee
\lt|r'-1\rt|=\lt|\int_0^t v'( x,s)ds \rt|  \le  t \sup_{s\in[0,t]}\|v'(s)\|_1\le\mathcal{ M} T \le \frac{1}{2}. \eee
In the proof of Lemma \ref{mainlem}, the time $t>0$ is taken sufficiently small so that the bounds \ef{xr} are always true.

 \section{Energy estimates}
The purpose of this section is to  derive a bound for
 $$\sup_{[0,t]}  \lt(\lt\| \sqrt{x \sa }  \pl_t^{5} v  \rt\|_0^2  +    \lt\| \sa \pl_t^{4} v\rt\|_1^2 +   \lt \|  \pl_t^{4} v\rt\|_0^2\rt).$$
 It should be noted that the estimate $ \|  \pl_t^{4} v\|_0$ is needed because the solution,we seek, satisfies $v(\cdot, t)\in C^1(I)$. By the Sobolev
 embedding, one needs to estimate $\|v(\cdot, t)\|_2$. Due to the degeneracy of the equation, one time derivative of the solution is equivalent to
 the half of the spatial derivative.

We first derive a general equation for time derivatives. Taking the $(k+1)$-th time derivative of equation $\ef{pe1-3}_1$ gives
\be\label{e1-11}\begin{split}
&x \sa \pl_t^{k+2} v   -2\left\{ \sa^2 \lt[   \frac{x^3} {r^3 r'^2} \frac{\pl_t^{k } v } {x}  +   \frac{x^2 } {r^2 r'^3}\pl_t^{k } v'\rt]\right\}'    + 2 \frac{\sa^2}{x}\lt[ 3 \frac{x^4} {r^4 r' } \frac{\pl_t^{k } v } {x}  +   \frac{x^3\pl_t^{k } v'}{r^3 r'^2} \rt] \\
= & \frac{2\mu}{x}\left[(x\sa)^2\left(\frac{\partial^{k+1}_t v}{x}\right)'\right]'+ 2\lt\{  \sa^2\lt[   I_{11}   +   I_{12}\right]\rt\}'- 2 \frac{\sa^2}{x}\lt[ 3 I_{21}   +  I_{22}  \rt]
-  \phi \sa x^2 \pl_t^{k+1} \lt(\frac{x^2}{r^2}\rt) ,
\end{split}\ee
where
\be\label{e1-6}\begin{split}
 &I_{11}=   \pl_t^{k } \lt( \frac{x^3   } {r^3 r'^2} \frac{v}{x}\rt)-  \frac{x^3} {r^3 r'^2} \frac{\pl_t^{k } v } {x} = \sum\limits_{\aa=0}^{k-1} C_\aa^{k-1} \pl_t^{k -\aa}\left( \frac{x^3}{r^3r'^2}\right)\lt(\frac{\pl_t^{\aa} v}{x}\rt)  ,\\
 &I_{12}=   \pl_t^{k } \lt( \frac{x^2  v' } {r^2 r'^3}\rt)-  \frac{x^2 } {r^2 r'^3}\pl_t^{k } v'
 =\sum\limits_{\aa=0}^{k-1} C_\aa^{k-1} \ \pl_t^{k -\aa}\left( \frac{x^2}{r^2r'^3}\right)(\pl_t^{\aa} v') ,\\
 &I_{21}=  \pl_t^{k } \lt( \frac{x^4   } {r^4 r' } \frac{v}{x}\rt)-  \frac{x^4} {r^4 r' } \frac{\pl_t^{k } v } {x} = \sum\limits_{\aa=0}^{k-1} C_\aa^{k-1} \pl_t^{k -\aa}\left( \frac{x^4}{r^4r' }\right)\lt(\frac{\pl_t^{\aa} v}{x}\rt)  ,\\
&I_{22}=  \pl_t^{k }\left( \frac{x^3v'}{r^3r'^2}\right)-\frac{x^3\pl_t^{k } v'}{r^3 r'^2} =\sum\limits_{\aa=0}^{k-1} C_\aa^{k-1} \pl_t^{k -\aa}\left( \frac{x^3}{r^3r'^2}\right)(\pl_t^{\aa} v').
 \end{split}\ee
Here and thereafter, $C_{\alpha}^{k-1}=(k-1)!/[(k-1-\alpha)!\alpha!]$.

Multiplying \ef{e1-11} with $k=4 $ by $\pl_t^{5} v$ and integrating the resulting equation  with respect to space and time yield, by virtue of integration by parts, that
\be\label{energy}\begin{split}
&\lt.\int \left\{ \frac{x \sa }{2} \left(\pl_t^{5} v \right)^2 + \frac{ x^2 }{r^2 r'} \left[ {\frac{1}{r'^2}(\sa \pl_t^{4} v')^2+ 3\frac{ x^2}{r^2} \lt(\frac{\sa}{x}\pl_t^{4} v\rt)^2 +2\frac{ x}{rr'}\lt(\frac{\sa}{x}\pl_t^{4} v\rt)(\sa\pl_t^{4} v')}   \right]         \right\} dx \rt|_0^t\\
& +2\mu\int_0^t\int    \lt[x\sa \lt(\frac{\pl_t^5 v}{x}\rt)'\rt]^2dxds\\
=& \int_0^t\int \left\{ \pl_t \left(\frac{x^2  }{r^2r'^3} \right)(\sa\pl_t^{4} v')^2 +3\pl_t \left(\frac{x^4  }{r^4r'}\right)  \lt(\frac{\sa}{x}\pl_t^{4} v \rt)^2  +2\pl_t\left(\frac{x^3 }{r^3r'^2}\right) \lt(\frac{\sa}{x}\pl_t^{4} v\rt)(\sa\pl_t^{4} v')    \right\}\\
&\times dxds - 2 \int_0^t\int   \lt[  \sa^2(I_{11}+I_{12})(\pl_t^{5} v')   +  (\sa/x) \sa(3I_{21}+I_{22}) (\pl_t^{5} v) \rt]dxds\\
&   -     \int_0^t \int  \phi   { x^2 \sa} \pl_t^{5}\lt(x^2/r^2\rt) (\pl_t^{5} v)   dxds\\
  =: &   J_1 -2 J_2 -  J_3 .
\end{split}\ee

In order to estimate the terms on the right-hand side of \ef{energy}, we notice that for all nonnegative integers $m$ and $n$,
\be\label{14}\begin{split}
 \lt|\pl_t^{k+1}\lt(\frac{x^m}{r^mr'^n}\rt)\rt| \le C \mathfrak{J}_k,\ \
 k=0,\cdots,4,
 \end{split}\ee
which follows from simple calculations and the a priori bounds \ef{xr}. Here
 \bee\label{}\begin{split}
&\mathfrak{J}_{0}= |{x^{-1}} { v}| + |  v'|,
 \ \
  \mathfrak{J}_{1}=|{x^{-1}} {\pl_t v}| + |\pl_t v'| +  \mathfrak{J}_{0}^2,
 \ \
 \mathfrak{J}_2=|{x^{-1}} {\pl_t^{2}v}  | +  | \pl_t^{2}v' |  +\mathfrak{J}_{1}\mathfrak{J}_{0},
\\
&\mathfrak{J}_{3}=|{x^{-1}} {\pl_t^{3}v}  | +  | \pl_t^{3}v' | + \mathfrak{J}_{2}\mathfrak{J}_{0}+\mathfrak{J}_{1}^2,
 \ \
 \mathfrak{J}_{4}=  |{x^{-1}} {\pl_t^{4 }v}  | +  | \pl_t^{4}v' | + \mathfrak{J}_{3}\mathfrak{J}_{0}+\mathfrak{J}_{2}\mathfrak{J}_{1}.
 \end{split}\eee
It follows from \ef{norm}, \ef{egn}, the H$\ddot{o}$lder inequality and $\|(\sa,\rho_0)\|_{L^\iy}\le  C$ that
\be\label{jk} \begin{split}
&\lt\|\mathfrak{J}_0\rt\|_{L^\iy}\le\lt\| {x^{-1}} {  v} \rt\|_{L^\iy} + \lt\|      v'\rt\|_{L^\iy}\le  C E^{1/2},\\
&\lt\|\mathfrak{J}_{1}\rt\|_{L^p}\le   \lt\| {x^{-1}} {\pl_t v} \rt\|_{L^\iy} + \lt\|\pl_t v'\rt\|_{L^p} + \lt\| \mathfrak{J}_{0} \rt\|_{L^\iy}^2  \le CP(E^{1/2}),\\
& \lt\|\mathfrak{J}_{2}\rt\|_{0}\le   \lt\| {x^{-1}} {\pl_t^2 v} \rt\|_0 + \lt\|\pl_t^2 v'\rt\|_0 + \lt\| \mathfrak{J}_{1} \rt\|_0 \lt\| \mathfrak{J}_{0} \rt\|_{L^\iy}   \le CP(E^{1/2}),
 \\
& \lt\|\sa  \mathfrak{J}_{2}\rt\|_{L^p}\le   C \lt\|   {\pl_t^2 v} \rt\|_{L^\iy} + \lt\|\sa \pl_t^2 v'\rt\|_{L^\iy} + C \lt\| \mathfrak{J}_{1} \rt\|_{L^p} \lt\| \mathfrak{J}_{0} \rt\|_{L^\iy}  \le CP(E^{1/2}),
 \\
&\lt\|\sa\mathfrak{J}_{3}\rt\|_{L^p}\le   C \lt\|   {\pl_t^3 v} \rt\|_{L^p} + \lt\|\sa \pl_t^3 v'\rt\|_{L^p} + \lt\| \sa\mathfrak{J}_{2} \rt\|_{L^p} \lt\| \mathfrak{J}_{0} \rt\|_{L^\iy}     + C \lt\| \mathfrak{J}_{1}\rt\|_{L^{2p}}^2  \le CP(E^{1/2}),\\
&\lt\|\sa\mathfrak{J}_{4}\rt\|_{0}\le  C \lt\|   {\pl_t^4 v} \rt\|_{0} + \lt\|\sa \pl_t^4 v'\rt\|_{0} + \lt\| \sa\mathfrak{J}_{3} \rt\|_0 \lt\| \mathfrak{J}_{0} \rt\|_{L^\iy} + \lt\| \sa\mathfrak{J}_{2} \rt\|_{L^4} \lt\| \mathfrak{J}_{1} \rt\|_{L^4} \le CP(E^{1/2}),
  \end{split}\ee
for any $p\in(1,\iy)$. Here and thereafter $P(\cdot)$ denotes a generic polynomial function. In particular, we have for $m\ge 1$ and $k=0,\cdots,4$,
\be\label{ik}\begin{split}
 \lt|\pl_t^{k+1}\lt(\frac{x^m}{r^m }\rt)\rt| \le C  \mathcal{I}_k \ \ {\rm satisfying} \  \ \lt\|x \mathcal{{I}}_{4}\rt\|_0  \le  CP(E^{1/2}),
 \end{split}\ee
where $\mathcal{I}_k$ equals $ {\mathfrak{J}}_k$ modular the terms involving spatial derivatives such as $\pl_t^i v'$ $(i=1,2,3,4)$. Similarly, one can use \ef{tnorm} and \ef{tegn}  to show that for $p\in(1,\iy)$,
\be\label{tjk} \begin{split}
&\lt\|\mathfrak{J}_0(t)\rt\|_{L^p}^2 +
 \lt\|\mathfrak{J}_{1}(t)\rt\|_{0}^2+  \lt\|\sa \mathfrak{J}_{2}(t)\rt\|_{L^p}^2+ \lt\|\sa\mathfrak{J}_{3}(t)\rt\|_{0}^2\le  M_0 + CtP\lt(\sup_{[0,t]}E\rt),
  \end{split}\ee
\be\label{pjk} \begin{split}
&\lt\|\mathcal{I}_0(t)\rt\|_{L^\iy}^2 +
 \lt\|\mathcal{I}_{1}(t)\rt\|_{0}^2+  \lt\|x \mathcal{I}_{1}(t)\rt\|_{L^\iy}^2+ \lt\|  \mathcal{I}_{2}(t)\rt\|_{0}^2 + \lt\|x \mathcal{I}_{3}(t)\rt\|_{0}^2\le  M_0 + CtP\lt(\sup_{[0,t]}E\rt) .
  \end{split}\ee

Next, we estimate the terms  on the right-hand side of \ef{energy}.  For $J_1$, it follows from \ef{14}, $\ef{jk}_1$ and the Cauchy inequality that
\be\label{J1}\begin{split}
  J_1
   \le & C
 \int_0^t \lt\{\lt\|\mathfrak{J}_{0}\rt\|_{L^\iy} \int  \lt[(\sa \pl_t^{4} v')^2
 +   \lt(\frac{\sa}{x}\pl_t^{4} v \rt)^2 \rt] dx \rt\}ds
 \le    C t \lt( \sup_{[0,t]} E ^{3/2}\rt).
\end{split}\ee
For $J_2$, an integration by parts leads to
\be\label{J2}\begin{split}
 J_2 = &
\lt. \int     \lt[ \sa  (I_{11}+I_{12})(\sa\pl_t^{4} v')   +   \sa (3I_{21}+I_{22}) \lt(\frac{\sa}{x}\pl_t^{4 } v\rt) \rt]dx \rt|_0^t \\
& -\int_0^t \int   \lt[ \sa (\pl_t  I_{11}+\pl_t I_{12})(\sa\pl_t^{4 } v')   +  \sa  (3\pl_t  I_{21}+\pl_t I_{22} ) \lt(\frac{\sa}{x}\pl_t^{4 } v\rt)  \rt]dxds\\
=:  & J_{21} - J_{22}.
\end{split}\ee
For $J_{22}$, noting from \ef{e1-6} and \ef{14} that
\bee\label{}\begin{split}
\lt| \pl_t I_{11}\rt|
= \sum\limits_{\aa=0}^{4}  C_{ \aa } \lt| \pl_t^{5-\aa}\left( \frac{x^3}{r^3r'^2}\right)\lt(\frac{\pl_t^{\aa} v}{x}\rt)\rt|
\le C \sum\limits_{\aa=0}^{4} \lt| \mathfrak{J}_{4-\aa} \lt(\frac{\pl_t^{\aa} v}{x}\rt)\rt| ;
\end{split}\eee
we can then obtain, using \ef{egn}, \ef{jk}  and the H${\rm \ddot{o}}$lder inequality, that
\be\label{j21}\begin{split}
 \lt\|\sa   \pl_t I_{11}\rt\|_0
\le & C   \lt\| \mathfrak{J}_{0}\rt\|_{L^\iy} \lt\|\pl_t^4 v \rt\|_0 +  C\lt\| \mathfrak{J}_{1}\rt\|_{L^4} \lt\|\pl_t^{3} v \rt\|_{L^4} + C\lt\| \mathfrak{J}_{2}\rt\|_{0} \lt\|\pl_t^{2} v \rt\|_{L^\iy} \\
&+  C\lt\| \sa \mathfrak{J}_{3}\rt\|_{0} \lt\|x^{-1}\pl_t  v \rt\|_{L^\iy}
+C\lt\|\sa\mathfrak{J}_{4}\rt\|_{0} \lt\|x^{-1}  v \rt\|_{L^\iy}\\
\le & C P(E^{1/2}) E^{1/2}.
\end{split}\ee
  Similarly, one can show that
\be\label{j22}\begin{split}
\lt\|\sa   \pl_t I_{21}\rt\|_0 \le C P(E^{1/2}) E^{1/2}.
\end{split}\ee
It follows from \ef{14}, \ef{egn}, \ef{jk} and the H${\rm \ddot{o}}$lder inequality that
\be\label{j23}\begin{split}
&\lt\|\sa   \pl_t I_{12}\rt\|_0 +\lt\|\sa \pl_t I_{22}\rt\|_0 \\
 \le & C   \lt\| \mathfrak{J}_{0}\rt\|_{L^\iy} \lt\|\sa\pl_t^{4} v' \rt\|_0 +  C\lt\| \mathfrak{J}_{1}\rt\|_{L^4} \lt\|\sa\pl_t^{3} v' \rt\|_{L^4} + C\lt\| \mathfrak{J}_{2}\rt\|_{0} \lt\|\sa\pl_t^{2} v' \rt\|_{L^\iy} \\
&+  C\lt\| \sa \mathfrak{J}_{3}\rt\|_{L^4} \lt\|\pl_t  v' \rt\|_{L^4}
+C\lt\|\sa\mathfrak{J}_{4}\rt\|_{0} \lt\|    v' \rt\|_{L^\iy}\\
\le & C P(E^{1/2}) E^{1/2}.
\end{split}\ee
Therefore, it follows from \ef{J2}-\ef{j23} and the H$\ddot{o}$lder inequality that
\be\label{J22}\begin{split}
\lt|J_{22}\rt | \le &  C\int_0^t  \lt[\lt(\lt\|\sa  \pl_t I_{11}\rt\|_0 + \lt\|\sa   \pl_t I_{12}\rt\|_0 \rt) \lt\|\sa  \pl_t^{4} v' \rt\|_0 \rt.\\
&+  \lt.\lt(\lt\|\sa  \pl_t I_{21}\rt\|_0 + \lt\|\sa    \pl_t I_{22}\rt\|_0 \rt) \lt\| (\sa/x) \pl_t^{4} v  \rt\|_0 \rt] ds \le C t P \lt(\sup\limits_{  [0,t]}E\rt).
\end{split}\ee
The term $J_{21}$  can be estimated as
\bee\label{}\begin{split}
  |J_{21}| \le &  M_0 + \ea \lt(\lt\|\sa \pl_t^{4} v'(t)\rt\|_0^2 + \lt\| (\sa/x) \pl_t^{4 } v(t)\rt\|_0^2  \rt) \\
  &+ C(\ea)  \lt\| \sa  (|I_{11}|+|I_{12}|+|I_{21}|+|I_{22}|)(t)\rt\|_0^2 \\
  \le & M_0 + \ea \lt(\lt\|\sa \pl_t^{4} v'(t)\rt\|_0^2 + \lt\|(\sa/x) \pl_t^{4 } v(t)\rt\|_0^2  \rt) \\
   &+ C(\ea) \sum_{\aa=0}^3 \lt( \lt\| \frac{\sa}{x} \mathfrak{J}_{3-\aa} (t)\pl_t^\aa v  (t)\rt\|_0^2 + \lt\| \sa \mathfrak{J}_{3-\aa}(t) \pl_t^\aa v'  (t)\rt\|_0^2  \rt),
\end{split}\eee
where  $\ea$ is a small positive constant to be determined later. Here we have used \ef{wsv}, \ef{14}, the Holder inequality and the Cauchy inequality.
 By virtue of \ef{norm}, \ef{egn}, \ef{jk} and \ef{tjk}, we obtain
\bee\label{}\begin{split}
&\lt\|\mathfrak{J}_0 (t) \pl_t^3 v (t) \rt\|_0+ \lt\| \mathfrak{J}_0 (t) \sa \pl_t^3 v'  (t) \rt\|_0 \\
= & \lt\| \mathfrak{J}_0 (t) \lt(\pl_t^3 v (0)+\int_0^t \pl_t^4 v (s)ds \rt)  \rt\|_0+ \lt\|\mathfrak{J}_0 (t) \lt(\sa \pl_t^3  v '(0)+\int_0^t \sa \pl_t^4 v'(s)ds \rt)   \rt\|_0\\
\le & \lt\| \mathfrak{J}_0 (t) \rt\|_{L^4} \lt(\lt\| \pl_t^3 v(0)  \rt\|_ {L^4}  +\lt\|\sa \pl_t^3 v'(0)  \rt\|_{L^4} \rt) \\
&+ \int_0^t   \lt( \lt\| \pl_t^4 v (s)  \rt\|_0+ \lt\|\sa \pl_t^4 v' (s) \rt\|_0\rt)  ds \lt\| \mathfrak{J}_0 (t)\rt\|_{L^\iy}  \\
\le & M_0 + CtP\lt(\sup_{[0,t]}E\rt).
\end{split}\eee
 Similarly, one can show that
\bee\label{}\begin{split}
&\sum_{\aa=0}^2 \lt( \lt\| \frac{\sa}{x} \mathfrak{J}_{3-\aa} (t)\pl_t^\aa v (t) \rt\|_0  + \lt\| \sa \mathfrak{J}_{3-\aa}(t) \pl_t^\aa v' (t) \rt\|_0   \rt)\\
\le & \lt\| \mathfrak{J}_1 \rt\|_{0} \lt(\lt\| \pl_t^2 v(0)  \rt\|_ {L^\iy}  +\lt\|\sa \pl_t^2 v'(0)  \rt\|_{L^\iy} \rt) + \int_0^t   \lt( \lt\| \pl_t^3 v (s)  \rt\|_{L^4}+ \lt\|\sa \pl_t^3 v' (s) \rt\|_{L^4}\rt)  ds \lt\| \mathfrak{J}_1 \rt\|_{L^4}\\
&+\lt\| \sa \mathfrak{J}_2 \rt\|_{L^4} \lt(\lt\|\frac{ \pl_t  v(0)} {x} \rt\|_ {L^4}  +\lt\|  \pl_t  v'(0)  \rt\|_{L^4} \rt) + \int_0^t   \lt( \lt\| \frac{\pl_t^2 v (s) }{x} \rt\|_0+ \lt\|\pl_t^2 v' (s) \rt\|_0\rt)  ds \lt\| \sa  \mathfrak{J}_2 \rt\|_{L^\iy}\\
&+\lt\| \sa \mathfrak{J}_3 \rt\|_{0} \lt(\lt\|\frac{    v(0)} {x} \rt\|_ {L^\iy}  +\lt\|    v'(0)  \rt\|_{L^\iy} \rt) + \int_0^t   \lt( \lt\|\frac{ \pl_t  v (s)  }{x} \rt\|_{L^\iy}+ \lt\|  \pl_t v' (s) \rt\|_{L^4}\rt)  ds \lt\|\sa \mathfrak{J}_3 \rt\|_{L^4}\\
\le & M_0 + CtP\lt(\sup_{[0,t]}E\rt).
\end{split}\eee
Therefore, we have arrived at
\be\label{J21}\begin{split}
  |J_{21}| \le  C(\ea)\lt[ M_0 + CtP\lt(\sup_{[0,t]}E\rt)\rt] + \ea \lt(\lt\|\sa \pl_t^{4} v'(t)\rt\|_0^2 + \lt\|(\sa/x) \pl_t^{4 } v(t)\rt\|_0^2  \rt)  .
\end{split}\ee
It remains to bound  $J_3$.   Note from \ef{ik}  that
\be\label{J3}\begin{split}
\lt|J_3\rt|  \le  &  \int_0^t \lt\| \phi\rt\|_{L^\iy}\lt\| x \pl_t^{5}\lt(x^2/ r^2\rt) (s) \rt\|_{0}  \lt\| x\sa  \pl_t^{5} v(s) \rt\|_0     ds \\
\le &  C \lt\|\rho_0\rt\|_{L^\iy} \lt(\sup_{s\in[0,t]}  \lt\| x\sa  \pl_t^5 v (s)\rt\|_0 \rt)   \int_0^t  \lt\| x \pl_t^5\lt(x^2/r^2\rt)(s) \rt\|_{0} ds   \\
\le &  C (\ea)   \lt(\int_0^t  \lt\| x \pl_t^5\lt(x^2/r^2\rt)(s) \rt\|_{0} ds \rt)^2+ \ea \lt(\sup_{[0,t]}  \lt\| x\sa  \pl_t^5 v \rt\|_0 \rt)^2 \\
 \le & C(\ea) t P\lt(\sup_{[0,t]} E\rt) + \ea  \sup_{[0,t]}  \lt\| x\sa  \pl_t^{5} v \rt\|_0^2 ,
\end{split}\ee
where $\ea>0$ is a small constant to be determined later.

In view of \ef{energy}, \ef{J1}, \ef{J2}, \ef{J22}, \ef{J21} and \ef{J3}, we see that
\bee\begin{split}
&\lt.\int \left\{ \frac{x \sa }{2} \left(\pl_t^{5} v \right)^2 + \frac{ x^2 }{r^2 r'} \left[ {\frac{1}{r'^2}(\sa \pl_t^{4} v')^2+ 3\frac{ x^2}{r^2} \lt(\frac{\sa}{x}\pl_t^{4} v\rt)^2 +2\frac{ x}{rr'}\lt(\frac{\sa}{x}\pl_t^{4} v\rt)(\sa\pl_t^{4} v')}   \right]         \right\} dx \rt|_0^t\\
& +2\mu\int_0^t\int    \lt[x\sa \lt(\frac{\pl_t^5 v}{x}\rt)'\rt]^2dxds\\
\le & C(\ea) \lt[ M_0+ C t P\lt(\sup_{[0,t]} E\rt) \rt]+ \ea \lt(\lt\|\sa \pl_t^{4} v'\rt\|_0^2 + \lt\| (\sa/x)\pl_t^{4 } v\rt\|_0^2  +  \sup_{[0,t]}  \lt\| x\sa  \pl_t^{5} v \rt\|_0^2 \rt).
\end{split}\eee
Since  $\lt\|  \left(\sqrt{x \sa } \pl_t^{5} v \right)(\cdot, 0)\rt\|_0^2$ can be bounded by $M_0$ due to \ef{e1-11} with $k=3$, and
\bee\begin{split}
&\frac{ x^2 }{r^2 r'} \left[ {\frac{1}{r'^2}(\sa \pl_t^{4} v')^2+ 3\frac{ x^2}{r^2} \lt(\frac{\sa}{x}\pl_t^{4} v\rt)^2 +2\frac{ x}{rr'}\lt(\frac{\sa}{x}\pl_t^{4} v\rt)(\sa\pl_t^{4} v')}   \right] \\
=&\frac{ x^2 }{r^2 r'} \left[ \frac{1}{2r'^2}(\sa \pl_t^{4} v')^2+ \frac{ x^2}{r^2} \lt(\frac{\sa}{x}\pl_t^{4} v\rt)^2 +\lt(\frac{1}{\sqrt{2}r' }(\sa \pl_t^{4} v')+\sqrt{2}\frac{ x }{r } \lt(\frac{\sa}{x}\pl_t^{4} v\rt) \rt)^2   \right] \\
\ge & \frac{ x^2 }{r^2 r'} \left[ \frac{1}{2r'^2}(\sa \pl_t^{4} v')^2+ \frac{ x^2}{r^2} \lt(\frac{\sa}{x}\pl_t^{4} v\rt)^2   \right]\\
\ge & C \lt[(\sa \pl_t^{4} v')^2+   (\sa x^{-1} \pl_t^{4} v)^2\rt],
\end{split}\eee
where the a priori lower bounds for $1/r'$ and $x/r$ were used;  then we have
\bee\begin{split}
 &\lt\| \sqrt{x \sa }  \pl_t^{5} v (t) \rt\|_0^2 +    \lt\| \sa \pl_t^{4} v'(t)\rt\|_0^2 +    \lt\| \sa x^{-1} \pl_t^{4} v(t)\rt\|_0^2 + \mu\int_0^t\int    \lt[x\sa \lt(\frac{\pl_t^5 v}{x}\rt)'\rt]^2dxds\\
 \le &  C(\ea) \lt[ M_0+ C t P\lt(\sup_{[0,t]} E\rt) \rt]+ C \ea \lt(\lt\|\sa \pl_t^{4} v'(t)\rt\|_0^2 + \lt\|\frac{\sa}{x}\pl_t^{4 } v(t)\rt\|_0^2  +  \sup_{[0,t]}  \lt\| x\sa  \pl_t^{5} v \rt\|_0^2 \rt),
\end{split}\eee
which implies, by choosing  $\ea$ suitably small, that
\be\begin{split}
 &\sup_{[0,t]}  \lt(\lt\| \sqrt{x \sa }  \pl_t^{5} v  \rt\|_0^2  +    \lt\| \sa \pl_t^{4} v'\rt\|_0^2 +   \lt \|  (\sa/x) \pl_t^{4} v\rt\|_0^2\rt)+\mu  \int_0^t \lt\|  x\sa \lt(\frac{\pl_t^5 v}{x}\rt)'(s)   \rt\|^2 ds \\
 \le &  M_0+ C t P\lt(\sup_{[0,t]} E\rt).
\end{split}\ee
The weighted Sobolev embedding \ef{weightedsobolev} implies
\bee\begin{split}
 \lt\| \pl_t^{4} v\rt\|_0^2  \le & C \lt(\lt\| \sa \pl_t^{4} v\rt\|_0^2 +   \lt\| \sa \pl_t^{4} v' \rt\|_0^2\rt) = C \lt(\lt\| x(\sa /x)   \pl_t^{4} v\rt\|_0^2 +   \lt\| \sa \pl_t^{4} v' \rt\|_0^2\rt)\notag\\
 \le & C \lt(\lt\| (\sa /x)   \pl_t^{4} v \rt\|_0^2 +   \lt\| \sa \pl_t^{4} v' \rt\|_0^2\rt),
\end{split}\eee
and we then obtain that
\bee\begin{split}
 &\sup_{[0,t]}  \lt(\lt\| \sqrt{x \sa }  \pl_t^{5} v  \rt\|_0^2  +    \lt\| \sa \pl_t^{4} v'\rt\|_0^2 +   \lt \|  \pl_t^{4} v\rt\|_0^2\rt)+\mu  \int_0^t \lt\|  x\sa \lt(\frac{\pl_t^5 v}{x}\rt)'(s)   \rt\|^2 ds \\
 \le& M_0+ C t P\lt(\sup_{[0,t]} E\rt),
\end{split}\eee
or equivalently
\be\label{et4}\begin{split}
 &\sup_{[0,t]}  \lt(\lt\| \sqrt{x \sa }  \pl_t^{5} v  \rt\|_0^2  +    \lt\| \sa \pl_t^{4} v\rt\|_1^2 +   \lt \|  \pl_t^{4} v\rt\|_0^2\rt)+\mu  \int_0^t \lt\|  x\sa \lt(\frac{\pl_t^5 v}{x}\rt)'(s)   \rt\|^2 ds \\
 \le& M_0+ C t P\lt(\sup_{[0,t]} E\rt).
\end{split}\ee

\section{Elliptic estimates}

In order to   estimate  the derivatives in the normal direction (the spatial derivatives in Lagrangian coordinates) which can not be obtained by  energy estimates as in the last section,
we employ the equation to perform the elliptic estimates.  Since the degeneracy of the equation near the origin $x=0$ and the boundary $x=1$ is of different orders, for example, in equation \ef{e1-3}, the coefficient of $\pl_t v$ is of the order $x^2$ as $x\to 0$, and of the order $(1-x)$ as $x\to 1$, we  separate the interior estimates and the estimates near the  boundary by choosing suitable cut-off functions.  To this end, we first identify the leading terms and lower order terms of the equation. Notice that
\be\label{}\begin{split}
- \left\{ \sa^2 \lt[   \frac{x^3} {r^3 r'^2} \frac{\pl_t^{k } v } {x}  +   \frac{x^2 } {r^2 r'^3}\pl_t^{k } v'\rt]\right\}'   +  \frac{\sa^2}{x}\lt[ 3 \frac{x^4} {r^4 r' } \frac{\pl_t^{k } v } {x}  +   \frac{x^3}{r^3 r'^2} \pl_t^{k } v'\rt]=-\sa(\mathfrak{{H}}_0+\mathfrak{H}_1+\mathfrak{H}_2),
\end{split}\ee
where
\be\label{h12}\begin{split}
\mathfrak{H}_0=&\sa\pl_t^k v''+\sa\lt(\frac{\pl_t^k v}{x}\rt)'+\lt[2\sa'-\frac{\sa}{x}\rt]\pl_t^k v'+\lt[2\sa'-3\frac{\sa}{x}\rt]\frac{\pl_t^k v}{x}
=  H_0 +4  \lt(\frac{\sa}{x}\rt)' \pl_t^k v,\\
\mathfrak{H}_1 = &  \lt\{
       2\sa' \lt( \frac{x^3  } {r^3 r'^2}-1\rt)
- \frac{3\sa  }{x } \lt(\frac{ x^4 } {r^4 r' }-1\rt) \rt\}\frac{\pl_t^{k } v}{x}  \\
& +\lt\{
       2\sa' \lt( \frac{x^2  } {r^2 r'^3}-1\rt)
- \frac{ \sa  }{x } \lt(\frac{ x^3 } {r^3 r'^2 }-1\rt) \rt\} {\pl_t^{k } v'} ,  \\
\mathfrak{H}_2= &  \sa  \lt[   \lt(\frac{x^3} {r^3 r'^2} \rt)'\frac{\pl_t^{k } v } {x}  +  \lt( \frac{x^2 } {r^2 r'^3}\rt)'\pl_t^{k } v' + \lt( \frac{x^3  } {r^3 r'^2}-1\rt)\lt(\frac{\pl_t^{k } v}{x}\rt)'+ \lt( \frac{x^2  } {r^2 r'^3}-1\rt) {\pl_t^{k } v''}\rt]
\end{split}\ee
and
\be\label{H0}
H_0=\sa  \pl_t^{k } v''+2\sa'  \pl_t^{k } v'  -2\sa'
  { \pl_t^{k } v }/{x}=\frac{1}{x\sa}\lt[(x\sa)^2\lt(\frac{\pl_t^k v}{x}\right)'\rt]'.
  \ee
We can then rewrite \ef{e1-11} as
\be\label{eie}\begin{split}
 H_0+\mu \pl_t H_0 = &
   \frac{1}{2}x   \pl_t^{k+2} v - 4  \lt(\frac{\sa}{x}\rt)' \pl_t^k v -  \mathfrak{H}_1- \mathfrak{H}_2 -\frac{1}{\sa} \left[  \sa^2  (I_{11}+I_{12}) \right]' \\
   & +    \lt(\frac{\sa}{x}\rt)  (3I_{21}+I_{22})
   + \frac{1}{2}\phi   { x^2  } \pl_t^{k+1}\lt(\frac{x^2}{r^2}\rt) =: \mathcal{G},
\end{split}\ee
where $I_{11}$, $I_{12}$, $I_{21}$ and $I_{22}$ are given by \ef{e1-6}.

 In order to obtain estimates independent of the regularization parameter $\mu$, we will also need the following lemma, whose proof can be found in \cite{10}:
\begin{lem}\label{lem1}  Let $\mu>0$ and $g\in L^{\infty}(0, T; H^s(I))$ be given,  and let   $f\in H^1(0, T; H^s(I))$ be such that
$$ f+\mu f_t=g, \qquad {\rm in~} (0, T)\times I.$$
Then
\be\label{lemma}
\|f\|_{L^{\infty}(0, T; H^s(I))}\le C \max \lt\{\|f(0)\|_s, \|g\|_{L^{\infty}(0, T; H^s(I))}\rt\}.\ee
\end{lem}
As an immediate consequence of \ef{eie} and \ef{lemma}, we  see that for any smooth function $\ba(x)$,
\be\label{ik3''}
\sup_{[0, t]}\|\beta H_0\|_0\le C\lt( \|\beta H_0(0)\|_0+\sup_{[0, t]} \|\beta \mathcal{G}\|_0\rt),\ee
\be\label{ik3'''}
\sup_{[0, t]}\|\beta H_0'\|_0\le C\lt( \|\beta H_0'(0)\|_0+\sup_{[0, t]} \|\beta \mathcal{G}'\|_0\rt).\ee
{Clearly, the weighted norm of $\pl_t^k v''$ (or $\pl_t^k v'''$)  can be derived from the corresponding weighted norm of $\pl_t^{k+2} v$ (or $\pl_t^{k+2} v'$).
Based on the energy estimate \ef{et4}, we can then obtain the estimates of $\pl_t^3 v''$ and $\pl_t^2 v''$ associated with weights. Furthermore, with the estimates of spatial derivatives of $\pl_t^3 v$ and $\pl_t^2 v$, one can get the weighted estimates of higher-order spatial derivatives of $\pl_t v$ and $v$.}

\subsection{Elliptic estimates -- Interior Estimates}
For the elliptic estimates, since the degeneracy of the equation near the origin $x=0$ and the boundary $x=1$ is of different orders, we will first choose a suitable cut-off function to separate the interior and boundary estimates. {The key  is to match the interior and boundary norms in the intermediate region.}

\subsubsection{Interior cut-off functions}
The interior cut-off function $\za(x)$  is chosen to satisfy
\be\label{delta}
 \za=1 \  \ {\rm on} \ \ [0,\da], \ \ \za=0 \ \ {\rm on} \ \  [2\da, 1], \ \  |\za'|\le s_0/\da,
\ee
for some constant $s_0$, where $\da$ is a constant to be chosen so that the estimates \ef{v''} and \ef{vt'''} below hold for all $k=0,1,2,3$. The choice of $\da$ will depend on the initial
density $\rho_0$.  Since
\bee
\sa'(x)=\rho_0(x)-x\rho_0'(x), \ \  \sa'(0)=\rho_0(0)>0,
\eee
there exists a constant $\da_0$ (depending only on $\rho_0(x)$) such that for all $ x\in[0,\da_0]$,
\be
m_0\le \rho_0(x)\le 3m_0, \ \ m_0\le \sa'(x)\le 3m_0, \ \ {\rm where} \ \ m_0=\rho_0(0)/2;
\ee
and then
\be
m_0x\le \sa(x)\le 3m_0 x, \ \ x\in [0,\da_0].
\ee
Set $m_1=\max_{0\le x\le \da_0} \lt\{  \lt|  \rho_0'(x)\rt| ,  \lt|  \rho_0''(x)\rt|  \rt\}$. Then for all $ x\in[0,\da_0]$,
\be
 \ \   \lt|\sa(x)/x-\sa'(x)\rt|   =\lt| x\rho_0'(x)\rt| \le m_1 x, \ \  \lt| \sa''(x)\rt|  \le 3m_1.
\ee

{\bf Analysis for $H_0$}. To this end, we rewrite $H_0$ as
\be\label{tl1}
H_0=\sa  f''+2\sa'  f'  -2\sa'
 \frac{ f }{x}, \ \ {\rm where} \ \ f=\pl_t^{k } v.
 \ee
Multiplying $H_0$ by the cut-off function $\zeta$ with $\da\in[0,\da_0/2]$, one may get
\bee\label{}\begin{split}
 \lt\| \zeta  {H}_{0 }\rt\|_0^2 =&
 \lt\| \za \sa   f'' \rt\|_0^2 + 4 \lt\| \zeta \sa'   f' \rt\|_0^2 + 4\lt\| \zeta  \sa'  \lt(
 \frac{ f }{x} \rt) \rt\|_0^2
  +4 \int \za \sa   f'' \zeta \sa'   f' dx
   \\
 &-4 \int  \za \sa   f'' \zeta  \sa'  \lt(
 \frac{ f }{x} \rt) dx  - 8 \int \zeta \sa'   f' \zeta  \sa'  \lt(
 \frac{ f }{x} \rt) dx.
 \end{split}\eee
Observing that
\bee\label{}\begin{split}
2 \int \za \sa   f'' \zeta \sa'   f' dx =&   -  \lt\| \zeta \sa'   f' \rt\|_0^2   -  \int \lt(  \za^2 \sa' \rt)' \sa \lt| f' \rt|^2 dx\\
\ge & -  \lt\| \zeta \sa'   f' \rt\|_0^2 - C(m_0,s_0) \int_\da^{2\da} \lt| f' \rt|^2 dx -C(m_0,m_1)  \da \lt\| \za   f'  \rt\|_0^2
 \end{split}\eee
and
\bee\label{}\begin{split}
    - \int  \za \sa   f'' \zeta  \sa'  \lt(
 \frac{ f }{x} \rt) dx
 = &   \int \lt(\za^2 \sa'\rt)' \sa \lt(
 \frac{ f }{x} \rt) f' dx +   \int  \za^2 \sa' \lt(\sa' -\frac{\sa}{x}\rt)    \lt(
 \frac{ f }{x} \rt) f' dx \\
 &+  \lt\| \zeta \sa'   f' \rt\|_0^2 +   \int \zeta^2 \sa' \lt( \frac{ \sa }{x} -\sa'\rt)     \lt| f' \rt|^2 dx\\
 \ge &   \lt\| \zeta \sa'   f' \rt\|_0^2  -C(m_0,s_0) \int_\da^{2\da} \lt(\lt| f' \rt|^2 + \lt|
 \frac{ f }{x} \rt|^2
 \rt) dx
 \\&-C(m_0, m_1)     \da \lt[ \lt\| \zeta     \lt(
 \frac{ f }{x} \rt) \rt\|_0^2  + \lt\| \zeta     f' \rt\|_0^2\rt],
 \end{split}\eee
 we have,  using the fact $\sa'(x)\ge m_0$ on $[0,2\da]$, that
\bee\label{}\begin{split}
 \lt\| \zeta  {H}_{0 }\rt\|_0^2 \ge &
 \lt\| \za \sa   f'' \rt\|_0^2 + \frac{2}{3} \lt\| \zeta \sa'   f' \rt\|_0^2 + \lt\| \zeta  \sa'  \lt(
 \frac{ f }{x} \rt) \rt\|_0^2\\
 & +\lt( \frac{16}{3} \lt\| \zeta \sa'   f' \rt\|_0^2 + 3\lt\| \zeta  \sa'  \lt(
 \frac{ f }{x} \rt) \rt\|_0^2 - 8 \int \zeta \sa'   f' \zeta  \sa'  \lt(
 \frac{ f }{x} \rt) dx\rt) \\
 &-C(m_0,s_0) \int_\da^{2\da} \lt( \lt| f' \rt|^2 + \lt|
 \frac{ f }{x} \rt|^2
 \rt) dx  -C(m_0, m_1)     \da \lt[ \lt\| \zeta     \lt(
 \frac{ f }{x} \rt) \rt\|_0^2  + \lt\| \zeta     f' \rt\|_0^2\rt]\\
  \ge &
 \lt\| \za \sa   f'' \rt\|_0^2 + \frac{2}{3} m_0^2 \lt\| \zeta     f' \rt\|_0^2 +  m_0^2 \lt\| \zeta    \lt(
 \frac{ f }{x} \rt) \rt\|_0^2-C(m_0,s_0) \int_\da^{2\da} \lt( \lt| f' \rt|^2 + \lt|
 \frac{ f }{x} \rt|^2
 \rt) dx\\
& -C(m_0, m_1)     \da \lt[ \lt\| \zeta     \lt(
 \frac{ f }{x} \rt) \rt\|_0^2  + \lt\| \zeta     f' \rt\|_0^2\rt].
 \end{split}\eee
Therefore, there exists a positive constant $\da_1=\da_1(m_0,m_1)$  such that if $\da \le \min\{\da_0/2, \da_1\}$,
\bee\label{}\begin{split}
 \lt\| \zeta  {H}_{0 }\rt\|_0^2 \ge &
 \lt\| \za \sa   f'' \rt\|_0^2 + \frac{1}{3} m_0^2 \lt\| \zeta     f' \rt\|_0^2 +  \frac{1}{2} m_0^2 \lt\| \zeta    \lt(
 \frac{ f }{x} \rt) \rt\|_0^2-C(m_0,s_0) \int_\da^{2\da} \lt( \lt| f' \rt|^2 + \lt|
 \frac{ f }{x} \rt|^2
 \rt) dx;
 \end{split}\eee
or equivalently
\be\label{v''}\begin{split}
     &\lt\| \za \sa    {\pl_t^{k } v''}    \rt\|_0^2+  \lt\| \za   {\pl_t^{k } v'}    \rt\|_0^2
    +   \lt\| \za   \lt( \frac{\pl_t^{k } v}{x}   \rt)  \rt\|_0^2\\
\le & C(m_0) \lt\|\za H_{ 0  }   \rt\|_0^2   + C(m_0,s_0) \int_\da^{2\da} \lt[ ({\pl_t^{k } v'})^2 +  \lt( \frac{\pl_t^{k } v}{x}   \rt)^2 \rt]dx.
 \end{split}\ee

{\bf Analysis for $H_0'$}. To estimate $H_0'$,  one needs also to compute the 1st spatial derivative of $H_0$. Clearly,
 \be\label{i10}\begin{split}
 H_0'+ 2\sa''\lt(\frac{f}{x}-f'\rt)= \sa   f'''  + 3 \sa'   f''   -2\sa'  \lt(
 \frac{ f }{x} \rt)'
 =:   \widetilde{H}_{0  }  , \ \ {\rm where} \ \ f=\pl_t^{k } v.
 \end{split}\ee
For any function $ \mathfrak{f}=\mathfrak{f}(x,t)$, it holds that
\be\label{formula}
\pl_x^j \mathfrak{f} =\pl_x^j \lt(x \frac{\mathfrak{f}}{x}  \rt) = x\pl_x^j \lt(\frac{\mathfrak{f}}{x}  \rt)  +j\pl_x^{j-1}\lt(\frac{\mathfrak{f}}{x}  \rt) , \ \  \ \  j=1,2,3;
\ee
so $\widetilde{H}_{0}$  can be rewritten as
\bee\label{}\begin{split}
  \widetilde{H}_{ 0  }    =  \sa   x g'' +3 \lt(\sa x\rt)'  g'
    +   4\sa'   g , \ \  {\rm where}
\ \  g=\lt(\frac{f}{x}  \rt)'=\lt(\frac{\pl_t^{k } v}{x}  \rt)'.
 \end{split}\eee
Thus,
\be\label{i1}\begin{split}
  {   \widetilde{ H}_{ 0  }-3  \lt(\sa'x-\sa\rt)  g'  }
    = &  \sa   x g'' +6   \sa  g'
    +   4   \sa'   g. \ \
 \end{split}\ee
Multiplying this equality  by the cut-off function $\za$ with $\da\in[0,\da_0]$ and taking the $L^2$-norm of the product yield
\be\label{hi1}\begin{split}
     \lt\|\za \widetilde{H}_{ 0  } +3\za \lt(\sa'x-\sa\rt)  g' \rt\|_0^2
     =
   \lt\| \za \sa   x g'' \rt\|_0^2+36 \lt\| \za \sa  g'\rt\|_0^2
    +  16 \lt\| \za \sa'   g \rt\|_0^2  \\
    + 12\int \za \sa   x g'' \za \sa  g' dx
      + 8 \int \za \sa   x g''    \za \sa'   g  dx + 48 \int   \za \sa  g'   \za \sa'   g   dx.
 \end{split}\ee
The last three terms on the right-hand side of \ef{hi1} can be bounded as follows:
\bee\label{hi2}\begin{split}
  -2 \int \za \sa   x g'' \za \sa  g' dx =   \int \lt( \za^2 \sa^2 x  \rt)'\lt|g'\rt|^2 dx =&  3 \lt\| \za \sa  g'\rt\|_0^2 + 2 \int   \za \za '  x  \lt|\sa g'\rt|^2 dx\\
  &+ 2 \int \za^2 \sa (\sa'x -\sa ) \lt| g'\rt|^2 dx,
 \end{split}\eee
\bee\label{hi3}\begin{split}
  & \int \za \sa   x g''    \za \sa'   g  dx =  \int \za^2  \sa  \lt(\sa' x -\sa \rt) g''    g  dx + \int \za^2  \sa^2    g''    g  dx \\
  =& \int \za^2  \sa  \lt(\sa' x -\sa \rt)     g g''  dx - 2 \int \za \za'  \sa^2    g  g' dx -2 \int \za^2    \sa \sa'  g  g' dx -   \lt\| \za \sa  g'\rt\|_0^2
 \end{split}\eee
and
\bee\label{hi4}\begin{split}
 -2 \int \za^2    \sa \sa'  g  g' dx =  \lt\| \za \sa'   g \rt\|_0^2 + 2 \int \za \za'    \sa \sa'  g^2   dx
 + \int \za^2     \sa \sa''  g^2   dx.
\end{split}\eee
It then follows from \ef{hi1} that
\bee\label{}\begin{split}
  & \lt\| \za \sa   x g'' \rt\|_0^2+10 \lt\| \za \sa  g'\rt\|_0^2 \\
= & \lt\|\za \widetilde{H}_{ 0  } +3\za \lt(\sa'x-\sa\rt)  g' \rt\|_0^2
       + 12\lt[ \int   \za \za '  x  \lt|\sa g'\rt|^2 dx +  \int \za^2 \sa (\sa'x -\sa ) \lt| g'\rt|^2 dx   \rt]  \\
&      - 8 \lt[ \int \za^2  \sa  \lt(\sa' x -\sa \rt)     g g''  dx - 2 \int \za \za'  \sa^2    g  g' dx  \rt]
      + 24 \lt[2 \int \za \za'    \sa \sa'  g^2   dx +\int \za^2     \sa \sa''  g^2   dx\rt]\\
\le & 2\lt\|\za \widetilde{H}_{ 0  }   \rt\|_0^2 + C(m_0,m_1) \da \lt[  \lt\| \za \sa   x g'' \rt\|_0^2 +   \lt\|  {\za}  \sa  g'\rt\|_0^2 +   \lt\|  {\za}     g \rt\|_0^2 \rt]+  C(m_0,s_0) \int_\da^{2\da} \lt[(\sa g')^2 + g^2   \rt] dx.
 \end{split}\eee
Therefore, there exists a constant $\da_2=\da_2(m_0,m_1)$ such that for   $\da\le  \min\{\da_0/2, \da_2\}$, it holds that
\be\label{hha}\begin{split}
  & \frac{1}{2}\lt\| \za \sa   x g'' \rt\|_0^2+ 5\lt\| \za \sa  g'\rt\|_0^2 \\
\le & 2\lt\|\za \widetilde{H}_{ 0  }   \rt\|_0^2  +C(m_0,m_1) \da \lt\|  {\za}     g \rt\|_0^2  +C(m_0,s_0) \int_\da^{2\da} \lt[(\sa g')^2 + g^2   \rt] dx.
 \end{split}\ee
To handle the term  $\lt\| \za    g \rt\|_0^2$, we need an additional estimate which follows from \ef{i1}, that is
\bee\label{}\begin{split}
     \lt\| \za \sa'   g \rt\|_0^2 \le
      & C(m_0) \lt[ \lt\|\za \widetilde{H}_{ 0 }   \rt\|_0^2 +\lt\| \za \sa   x g'' \rt\|_0^2+  \lt\| \za \sa  g'\rt\|_0^2 \rt] \\
  \le & C(m_0) \lt\|\za \widetilde{H}_{ 0 }   \rt\|_0^2  +C(m_0,m_1) \da \lt\|  {\za}     g \rt\|_0^2  +C(m_0,s_0) \int_\da^{2\da} \lt[(\sa g')^2 + g^2   \rt] dx ,
 \end{split}\eee
where we have used \ef{hha}. Hence, it holds that
\bee\label{}\begin{split}
  &  \lt\| \za \sa   x g'' \rt\|_0^2+  \lt\| \za \sa  g'\rt\|_0^2
    +   \lt\| \za \sa'   g \rt\|_0^2 \\
\le  &   C(m_0) \lt\|\za \widetilde{H}_{ 0 }   \rt\|_0^2  +C(m_0,m_1) \da \lt\|  {\za}     g \rt\|_0^2  +C(m_0,s_0) \int_\da^{2\da} \lt[(\sa g')^2 + g^2   \rt] dx .
 \end{split}\eee
Thus, there exists a constant $\da_3=\da_3(m_0,m_1)$ such that
\bee\label{}\begin{split}
  &  \lt\| \za \sa   x g'' \rt\|_0^2+  \lt\| \za \sa  g'\rt\|_0^2
    +   \frac{1}{2}\lt\| \za     g \rt\|_0^2
\le     C(m_0) \lt\|\za \widetilde{H}_{ 0 }   \rt\|_0^2   +C(m_0,s_0) \int_\da^{2\da} \lt[(\sa g')^2 + g^2   \rt] dx ,
 \end{split}\eee
provided $\da\le \min\{\sa_0/2, \da_2,\da_3\}$; where we have used the fact $\sa'(x)\ge m_0$ on $[0,\da_0]$.  It then follows from \ef{formula} and \ef{i10} that
\be\label{vt'''}\begin{split}
     &\lt\| \za \sa    {\pl_t^{k } v}   ''' \rt\|_0^2+  \lt\| \za  {\pl_t^{k } v}   ''\rt\|_0^2
    +   \lt\| \za    \lt( \frac{\pl_t^{k } v}{x}   \rt)' \rt\|_0^2 \\
\le & C(m_0)   \lt\|\za \widetilde{H}_{ 0 }   \rt\|_0^2     + C( m_0,s_0) \int_\da^{2\da} \lt[ \lt| {\pl_t^{k } v}   '' \rt|^2+ \lt|\lt( \frac{\pl_t^{k } v}{x}   \rt)'\rt|^2    \rt] dx \\
\le & C(m_0)   \lt\|\za  {H}_{ 0 }'   \rt\|_0^2     + C(m_0,m_1) \lt(
 \lt\| \za  {\pl_t^{k } v}   '\rt\|_0^2
    +   \lt\| \za     \frac{\pl_t^{k } v}{x}     \rt\|_0^2\rt)
\\
& +C( m_0,s_0) \int_\da^{2\da} \lt[ \lt|  {\pl_t^{k } v}    '' \rt|^2+ \lt|\lt( \frac{\pl_t^{k } v}{x}   \rt)'\rt|^2    \rt] dx.
 \end{split}\ee

{\bf A Choice of $\da$}. Choose
\be\label{da}\da=\min\{\da_0/2,\da_1,\da_2,\da_3\},\ee
then the estimates \ef{v''} and \ef{vt'''} hold for all $k=0,1,2,3$.
\subsubsection{Interior estimates for $\pl_t^3 v$ and $\pl_t^2 v$}
 Consider equation \ef{eie} with $k=3$, that is
\be\label{ik3}\begin{split}
 H_0+\mu \pl_t H_0 = &
   \frac{1}{2}x   \pl_t^{5} v - 4  \lt(\frac{\sa}{x}\rt)' \pl_t^3 v -  \mathfrak{H}_1- \mathfrak{H}_2 -\frac{1}{\sa} \left[  \sa^2  (I_{11}+I_{12}) \right]' \\
&    +    \lt(\frac{\sa}{x}\rt)  (3I_{21}+I_{22}) + \frac{1}{2} \phi   { x^2  } \pl_t^{4}\lt(\frac{x^2}{r^2}\rt).\\
\end{split}\ee
In order to bound $\lt\|\za H_0\rt\|$ by applying \ef{ik3''} with $\beta=\za$ given by \ef{delta}, we need to estimate the $L^2$-norm of the right-hand side of \ef{ik3} term by term. For this purpose, we first derive some estimates which will be used later. In addition to \ef{egn}, \ef{tnorm} and \ef{tegn},  we have some interior bounds:
\be\label{3egn}\begin{split}
 &  \lt\| \lt( \za \pl_t v',  \ \za \sa \pl_t v'', \ \za \pl_t^3 v, \ \za \sa \pl_t^3 v'
\rt)(\cdot,t)\rt\|_{L^\iy}\le C\sqrt{E(t)},\\
& \lt\|\lt(\za  {\pl_t^2 v} \rt)(\cdot,t)\rt\|_1^2+ \lt\|\lt(\za \sa  {\pl_t^2 v}, \za v \rt)(\cdot,t)\rt\|_2^2  +  \lt\|\lt(\za \sa  {  v} \rt)(\cdot,t)\rt\|_3^2
 \le  M_0 + CtP\lt(\sup_{[0,t]}E\rt),\\
&\lt\| \lt(   \za v',  \ \za \sa   v'', \ \za \pl_t^2 v, \ \za \sa \pl_t^2 v'
\rt)(\cdot,t)\rt\|_{L^\iy}^2 \le M_0 + CtP\lt(\sup_{[0,t]}E\rt);
\end{split}\ee
which implies
\bee\label{}\begin{split}
& \lt\|\za \mathfrak{J}_{1}\rt\|_{L^\iy}\le   \lt\|  {x^{-1}} {\pl_t v} \rt\|_{L^\iy} + \lt\|\za \pl_t v'\rt\|_{L^\iy} + \lt\| \mathfrak{J}_{0} \rt\|_{L^\iy}^2  \le CP(E^{1/2}),\\
&\lt\|\za \mathfrak{J}_0(t)\rt\|_{L^\iy}^2 \le 2 \lt(\lt\| {x^{-1}} {  v}  \rt\|_{L^\iy}^2 + \lt\|  \za    v' \rt\|_{L^\iy}^2 \rt)\le M_0 + CtP\lt(\sup_{[0,t]}E\rt),\\
 & \lt\|\za \mathfrak{J}_{2}(t)\rt\|_{0}^2
\le
C\lt(\lt\| {x^{-1}} {\pl_t^2 v}  \rt\|_0^2 + \lt\|\za \pl_t^2 v' \rt\|_0^2 + \lt\| \mathfrak{J}_{1}  \rt\|_0 ^2\lt\| \za \mathfrak{J}_{0}  \rt\|_{L^\iy}^2\rt)
 \le  M_0 + CtP\lt(\sup_{[0,t]}E\rt).
  \end{split}\eee
This, together with  \ef{jk} and \ef{tjk}, yields that for $p\in(1,\iy)$,
 \be\label{3jk} \begin{split}
 &\lt\|\mathfrak{J}_0\rt\|_{L^\iy}  +
 \lt\|\mathfrak{J}_{1}\rt\|_{L^p} + \lt\|\za \mathfrak{J}_{1}\rt\|_{L^\iy}+  \lt\|  \mathfrak{J}_{2}\rt\|_{0}  \le  CP(E^{1/2}) , \ \ \\
&\lt\|\mathfrak{J}_0(t)\rt\|_{L^p}^2+
 \lt\|\za \mathfrak{J}_{0}(t)\rt\|_{L^\iy}^2 +
 \lt\|\mathfrak{J}_{1}(t)\rt\|_{0}^2+  \lt\|\za \mathfrak{J}_{2}(t)\rt\|_{0}^2 \le  M_0 + CtP\lt(\sup_{[0,t]}E\rt) .
  \end{split}\ee
In a similar way as the derivation of \ef{14}, we have that for  nonnegative integers  $m$ and $n$,
\be\label{314}\begin{split}
 \lt|\pl_t^{k+1}\lt(\frac{x^m}{r^mr'^n}\rt)'\rt| \le C \mathfrak{L}_k,\ \
 k=0,1,2;
 \end{split}\ee
where
 \bee\label{}\begin{split}
& \mathfrak{L}_{0}  = \lt|v''\rt| + \lt| \lt(\frac{ v}{x}  \rt)'\rt|+ \mathcal{R}_0 \mathfrak{J}_{0}, \ \ {\rm with} \ \ \mathcal{R}_0 =\lt|\lt(\frac{r}{x}\rt)'\rt|+\lt|r''\rt| ,
 \\
&  \mathfrak{L}_{1}  =\lt|    \pl_t v''\rt|+ \lt| \lt(\frac{\pl_t v}{x}  \rt)'\rt| +   \mathcal{R}_0 \mathfrak{J}_{1}
+ \mathfrak{L}_{0}  \mathfrak{J}_{0},  \\
&   \mathfrak{L}_{2}  = \lt|\pl_t^2 v'' \rt|+ \lt| \lt(\frac{\pl_t^2 v}{x}  \rt)'\rt|
  + \mathcal{R}_0 \mathfrak{J}_{2} + \mathfrak{L}_{0}  \mathfrak{J}_{1}+\mathfrak{L}_{1}  \mathfrak{J}_{0}.
 \end{split}\eee
It can be checked (see the Appendix) that the following estimates hold:
\be\label{3lk} \begin{split}
&\lt\| \mathcal{R}_0 \rt\|_{0}  +
 \lt\|  \sa \mathcal{R}_0\rt\|_{L^\iy } \le C t \sup_{[0,t]} \sqrt{E} ,
  \\
&\lt\|  \mathfrak{L}_0\rt\|_{0}^2+ \lt\|\sa \mathfrak{L}_0\rt\|_{L^\iy }^2 + \lt\|  \sa \mathfrak{L}_{1}\rt\|_{L^p}^2+
 \lt\|\za \sa \mathfrak{L}_{1}\rt\|_{L^\iy}^2+  \lt\| \sa \mathfrak{L}_{2}\rt\|_{0}^2 \le  C P\lt(E(t)\rt) + C t P\lt(\sup_{[0,t]} E\rt),
 \\
&\lt\| \za \mathfrak{L}_0\rt\|_{0}^2+ \lt\| \sa \mathfrak{L}_0\rt\|_{L^p}^2 + \lt\|\za \sa \mathfrak{L}_0\rt\|_{L^\iy}^2 +
 \lt\|\sa \mathfrak{L}_{1}\rt\|_{0}^2+  \lt\|\za \sa \mathfrak{L}_{2}\rt\|_{0}^2 \le  M_0 + CtP\lt(\sup_{[0,t]}E\rt),
  \end{split}\ee
with $\|\cdot\|$ denoting $\| \cdot(t)\|$.

Next, we will bound $\|\za H_0 \|$ by the terms on the right-hand side of \ef{ik3}. It follows from \ef{tnorm}, \ef{et4} and the lower bound of $\rho_0(x)$ in the interior region  that
\be\label{t31}\begin{split}
  \lt\|
   \za\lt(\frac{1}{2}x   \pl_t^{5} v - 4  \lt(\frac{\sa}{x}\rt)' \pl_t^3 v\rt)(t)\rt\|_0^2 \le C \lt\| \za \sqrt{x \frac{\sa}{\rho_0}}   \pl_t^{5} v(t)\rt\|_0^2 +C  \lt\| \pl_t^3 v (t)\rt\|_0^2 \\
   \le   C \lt\|   \sqrt{x  {\sa} }  \pl_t^{5} v(t)\rt\|_0^2 +C  \lt\| \pl_t^3 v (t)\rt\|_0^2 \le
  M_0+ C t P\lt(\sup_{[0,t]} E\rt).
    \end{split}\ee
For $\mathfrak{H}_1$, noting from \ef{egn} that
  \be\label{r0}\begin{split}
 \lt\| \frac{ x }{r(x,t)  }-1 \rt\|_{L^\iy}   +
    \lt\| \frac{1 }{r'(x,t)}-1\rt\|_{L^\iy}
    \le     \lt\| \frac{ x }{r  } \lt(1-\frac{r}{x}\rt) \rt\|_{L^\iy}  +
    \lt\| \frac{1 }{r'}\lt(1-r'\rt)\rt\|_{L^\iy}  \\
    \le   C \int_0^t \lt(\lt\| \frac{v}{x}\rt\|_{L^\iy}+\lt\|v'\rt\|_{L^\iy}\rt)ds \le C t \lt(\sup_{[0,t]} \sqrt{E}\rt) ,
     \end{split}\ee
  we have
\be\label{t32}\begin{split}
  \lt\| \za \mathfrak{ H}_1 (t) \rt\|_0^2 \le   C  \lt\{ \lt\|     \frac{ x^4 } {r^4 r' }-1 \rt\|_{L^\iy}^2 +
    \lt\| \frac{x^3  } {r^3 r'^2}-1\rt\|_{L^\iy}^2 +\lt\|\frac{x^2  } {r^2 r'^3}-1\rt\|_{L^\iy}^2
 \rt\} \\
  \times \lt\{ \lt\|  \frac{\pl_t^{3 } v}{x} \rt\|_0^2 + \lt\|    \za \pl_t^{3 } v' \rt\|_0^2   \rt\}
 \le    C t P\lt(\sup_{[0,t]} E\rt).
    \end{split}\ee
For $\mathfrak{H}_2$, it follows from  $\ef{norm}$, \ef{r0} and $\ef{3lk}_1$ that
\be\label{t33}\begin{split}
  \lt\| \za \mathfrak{H}_2  (t)\rt\|_0^2 \le &\lt\{
    \lt\| \frac{x^3  } {r^3 r'^2}-1\rt\|_{L^\iy}^2 +\lt\|\frac{x^2  } {r^2 r'^3}-1\rt\|_{L^\iy}^2
 \rt\}  \lt\{ \lt\| \za  \sa \lt(\frac{\pl_t^{3 } v}{x} \rt)'\rt\|_0^2 + \lt\|  \za  \sa\pl_t^{3 } v'' \rt\|_0^2   \rt\}  \\
 &+ C\lt\| \sa \mathcal{R}_0 \rt\|_{L^\iy}^2 \lt(\lt\|   {\pl_t^{3} v }/{x}\rt\|_{0}^2  + \lt\|  \za \pl_t^{3} v'    \rt\|_{0}^2  \rt)
 \le     C t P\lt(\sup_{[0,t]} E\rt),
    \end{split}\ee
since
$$ \lt\| \za  \sa \lt(\frac{\pl_t^{3 } v}{x} \rt)'\rt\|_0
\le C \lt\| \za  x \lt(\frac{\pl_t^{3 } v}{x} \rt)'\rt\|_0
= C \lt\|  \za  \pl_t^{3 } v'-  \za \lt(\frac{\pl_t^{3 } v}{x} \rt)\rt\|_0.$$
Next, we will handle the terms involving $I_{11}$ and $I_{12}$ as follows,
\bee\label{}\begin{split}
 & \lt\| \za \frac{1}{\sa} \left[ \sa^2  (I_{11}+I_{12}) \right]'   \rt\|_0^2
  \le C  \lt\|\za (I_{11}+I_{12}) \rt\|_0^2 + C\lt\|\za \sa ( I_{11}+ I_{12})' \rt\|_0^2
  \\
  \le & C \sum_{\aa=0}^2     \lt\|\za   \mathfrak{J}_{2-\aa} \lt( \lt| {\pl_t^\aa v}  /x    \rt|+|\pl_t^\aa v'|    \rt)\rt\|_0^2  +C \sum_{\aa=0}^2      \lt\|\za  \sa \mathfrak{L}_{2-\aa} \lt( |\pl_t^\aa v/x|  +  | \pl_t^\aa v'|
  \rt) \rt\|_0^2 \\
  &+ C \sum_{\aa=0}^2   \lt\|\za   \mathfrak{J}_{2-\aa} \lt( \lt| \sa \lt( {\pl_t^\aa v}  /x \rt)'   \rt| + |\sa \pl_t^\aa v''|
  \rt)\rt\|_0^2   \\
  \le & C \sum_{\aa=0}^2   \lt\{   \lt\|\za   \mathfrak{J}_{2-\aa} \lt( \lt| {\pl_t^\aa v}  /x    \rt|+|\pl_t^\aa v'|  +  |\sa \pl_t^\aa v''|\rt)\rt\|_0^2  +      \lt\|\za  \sa \mathfrak{L}_{2-\aa} \lt( |\pl_t^\aa v/x|  +  | \pl_t^\aa v'|
  \rt) \rt\|_0^2 \rt\}.
 \end{split}\eee
Here we have used \ef{14} and \ef{314}. It follows  from \ef{norm}, \ef{egn}, \ef{tnorm}, $\ef{3egn}_1$ and $\ef{3jk}$ that
\be\label{abc}\begin{split}
 &\sum_{\aa=0}^2  \lt\|\za \mathfrak{J}_{2-\aa}(t)\lt( \lt| \pl_t^\aa v' (t)  \rt|+\lt| x^{-1}  {\pl_t^\aa v}  (t)    \rt| \rt) \rt\|_0
   \\
\le & \lt\| \za \mathfrak{J}_0 \rt\|_{L^\iy} \lt(\lt\| \frac{\pl_t^2 v(0) }{x} \rt\|_ {0}  +\lt\|  \pl_t^2 v'(0)  \rt\|_{0} \rt) + \int_0^t   \lt( \lt\|   \frac{\pl_t^3 v}{x}    \rt\|_{0}+ \lt\|\za \pl_t^3 v'  \rt\|_{0}\rt)  ds \lt\| \mathfrak{J}_0 \rt\|_{L^\iy}\\
&+\lt\|   \mathfrak{J}_1 \rt\|_{0} \lt(\lt\|\frac{ \pl_t  v(0)} {x} \rt\|_ {L^\iy}  +\lt\| \za \pl_t  v'(0)  \rt\|_{L^\iy} \rt) + \int_0^t   \lt( \lt\| \frac{\pl_t^2 v   }{x} \rt\|_0+ \lt\|\pl_t^2 v'   \rt\|_0\rt)  ds \lt\| \za  \mathfrak{J}_1 \rt\|_{L^\iy}\\
&+\lt\| \za \mathfrak{J}_2 \rt\|_0 \lt(\lt\| \frac{  v(0) }{x} \rt\|_ {L^\iy}  +\lt\|    v'(0)  \rt\|_{L^\iy}\rt) + \int_0^t   \lt( \lt\| \frac{\pl_t  v}{x}    \rt\|_{L^\iy}+ \lt\|\za \pl_t  v'  \rt\|_{L^\iy}\rt)  ds \lt\| \mathfrak{J}_2 \rt\|_{0}\\
\le & M_0 + CtP\lt(\sup_{[0,t]}E\rt).
\end{split}\ee
Similarly,
\bee\label{}\begin{split}
  \sum_{\aa=0}^2  \lt\|\za \mathfrak{J}_{2-\aa} (t)   \lt( \sa {\pl_t^\aa v''}\rt)(t)       \rt\|_0
\le   M_0 + CtP\lt(\sup_{[0,t]}E\rt);
\end{split}\eee
and
\bee\label{}\begin{split}
  \sum_{\aa=0}^2     \lt\|\za  \sa \mathfrak{L}_{2-\aa}(t) \lt( |\pl_t^\aa v/x|  +  | \pl_t^\aa v'|
  \rt)(t) \rt\|_0^2
   \le &   M_0 + CtP\lt(\sup_{[0,t]}E\rt).
 \end{split}\eee
Here we have used $\ef{3lk}_{2,3}$ to derive the last inequality. Therefore, it holds that
\be\label{t34}\begin{split}
 \lt\| \za \frac{1}{\sa} \left[ \sa^2  (I_{11}+I_{12}) \right]'(t) \rt\|_0^2 \le M_0 + CtP\lt(\sup_{[0,t]}E\rt).
 \end{split}\ee
In view of \ef{abc}, we obtain
 \be\label{t35}\begin{split}
    \lt\|\za \frac{\sa}{x}  (3I_{21}+I_{22})(t)\rt\|_0 ^2 \le &C \sum_{\aa=0}^2  \lt\|\za \mathfrak{J}_{2-\aa}(t)\lt( \lt| \pl_t^\aa v'   \rt|+\lt| \frac{ {\pl_t^\aa v}}{x}      \rt| \rt)(t) \rt\|_0\\
    \le & M_0+ C t P\lt(\sup_{[0,t]} E\rt).
    \end{split}\ee
Noting from \ef{ik} and \ef{pjk} that
\bee\label{phi3}\begin{split}
   \lt\|     \phi   { x^2  } \pl_t^{4}\lt(\frac{x^2}{r^2}\rt)(t)
 \rt\|_0^2
  \le & C \lt\|    x  {\mathcal{I}}_3(t)
 \rt\|_0^2
 \le   M_0+ C t P\lt(\sup_{[0,t]} E\rt),
    \end{split}\eee
one then derives from \ef{ik3''}, \ef{ik3}, \ef{t31}, \ef{t32}, \ef{t33}, \ef{t34}-\ef{t35} that
 \bee\label{}\begin{split}
\sup_{[0, t]}\lt\|\za H_0\rt\|_0^2 \le M_0+ C t P\lt(\sup_{[0,t]} E\rt).
\end{split}\eee
In view of \ef{v''} and \ef{tnorm}, we can therefore obtain, for any $s\in [0, t]$,
\bee\label{}\begin{split}
     & \lt\| \za \sa    {\pl_t^3 v''}   (s) \rt\|_0^2+  \lt\| \za   {\pl_t^3 v'}  (s)   \rt\|_0^2
    +   \lt\| \za    \frac{\pl_t^{3} v (s)}{x}     \rt\|_0^2  \\
\le & C(m_0) \lt\|\za H_{ 0  }  (s)  \rt\|_0^2    + C(m_0,s_0) \int_\da^{2\da} \lt[ ({\pl_t^{3} v'(s)})^2 +  \lt( \frac{\pl_t^{3 } v(s)}{x}   \rt)^2 \rt] dx\\
 \le & C \sup_{[0, t]}\lt\|\za H_0\rt\|_0^2  +  C \int_\da^{2\da} \lt[  (\sa {\pl_t^{3} v'(s)})^2
+  \lt(  {\pl_t^{3 } v(s)}    \rt)^2 \rt]  dx\\
\le &M_0+ C t P\lt(\sup_{[0,t]} E\rt) +M_0+ C s P\lt(\sup_{[0,s]} E\rt)\\
\le &M_0+ C t P\lt(\sup_{[0,t]} E\rt),
 \end{split}\eee
where we used the fact that $\sa(x)\ge m_0 \da $ on $[\da,2\da]$. This, together with \ef{tnorm}, implies that
\be\label{hk3}\begin{split}
      \sup_{[0, t]}\lt(\lt\| \za \sa    {\pl_t^3 v }    \rt\|_2^2+  \lt\| \za   {\pl_t^3 v }    \rt\|_1^2
    +   \lt\| \za   \lt( \frac{\pl_t^{3} v}{x}   \rt)  \rt\|_0^2\rt)
\le &M_0+ C t P\lt(\sup_{[0,t]} E\rt).
 \end{split}\ee
It follows from \ef{hk3} and \ef{norm} that
\be\label{iry2}\begin{split}
     \sup_{[0, t]}\lt( \lt\| \za \sa    {\pl_t^2 v }    \rt\|_2^2+  \lt\| \za   {\pl_t^2 v }    \rt\|_1^2
    +   \lt\| \za   \lt( \frac{\pl_t^{2} v}{x}   \rt)  \rt\|_0^2\rt)
\le &M_0+ C t P\lt(\sup_{[0,t]} E\rt).
 \end{split}\ee

 \subsubsection{Interior estimates for $\pl_t  v$ and $ v$.}
Consider \ef{eie} with  $k=1$. The basic idea is to   apply \ef{ik3'''} with $\beta=\za$.  As before, we first list some useful estimates here and then deal with
 $\lt\|\za H_0'\rt\|_0$ later.
Note that for all nonnegative integers $m$ and $n$,
\be\label{114}\begin{split}
 \lt|\pl_t \lt(\frac{x^m}{r^mr'^n}\rt)''\rt| \le C \mathcal{Q},
 \end{split}\ee
where
\bee\label{}\begin{split}
\mathcal{Q}= \lt|\lt(\frac{v}{x}\rt)''\rt|+ \lt|v'''\rt| + \mathcal{R}_0 \mathfrak{L}_0   +  \lt(\mathcal{R}_1 + \mathcal{R}_0^2 \rt)\mathfrak{J}_0 \ \ {\rm with} \ \  \mathcal{R}_1= \lt|r'''\rt|+ \lt|\lt(\frac{r}{x}\rt)''\rt|.
 \end{split}\eee
It follows from \ef{norm}, \ef{tnorm}, \ef{formula} and $\ef{3egn}_2$ that
\bee\label{}\begin{split}
   &\lt\| \sa \lt(\frac{v}{x}\rt)''(t)\rt\|_0+ \lt\|\sa v'''(t)\rt\|_{0}
  \le   C \lt\|v''- 2\lt(\frac{v}{x}\rt)'\rt\|_0+  \lt\| \sa v'''\rt\|_0 \le C \sqrt{E(t)},
    \\
   &\lt\|\za \sa \lt(\frac{v}{x}\rt)''(t)\rt\|_0^2+ \lt\|\za \sa v'''(t)\rt\|_0^2
  \le   2C \lt\|\za v''(t)- 2\za \lt(\frac{v}{x}\rt)'(t)\rt\|_0^2+  2\lt\|\za \sa v'''(t)\rt\|_0^2\\
    &\qquad \qquad\qquad\qquad \qquad\qquad  \quad\le     M_0+ C t P\lt(\sup_{[0,t]} E\rt).
    \end{split}\eee
We then have, by $\ef{3jk}_1$ and $\ef{3lk}_{1,2}$, that
\be\label{rrr}\begin{split}
   \lt\|\sa \mathcal{R}_1(t)\rt\|_{0}
  \le & \int_0^t \lt(  \lt\|\sa \lt(\frac{v}{x}\rt)''\rt\|_0+ \lt\|\sa v'''\rt\|_{0} \rt)ds
  \le  C tP\lt(\sup_{[0,t]} \sqrt{E}\rt),
    \end{split}\ee
\be\label{eq}\begin{split}
   \lt\| \sa \mathcal{Q} (t)\rt\|_{0}
  \le &  C\lt\|\sa\lt(\frac{v}{x}\rt)''(t)\rt\|_0+ \lt\|\sa v'''(t)\rt\|_{0}
  + \lt\|\sa \mathcal{R}_0 (t) \rt\|_{L^\iy }\lt\|  \mathfrak{L}_0(t)\rt\|_0
  \\
  &+ \lt(  \lt\| \sa \mathcal{R}_1(t)\rt\|_{0} +  \lt\|\sa \mathcal{R}_0 (t) \rt\|_{L^\iy }\lt\|  \mathcal{R}_0 (t)\rt\|_0\rt)\lt\|\mathfrak{J}_0 (t) \rt\|_{L^\iy }\\
  \le & C \sqrt{E(t)}+ C tP\lt(\sup_{[0,t]} \sqrt{E}\rt),
    \end{split}\ee
\be\label{tq}\begin{split}
   \lt\| \za \sa \mathcal{Q}(t) \rt\|_{0}^2
  \le &  C \lt\|\za \sa \lt(\frac{v}{x}\rt)''(t)\rt\|_0^2 + C \lt\|\za \sa v'''(t)\rt\|_{0}^2
  + C \lt\|\sa \mathcal{R}_0 (t) \rt\|_{L^\iy }^2\lt\|  \mathfrak{L}_0(t)\rt\|_0^2
  \\
  &+ C \lt(  \lt\| \sa \mathcal{R}_1(t)\rt\|_{0} +  \lt\|\sa \mathcal{R}_0(t)  \rt\|_{L^\iy }\lt\|  \mathcal{R}_0 (t)\rt\|_0\rt)^2\lt\|\mathfrak{J}_0 (t) \rt\|_{L^\iy }^2\\
  \le &  M_0+ C t P\lt(\sup_{[0,t]} E\rt).
    \end{split}\ee

Now, we are ready to deal with $\lt\|\za H_0'\rt\|_0$. For $\mathfrak{H}_1$, it follows from \ef{norm}, \ef{r0}, \ef{egn}, $\ef{3egn}_1$ and $\ef{3lk}_1$ that
\be\label{l11}\begin{split}
  \lt\| \za \mathfrak{H}_1'(t) \rt\|_0^2 \le & C  \lt\{ \lt\|     \frac{ x^4 } {r^4 r' }-1 \rt\|_{L^\iy}^2 +
    \lt\| \frac{x^3  } {r^3 r'^2}-1\rt\|_{L^\iy}^2
 +\lt\|\frac{x^2  } {r^2 r'^3}-1\rt\|_{L^\iy}^2\rt\}\\
 &\times \lt(\lt\|  \frac{\pl_t   v}{x} \rt\|_1^2  + \lt\|    \pl_t  v' \rt\|_0^2+\lt\| \za  \pl_t  v'' \rt\|_0^2\rt) \\
 &+ C \lt\|\mathcal{R}_0\rt\|_0^2 \lt\{\lt\| \za  \pl_t  v' \rt\|_{L^\iy}^2    + \lt\| {\pl_t v}/ {x}    \rt\|_{L^\iy}^2\rt\}
 \le   C t P\lt(\sup_{[0,t]} E\rt).
    \end{split}\ee
For $\mathfrak{H}_2$, it follows from \ef{norm}, \ef{egn},  $\ef{3lk}_1$, \ef{3egn}  and \ef{rrr} that
\be\label{l12}\begin{split}
  \lt\| \za \mathfrak{H}_2'(t) \rt\|_0^2 \le &
  C \lt\|\mathcal{R}_0\rt\|_0^2 \lt\{\lt\|    \pl_t  v/x \rt\|_{L^\iy}^2+ \lt\| \za  \pl_t  v' \rt\|_{L^\iy}^2 + \lt\| \za \sa ( \pl_t  v/x)'  \rt\|_{L^\iy}^2 + \lt\| \za \sa \pl_t  v''  \rt\|_{L^\iy}^2  \rt\} \\
 & +C \lt(\lt\|\mathcal{R}_0\rt\|_0^2 \lt\|\sa\mathcal{R}_0\rt\|_{L^\iy}^2
 + \lt\|\sa \mathcal{R}_1\rt\|_0^2 \rt) \lt\{\lt\|   \pl_t  v/x \rt\|_{L^\iy}^2+ \lt\|\za   \pl_t  v' \rt\|_{L^\iy}^2   \rt\} \\
  &+C\lt\{
    \lt\| \frac{x^3  } {r^3 r'^2}-1\rt\|_{L^\iy}^2 +\lt\|\frac{x^2  } {r^2 r'^3}-1\rt\|_{L^\iy}^2
 \rt\}  \lt\{ \lt\| \za  \sa \lt(\frac{\pl_t  v}{x} \rt)''\rt\|_0^2
  + \lt\|  \za  \sa\pl_t  v''' \rt\|_0^2 \rt.\\
 &\lt. +\lt\|     \lt(\frac{\pl_t  v}{x} \rt)'\rt\|_0^2 + \lt\|  \za   \pl_t  v'' \rt\|_0^2   \rt\}
\le    C t P\lt(\sup_{[0,t]} E\rt) ,
    \end{split}\ee
since
$$\lt\| \za \sa ( \pl_t  v/x)'  \rt\|_{L^\iy}\le C\lt\| \za x ( \pl_t  v/x)'  \rt\|_{L^\iy} \le C\lt\| \za \lt(\pl_t v' - \pl_t v /x\rt) \rt\|_{L^\iy},$$
$$\lt\| \za  \sa \lt(\frac{\pl_t  v}{x} \rt)''\rt\|_0 \le C\lt\| \za  x \lt(\frac{\pl_t  v}{x} \rt)''\rt\|_0
\le C \lt\| \za \pl_t v '' - 2\za (\pl_t v /x)'\rt\|_0.$$
For the term involving $I_{11}$ and $I_{12}$,  we have from \ef{14}, \ef{314} and \ef{114} that
\bee\label{ehk}\begin{split}
  & \lt\| \za  \lt\{\frac{1}{\sa} \left[  \sa^2  (I_{11}+I_{12}) \right]' \rt\}' \rt\|_0\\
     \le&  C \lt\|\za (I_{11}+I_{12}) \rt\|_0 + C \lt\|\za (I_{11}+I_{12})' \rt\|_0 +\lt\|\za \sa (I_{11}+I_{12})'' \rt\|_0 \\
   \le & C  \lt\|\za \sa \mathcal{Q}  \lt(\lt|v/x\rt| + \lt| v'\rt|\rt) \rt\|_0 + C\lt\|\za \mathfrak{L}_0  \lt(|v/x|+\lt|v'\rt| + |\sa(v/x)'|+\lt|\sa v''\rt|\rt)  \rt\|_0 \\
    & + C\lt\| \za\mathfrak{J}_0  \lt(|v/x|+|v'|+|(v/x)'|+\lt|v''\rt| + |\sa(v/x)''| + \lt|\sa v'''\rt|\rt) \rt\|_0 .
    \end{split}\eee
Note that
\bee\label{}\begin{split}
 &\lt\|\za \sa \mathcal{Q} (t)\lt(\lt|v/x\rt| + \lt| v'\rt|\rt)(t)\rt\|_0^2 \\ \le & 2 \lt\|\za \sa \mathcal{Q}  (t)\rt\|_0^2  \lt\|v(t)/x\rt\|_{L^\iy}^2  + 2\lt(\lt\| \za  \sa\mathcal{Q}(t) \rt\|_0  \lt\| v'(0)\rt\|_{L^\iy}
  +  \lt\|  \sa \mathcal{Q} \rt\|_{0}  \int_0^t   \lt\| \za \pl_t v' \rt\|_{L^\iy} ds\rt)^2 \\
\le &  M_0+ C t P\lt(\sup_{[0,t]} E\rt),
        \end{split}\eee
where we have used \ef{egn}, \ef{tegn}, $\ef{3egn}_1$, \ef{eq} and \ef{tq};
\bee\label{}\begin{split}
 \lt\|\za \mathfrak{L}_0 (t)\lt(|v/x|+\lt|v'\rt| + |\sa(v/x)'|+\lt|\sa v''\rt|\rt)(t)\rt\|_0^2
 \le     M_0+ C t P\lt(\sup_{[0,t]} E\rt),
        \end{split}\eee
due to $\ef{3lk}_{2,3}$; and
\bee\label{}\begin{split}
\lt\| \za\mathfrak{J}_0(t) \lt(|v/x|+|v'|+|(v/x)'|+\lt|v''\rt| + |\sa(v/x)''| + \lt|\sa v'''\rt|\rt)(t)\rt\|_0^2 \le
M_0+ C t P\lt(\sup_{[0,t]} E\rt),
        \end{split}\eee
since
$$|\sa(v/x)''|\le C |x(v/x)''|=C |v''-2(v/x)'|$$
and
\bee\label{}\begin{split}
& \lt\|\za \mathfrak{J}_0(t) \lt(|v''|+\lt|\sa v'''\rt| \rt)(t)\rt\|_0^2
 \\ \le &  \lt[\lt\| \za  \mathfrak{J}_0 (t) \rt\|_{L^\iy}  \lt(\lt\| v''(0)\rt\|_{0} +\lt\|\sa v'''(0)\rt\|_{0}\rt)
   +  \lt\| \mathfrak{J}_0(t)  \rt\|_{L^\iy}  \int_0^t  \lt( \lt\| \za \pl_t v'' \rt\|_{0} +\lt\| \za \sa \pl_t v''' \rt\|_{0}\rt)ds \rt]^2\\
\le  &  M_0+ C t P\lt(\sup_{[0,t]} E\rt).
        \end{split}\eee
Then, we have arrived at
\be\label{l13}\begin{split}
  & \lt\| \za  \lt\{\frac{1}{\sa} \left[  \sa^2  (I_{11}+I_{12}) \right]' \rt\}' (t)\rt\|_0^2  \le M_0+ C t P\lt(\sup_{[0,t]} E\rt).
    \end{split}\ee
In a similar but easier way as for \ef{l13}, one can show
 \be\label{l14}\begin{split}
    \lt\|\za \lt[(\sa/x)(3I_{21}+I_{22})\rt]'(t)\rt\|_0^2    \le  M_0+ C t P\lt(\sup_{[0,t]} E\rt).
    \end{split}\ee
Finally, the last term in $\za \mathcal{G}'$ can be bounded as
\be\label{l15}\begin{split}
  \lt\| \za  \lt[ \phi   { x^2 } \pl_t^{2}\lt(\frac{x^2}{r^2}\rt)\rt]'(t)
 \rt\|_0^2
  \le & C   \lt\| \za  x  \pl_t^2 \lt(\frac{x^2}{r^2}\rt)
 \rt\|_0^2+ C\lt\| \za x^2 \pl_t^2 \lt(\frac{x^2}{r^2}\rt)'
 \rt\|_0^2\\
 \le &C   \lt\| \za  x  \mathfrak{J}_1
 \rt\|_0^2+ C\lt\| \za x^2  \mathfrak{L}_1
 \rt\|_0^2 \\
 \le &  C\lt\| \za     \mathfrak{J}_1
 \rt\|_0^2+ C\lt\| \za  (x^2/\sa) \sa  \mathfrak{L}_1
 \rt\|_0^2
 \le  M_0+ C t P\lt(\sup_{[0,t]} E\rt),
    \end{split}\ee
due to $\ef{3jk}_2$, $\ef{3lk}_3$ and the lower bound of $\rho_0$ in the interior region.

It follows from \ef{ik3'''}, \ef{eie}, \ef{hk3}, \ef{l11}, \ef{l12}, \ef{l13}-\ef{l15} that
 \bee\label{}\begin{split}
\sup_{[0,t]} \lt\|\za H_0'\rt\|_0^2 \le  M_0+ C t P\lt(\sup_{[0,t]} E\rt) +  C \sup_{[0,t]}\lt\|\za \pl_t^3 v\rt\|_1^2  \le  M_0+ C t P\lt(\sup_{[0,t]} E\rt).
\end{split}\eee
In view of \ef{vt'''} and \ef{tnorm}, we can then obtain
\bee\label{}\begin{split}
     &\sup_{[0,t]}\lt(\lt\| \za \sa    {\pl_t  v}  ''' \rt\|_0^2+  \lt\| \za   {\pl_t  v}   ''\rt\|_0^2
    +   \lt\| \za    \lt(  {\pl_t  v}/{x}   \rt)' \rt\|_0^2\rt) \\
\le & C \sup_{[0,t]}\lt[\lt\|\za  H_0'   \rt\|_0^2 + \lt\|\za {\pl_t v}/{x}\rt\|_0^2
+\lt\|\za \pl_t v'\rt\|_0^2
\rt] + C(\da) \sup_{[0,t]}\lt[\lt\|\sa \pl_t v''\rt\|_0^2 + \lt\| \pl_t v'\rt\|_0^2 +\lt\| \pl_t v\rt\|_0^2   \rt]\\
\le & M_0+ C t P\lt(\sup_{[0,t]} E\rt),
 \end{split}\eee
where we used the fact that $\sa(x)\ge m_0 \da $ on $[\da,2\da]$. This, together with \ef{tnorm} and $\ef{3egn}_2$ produces that
\be\label{hk1}\begin{split}
     \sup_{[0,t]}\lt(\lt\| \za \sa    {\pl_t  v}   \rt\|_3^2+  \lt\| \za   {\pl_t  v}   \rt\|_2^2
    +   \lt\| \za    \lt( \frac{\pl_t  v}{x}   \rt)  \rt\|_1^2\rt)
\le & M_0+ C t P\lt(\sup_{[0,t]} E\rt).
 \end{split}\ee
Then we can derive from \ef{hk1} and \ef{norm} that
\be\label{iry0}\begin{split}
   \sup_{[0,t]}\lt(  \lt\| \za \sa      v    \rt\|_3^2+  \lt\| \za   v  \rt\|_2^2
    +   \lt\| \za    \lt( \frac{  v}{x}   \rt)  \rt\|_1^2 \rt)
\le & M_0+ C t P\lt(\sup_{[0,t]} E\rt).
 \end{split}\ee
\subsection{Elliptic estimates -- boundary estimates}
 For the boundary estimates, we introduce a cut-off function $\chi(x)$ satisfying
\be\label{chi}
 \chi=1 \  \ {\rm on} \ \ [\da,1], \ \ \chi=0 \ \ {\rm on} \ \  [0, \da/2], \ \  |\chi'|\le s_0/\da,
\ee
for some constant $s_0$, where $\da$ is given by \ef{da}.  Let
\be\label{h-B}\begin{split}
B =&    \sa  \pl_t^{k } v''+2\sa'  \pl_t^{k } v'=H_0+2\sa'\pl^k_t v/x.
\end{split}\ee
Since for any function $h=h(x,t)$ and integer $i\ge 2$, it holds that
\be\label{hzz}\begin{split}
  &\lt\| \chi \sa   h'\rt\|_0^2 +  \lt\| \chi \sa' h \rt\|_0^2
 \le     \lt\| \chi \lt(\sa   h' + i \sa' h \rt)  \rt\|_0^2 + C \lt\|  \sa^{1/2}  h  \rt\|_0^2 ,
 \\
 & \lt\| \chi \sa^{3/2}   h'\rt\|_0^2 + \lt\| \chi \sa^{1/2} \sa' h \rt\|_0^2 \le
4 \lt\| \chi \sa^{1/2} \lt(\sa   h' + i \sa' h \rt)  \rt\|_0^2 + C \lt\|  \sa  h  \rt\|_0^2.
 \end{split}\ee
We can see that
\be\label{bht}\begin{split}
&\lt\| \chi \sa^{3/2} \pl_t^3 v''\rt\|_0^2 + \lt\| \chi \sa^{1/2} \sa' \pl_t^3 v' \rt\|_0^2 \le
4 \lt\| \chi \sa^{1/2} B \rt\|_0^2 + C \lt\|  \sa  \pl_t^3 v'  \rt\|_0^2 , \ \ k=3;
\\
  &\lt\| \chi \sa   \pl_t^2 v'' \rt\|_0^2 +  \lt\| \chi \sa' \pl_t^2 v' \rt\|_0^2
 \le     \lt\| \chi B  \rt\|_0^2 + C \lt\|  \sa^{1/2}  \pl_t^2 v' \rt\|_0^2 , \ \  k=2;
\\
&\lt\| \chi \sa^{3/2} \pl_t  v'''\rt\|_0^2 + \lt\| \chi \sa^{1/2} \sa' \pl_t  v'' \rt\|_0^2 \le
4 \lt\| \chi \sa^{1/2} \lt( B' - 2\sa'' \pl_t v'\rt) \rt\|_0^2 + C \lt\|  \sa  \pl_t  v''  \rt\|_0^2 , \ \ k=1;
\\
  &\lt\| \chi \sa    v''' \rt\|_0^2 +  \lt\| \chi \sa'   v'' \rt\|_0^2
 \le     \lt\| \chi \lt( B' - 2\sa''   v'\rt)  \rt\|_0^2 + C \lt\|  \sa^{1/2}  v'' \rt\|_0^2 , \ \  k=0.
 \end{split}\ee
Thus, we need to deal with $\|\sa^{1/2}\chi B \|_0$ when $k=3$, $\|\chi B\|_0$ for $k=2$, $\|\sa^{1/2}\chi B' \|_0$ when $k=1$ and $\|\chi B'\|_0$ for $k=0$. The proof of \ef{hzz} is left to the appendix.

\subsubsection{Boundary estimates for $\pl_t^2 v$}
To estimate $\|\chi B\|_0$ with $k=2$, we consider equation \ef{eie} with $k=2$. To this end,
we will first list some useful facts. Similar to  \ef{tnorm}, one can obtain also
\be\label{ho2}\begin{split}
 \lt\| \lt( \sa^{1/2}\pl_t^2 v', \ \sa^{3/2} \pl_t^2 v'', \ \sa^{1/2}    v'',    \ \sa^{3/2}  v'''
\rt)(\cdot,t)\rt\|_{0}^2
 \le  M_0 + CtP\lt(\sup_{[0,t]}E\rt).
 \end{split}\ee
Setting $\|\cdot\|=\|\cdot(t)\|$,  we can summarize from $\ef{jk}$, \ef{tjk}, $\ef{3lk}$, \ef{rrr} and \ef{r0} that
\be\label{lb1}\begin{split}
&\lt\| x/r-1 \rt\|_{L^\iy}  +
    \lt\| 1/r'-1\rt\|_{L^\iy} +\lt\| \mathcal{R}_0 \rt\|_{0}  +
 \lt\|  \sa \mathcal{R}_0\rt\|_{L^\iy } + \lt\|\sa \mathcal{R}_1\rt\|_0 \le C t P\lt( \sup_{[0,t]} \sqrt{E}\rt) ,
  \\
&\lt\|    \mathfrak{J}_0\rt\|_{L^\iy}^2+ \lt\|    \mathfrak{J}_1\rt\|_{L^4}^2+ \lt\|  \sa \mathfrak{L}_0\rt\|_{L^\iy}^2+ \lt\|   \sa \mathfrak{L}_1\rt\|_{L^4}^2 \le  C P\lt(E(t)\rt) + C t P\lt(\sup_{[0,t]} E\rt),
 \\
&\lt\|   \mathfrak{J}_0\rt\|_{L^4}^2+ \lt\|  \mathfrak{J}_1\rt\|_{0}^2 +\lt\| \sa \mathfrak{L}_0\rt\|_{L^4}^2+ \lt\|  \sa \mathfrak{L}_1\rt\|_{0}^2  \le  M_0 + CtP\lt(\sup_{[0,t]}E\rt).
\end{split}\ee

Next, we will deal with the terms on the right-hand side of \ef{eie}. It follows from \ef{tnorm} and \ef{et4} that
 \be\label{bt21}\begin{split}
&\lt\|\chi \lt(\frac{1}{2}x   \pl_t^{4} v  - 4  \lt(\frac{\sa}{x}\rt)' \pl_t^2 v   \rt)(t)\rt\|_0^2\\
 \le & C \lt\|  \pl_t^{4} v (t) \rt\|_0^2 + C(\da)\lt\|
  { \pl_t^{2} v }  (t)\rt\|_0^2 \le  M_0+ C t P\lt(\sup_{[0,t]} E\rt).
\end{split}\ee
For $\mathfrak{H}_1$ and $\mathfrak{H}_2$, by virtue of $\ef{lb1}_1$, \ef{norm}, and Hardy's inequality, one has
  \be\label{bt22}\begin{split}
  \lt\| \chi \mathfrak{ H}_1 (t)\rt\|_0^2 \le & C  \lt\{ \lt\|     \frac{ x^4 } {r^4 r' }-1 \rt\|_{L^\iy}^2 +
    \lt\| \frac{x^3  } {r^3 r'^2}-1\rt\|_{L^\iy}^2 +\lt\|\frac{x^2  } {r^2 r'^3}-1\rt\|_{L^\iy}^2
 \rt\}  \lt\|   {\pl_t^{2} v}  \rt\|_1^2   \\
 \le &  C t P\lt(\sup_{[0,t]} E\rt) ,
    \end{split}\ee
\be\label{bt23}\begin{split}
  \lt\| \chi \mathfrak{H}_2 (t)\rt\|_0^2 \le &C\lt\| \sa \mathcal{R}_0 \rt\|_{L^\iy}^2  \lt\|   {\pl_t^{2} v } \rt\|_{1}^2
   +\lt\{
    \lt\| \frac{x^3  } {r^3 r'^2}-1\rt\|_{L^\iy}^2 +\lt\|\frac{x^2  } {r^2 r'^3}-1\rt\|_{L^\iy}^2
 \rt\}   \\
 &\times \lt( \lt\|   \sa  {\pl_t^{2 } v}  \rt\|_2^2 +\lt\|    {\pl_t^{2 } v}  \rt\|_1^2  \rt)  \\
 \le &   C t P\lt(\sup_{[0,t]} E\rt).
    \end{split}\ee
For the term involving $I_{11}$ and $I_{12}$, we derive from \ef{14} and \ef{314} that
\be\label{bt24}\begin{split}
   \lt\| \chi \frac{1}{\sa} \left[ \sa^2  (I_{11}+I_{12}) \right]'(t) \rt\|_0^2 \
  \le   C \sum_{\aa=0,1}    \lt\{   \lt\|\chi   \mathfrak{J}_{1-\aa}(t) \lt( \lt| {\pl_t^\aa v}    \rt|+|\pl_t^\aa v'|  +  |\sa \pl_t^\aa v''|\rt)(t)\rt\|_0^2 \rt.\\
  \lt. +      \lt\|\chi  \sa \mathfrak{L}_{1-\aa} (t)\lt( |\pl_t^\aa v |  +  | \pl_t^\aa v'|
  \rt)(t) \rt\|_0^2 \rt\}
  \le  M_0 + CtP\lt(\sup_{[0,t]}E\rt).
 \end{split}\ee
Indeed, it follows from \ef{egn}, \ef{tegn} and \ef{lb1} that
\bee\label{}\begin{split}
   & \sum_{\aa=0,1}  \lt(   \lt\|  \mathfrak{J}_{1-\aa}(t)   {\pl_t^\aa v} (t) \rt\|_0^2 +  \lt\|  \sa  \mathfrak{L}_{1-\aa} (t) \pl_t^\aa v (t) \rt\|_0^2 \rt)\\
   \le & C \lt(\lt\|   \mathfrak{J}_{1 }(t) \rt\|_0^2 + \lt\|  \sa \mathfrak{L}_{1 }(t) \rt\|_0^2 \rt)  \lt\|  v (t)\rt\|_{L^\iy}^2 +  C \lt(\lt\|   \mathfrak{J}_{0 } (t)\rt\|_0^2 + \lt\|  \sa \mathfrak{L}_{0}(t) \rt\|_0^2 \rt)  \lt\| \pl_t v (t)\rt\|_{L^\iy}^2 \\
   \le &  M_0 + CtP\lt(\sup_{[0,t]}E\rt),
 \end{split}\eee
\bee\label{}\begin{split}
 &\sum_{\aa=0,1}   \lt\|  \mathfrak{J}_{1-\aa}(t)\lt( \lt| \pl_t^\aa v'   \rt|+\lt| \sa  {\pl_t^\aa v''}      \rt| \rt)(t) \rt\|_0
   \\
\le & \lt\|   \mathfrak{J}_0(t) \rt\|_{L^4} \lt(\lt\| {\pl_t  v'(0) }  \rt\|_ {L^4}  +\lt\| \sa \pl_t  v''(0)  \rt\|_{L^4} \rt) + \int_0^t   \lt( \lt\|   {\pl_t^2 v'}     \rt\|_{0}+ \lt\|\sa \pl_t^2   v''  \rt\|_{0}\rt)  ds \lt\| \mathfrak{J}_0 (t)\rt\|_{L^\iy}\\
&+\lt\|  \mathfrak{J}_1 (t)\rt\|_0 \lt(\lt\| {  v'(0) }  \rt\|_ {L^\iy}  +\lt\|   \sa v''(0)  \rt\|_{L^\iy}\rt) + \int_0^t   \lt( \lt\|  {\pl_t  v'}   \rt\|_{L^4}+ \lt\|\sa \pl_t  v'   \rt\|_{L^4}\rt)  ds \lt\| \mathfrak{J}_1 (t) \rt\|_{L^4}\\
\le & M_0 + CtP\lt(\sup_{[0,t]}E\rt)
\end{split}\eee
and
\bee\label{}\begin{split}
  \sum_{\aa=0,1}     \lt\|  \sa \mathfrak{L}_{1-\aa}  (t)    \pl_t^\aa v'(t)
  \rt\|_0^2 \le
         M_0 + CtP\lt(\sup_{[0,t]}E\rt),
 \end{split}\eee
so \ef{bt24} follows. Similarly, one can also obtain
 \be\label{bt25}\begin{split}
    \lt\|\chi  (\sa/x )  (3I_{21}+I_{22})(t)\rt\|_0 ^2   \le  M_0+ C t P\lt(\sup_{[0,t]} E\rt).
    \end{split}\ee
Finally, one has
\be\label{bt26}\begin{split}
   \lt\| \chi    \phi   { x^2 } \pl_t^{3}\lt(\frac{x^2}{r^2}\rt)(t)
 \rt\|_0^2
  \le  C \lt\|         \mathcal{I}_2(t)
 \rt\|_0^2
 \le   M_0+ C t P\lt(\sup_{[0,t]} E\rt).
    \end{split}\ee
Here \ef{ik} and \ef{pjk} have been used. Applying \ef{ik3''} with $k=2$ and $\beta=\chi$, with the help of \ef{eie}, \ef{bt21}-\ef{bt26} , we obtain
\bee\label{}\begin{split}
\sup_{[0,t]}\lt\|\chi H_0 \rt\|_0^2 \le M_0+ C t P\lt(\sup_{[0,t]} E\rt).
\end{split}\eee
In view of \ef{h-B} and \ef{tnorm}, one can thus get
 \bee\label{}\begin{split}
\sup_{[0,t]}\lt\|\chi B \rt\|_0^2 \le 2\sup_{[0,t]}\lt\|\chi H_0 \rt\|_0^2
  + C(\da)\sup_{[0,t]}\lt\|\pl_t^2 v\rt\|_0^2   \le M_0+ C t P\lt(\sup_{[0,t]} E\rt).
\end{split}\eee
It follows from this, $\ef{bht}_2$ and \ef{ho2} that
\bee\label{}\begin{split}
\sup_{[0,t]}\lt( \lt\| \chi \sa   \pl_t^2 v'' \rt\|_0^2 +  \lt\| \chi \sa' \pl_t^2 v' \rt\|_0^2 \rt)
 \le  \sup_{[0,t]}\lt(   \lt\| \chi B  \rt\|_0^2 + C \lt\|  \sa^{1/2}  \pl_t^2 v' \rt\|_0^2 \rt)
 \le  M_0+ C t P\lt(\sup_{[0,t]} E\rt).
 \end{split}\eee
This, together with \ef{tnorm}, yields that
\be\label{bde2}\begin{split}
 \sup_{[0,t]} \lt(\lt\| \chi \sa   \pl_t^2 v  \rt\|_2^2 +  \lt\| \chi   \pl_t^2 v  \rt\|_1^2 \rt)
 \le  M_0+ C t P\lt(\sup_{[0,t]} E\rt),
 \end{split}\ee
due to the estimate:
\bee\label{}\begin{split}
\sup_{[0,t]} \lt\| \chi     \pl_t^2 v '  \rt\|_0^2 \le \sup_{[0,t]}\lt( C   \lt\| \chi  \sa \pl_t^2 v'   \rt\|_0^2 + C \lt\| \chi  \sa' \pl_t^2 v'  \rt\|_0^2\rt)
 \le  M_0+ C t P\lt(\sup_{[0,t]} E\rt).
 \end{split}\eee

\subsubsection{Boundary estimates for $ v$}
Consider now \ef{eie} with $k=0$. Our goal is to bound $\lt\|\chi B'\rt\|_0$. It follows from \ef{tnorm} and \ef{bde2} that
 \be\label{bt01}\begin{split}
\lt\|\chi \lt( \frac{1}{2}x   \pl_t^{2} v  - 4  \lt(\frac{\sa}{x}\rt)'   v \rt)'(t)\rt\|_0^2 \le C \lt\| \chi \pl_t^{2} v  (t)\rt\|_1^2 + C \lt\|
  {   v }  (t)\rt\|_1^2 \le  M_0+ C t P\lt(\sup_{[0,t]} E\rt).
\end{split}\ee
For $\mathfrak{H}_1$ and $\mathfrak{H}_2$, it follows from \ef{norm}, \ef{egn} and $\ef{lb1}_1$    that
\be\label{bt02}\begin{split}
 \lt\| \chi \mathfrak{H}_1' (t)\rt\|_0^2 \le & C  \lt\{ \lt\|     \frac{ x^4 } {r^4 r' }-1 \rt\|_{L^\iy}^2 +
    \lt\| \frac{x^3  } {r^3 r'^2}-1\rt\|_{L^\iy}^2
 +\lt\|\frac{x^2  } {r^2 r'^3}-1\rt\|_{L^\iy}^2\rt\}  \lt\|   {   v}  \rt\|_2^2  \\
 &+ C \lt\|\mathcal{R}_0\rt\|_0^2 \lt\{\lt\|    v' \rt\|_{L^\iy}^2    + \lt\| v  \rt\|_{L^\iy}^2\rt\}
 \le  C t P\lt(\sup_{[0,t]} E\rt)
    \end{split}\ee
and
\be\label{bt03}\begin{split}
   \lt\| \chi \mathfrak{H}_2' (t)\rt\|_0^2 \le &
  C \lt\|\mathcal{R}_0\rt\|_0^2 \lt\{\lt\|    v \rt\|_{L^\iy}^2+ \lt\| v' \rt\|_{L^\iy}^2 +  \lt\|   \sa v''  \rt\|_{L^\iy}^2  \rt\} \\
 & +C \lt(\lt\|\mathcal{R}_0\rt\|_0^2 \lt\|\sa\mathcal{R}_0\rt\|_{L^\iy}^2
 + \lt\|\sa \mathcal{R}_1\rt\|_0^2 \rt) \lt\{\lt\|   v  \rt\|_{L^\iy}^2+ \lt\|      v' \rt\|_{L^\iy}^2   \rt\} \\
  &+C\lt\{
    \lt\| \frac{x^3  } {r^3 r'^2}-1\rt\|_{L^\iy}^2 +\lt\|\frac{x^2  } {r^2 r'^3}-1\rt\|_{L^\iy}^2
 \rt\}  \lt\{    \lt\|     \sa   v''' \rt\|_0^2  + \lt\|     v  \rt\|_2^2   \rt\} \\
\le  &  C t P\lt(\sup_{[0,t]} E\rt).
    \end{split}\ee
 Using \ef{tnorm}, one has
\be\label{bt04}\begin{split}
   \lt\| \chi \lt[   \phi   { x^2 } \pl_t \lt(\frac{x^2}{r^2}\rt)\rt]'(t)
 \rt\|_0^2 \le & C(\da) \lt\|       v(t)\rt\|_1^2
 \le   M_0+ C t P\lt(\sup_{[0,t]} E\rt).
    \end{split}\ee
It yields from \ef{ik3'''}, \ef{eie}, \ef{bt01}-\ef{bt04} that
\bee\label{hh}\begin{split}
\sup_{[0,t]}\lt\|\chi H'_0 \rt\|_0^2 \le M_0+ C t P\lt(\sup_{[0,t]} E\rt).
\end{split}\eee
In view of \ef{h-B} and \ef{tnorm}, one gets
 \bee\label{}\begin{split}
\sup_{[0,t]}\lt\|\chi B' \rt\|_0^2 \le 2\sup_{[0,t]}\lt\|\chi H_0 ' \rt\|_0^2
  + C(\da)\sup_{[0,t]}\lt\| v\rt\|_1^2   \le M_0+ C t P\lt(\sup_{[0,t]} E\rt).
\end{split}\eee
We can then obtain, using $\ef{bht}_4$, \ef{tnorm} and \ef{ho2}, that
\bee\label{}\begin{split}
\sup_{[0,t]}\lt( \lt\| \chi \sa    v''' \rt\|_0^2 +  \lt\| \chi \sa'   v'' \rt\|_0^2 \rt)
 \le    \sup_{[0,t]}\lt(  \lt\| \chi \lt( B' - 2\sa''   v'\rt)  \rt\|_0^2 + C \lt\|  \sa^{1/2}  v'' \rt\|_0^2 \rt)\\
 \le   C \sup_{[0,t]} \lt(\lt\| \chi   B' \rt\|_0^2  +   \lt\|  v'   \rt\|_0^2 +   \lt\|  \sa^{1/2}  v'' \rt\|_0^2 \rt)\le M_0+ C t P\lt(\sup_{[0,t]} E\rt).
 \end{split}\eee
This, together with \ef{tnorm}, yields
\be\label{bde0}\begin{split}
 \sup_{[0,t]}\lt(\lt\| \chi \sa     v  \rt\|_3^2 +  \lt\| \chi     v  \rt\|_2^2 \rt)
 \le  M_0+ C t P\lt(\sup_{[0,t]} E\rt),
 \end{split}\ee
since
\bee\label{}\begin{split}
\sup_{[0,t]} \lt\| \chi     v ''   \rt\|_0^2 \le     \sup_{[0,t]} \lt(C\lt\| \chi  \sa v''   \rt\|_0^2 + C \lt\| \chi  \sa'   v''   \rt\|_0^2\rt)
 \le  M_0+ C t P\lt(\sup_{[0,t]} E\rt).
 \end{split}\eee

\subsubsection{Boundary estimates for $\pl_t^3 v$}
Consider equation \ef{eie} with $k=3$.
As before, we list here some estimates which will be used later. First, it follows from \ef{iry2}, \ef{bde2}, \ef{iry0} and \ef{bde0} that
\be\label{f0}\begin{split}
\sup_{[0,t]}\lt(\lt\|   \sa     v  \rt\|_3^2 +  \lt\|    v  \rt\|_2^2 +\lt\| \sa \pl_t^2        v  \rt\|_2^2 +  \lt\| \pl_t^2   v  \rt\|_1^2 \rt)
\le  M_0+ C t P\lt(\sup_{[0,t]} E\rt).
\end{split}\ee
Moreover, we have the following estimates for $\pl_t v$ and $\pl_t^3 v$:
\be\label{ho3}\begin{split}
 & \lt\| \lt( \sa^{1/2}  \pl_t v',    \ \sa^{1/2} \pl_t^3 v
\rt)(\cdot,t)\rt\|_{L^\iy}^2 + \lt\| \lt( \sa^{3/2}  \pl_t v'',    \ \sa^{3/2} \pl_t^3 v'
\rt)(\cdot,t)\rt\|_{L^\iy}^2  \le C{E(t)};
\end{split}\ee
and those for $\mathfrak{J}$ and $\mathfrak{L}$:
\be\label{lb0}\begin{split}
  &\lt\| \sa^{1/2} \lt( \mathfrak{J}_1 ,  \sa \mathfrak{L}_1\rt)(\cdot,t)\rt\|_{L^\iy}^2   \le CP\lt(E(t)  \rt),\\
 & \lt\|    \lt(\mathfrak{J}_0,  \sa \mathfrak{L}_0\rt)(\cdot,t)\rt\|_{L^\iy}^2  + \lt\|    \lt( \mathfrak{J}_1, \sa \mathfrak{L}_1\rt) (\cdot,t)\rt\|_{0}^2 + \lt\|  \lt( \mathfrak{J}_2, \sa \mathfrak{L}_2\rt) (\cdot,t)\rt\|_{0}^2 \le
 M_0+ C t P\lt(\sup_{[0,t]} E\rt).
\end{split}\ee
The proofs of \ef{ho3} and \ef{lb0} will be given in the appendix.

We are now ready to do the estimates. First,  \ef{tnorm} and \ef{et4} imply that
 \be\label{bt31}\begin{split}
 \lt\|
   \chi \sa^{1/2}\lt(\frac{1}{2}x   \pl_t^{5} v - 4  \lt(\frac{\sa}{x}\rt)' \pl_t^3 v\rt)(t)\rt\|_0^2  \le & C(\da) \lt( \lt\| (x\sa)^{1/2} \pl_t^{5} v  (t)\rt\|_0^2 +  \lt\|
   \pl_t^{3} v  (t) \rt\|_0^2 \rt)\\
  \le & M_0+ C t P\lt(\sup_{[0,t]} E\rt).
\end{split}\ee
For $\mathfrak{H}_1$ and $\mathfrak{H}_2$, it follows from \ef{norm} and $\ef{lb1}_1$ that
\be\label{bt32}\begin{split}
  \lt\| \chi \sa^{1/2} \mathfrak{ H}_1 (t) \rt\|_0^2 \le   C  \lt\{ \lt\|     \frac{ x^4 } {r^4 r' }-1 \rt\|_{L^\iy}^2 +
    \lt\| \frac{x^3  } {r^3 r'^2}-1\rt\|_{L^\iy}^2 +\lt\|\frac{x^2  } {r^2 r'^3}-1\rt\|_{L^\iy}^2
 \rt\} \\
  \times \lt\{ \lt\| \sa^{1/2}  {\pl_t^{3 } v}  \rt\|_0^2 + \lt\|  \sa^{1/2} \pl_t^{3 } v' \rt\|_0^2   \rt\}
 \le    C t P\lt(\sup_{[0,t]} E\rt);
    \end{split}\ee
\be\label{bt33}\begin{split}
  \lt\| \chi \sa^{1/2}  \mathfrak{H}_2 (t)\rt\|_0^2   \le &C\lt\| \sa \mathcal{R}_0 \rt\|_{L^\iy}^2 \lt(\lt\|  \sa^{1/2} {\pl_t^{3} v } \rt\|_{0}^2  + \lt\| \sa^{1/2} \pl_t^{3} v'    \rt\|_{0}^2  \rt) \\
  &+C\lt\{
    \lt\| \frac{x^3  } {r^3 r'^2}-1\rt\|_{L^\iy}^2 +\lt\|\frac{x^2  } {r^2 r'^3}-1\rt\|_{L^\iy}^2
 \rt\}    \lt\|    \sa^{3/2} \lt(\pl_t^3 v, \pl_t^3 v', \pl_t^{3 } v'' \rt)\rt\|_0^2    \\
 \le &   C t P\lt(\sup_{[0,t]} E\rt).
    \end{split}\ee
For the term involving $I_{11}$ and $I_{12}$, one can derive from \ef{14} and \ef{314} that
\bee\label{}\begin{split}
& \lt\| \chi \sa^{1/2} \frac{1}{\sa} \left[ \sa^2  (I_{11}+I_{12}) \right]' \rt\|_0^2  \\
  \le & C \sum_{\aa=0}^2   \lt\{   \lt\|\chi \sa^{1/2}   \mathfrak{J}_{2-\aa} \lt( \lt| {\pl_t^\aa v}      \rt|+|\pl_t^\aa v'|  +  |\sa \pl_t^\aa v''|\rt)\rt\|_0^2  +      \lt\|\chi  \sa^{3/2} \mathfrak{L}_{2-\aa} \lt( |\pl_t^\aa v |  +  | \pl_t^\aa v'|
  \rt) \rt\|_0^2 \rt\}.
 \end{split}\eee
Note that
\bee\label{}\begin{split}
 &\sum_{\aa=0 }^2   \lt\|\sa^{1/2}  \mathfrak{J}_{2-\aa}(t)\lt( \lt| \pl_t^\aa v'   \rt|+\lt| \sa  {\pl_t^\aa v''}      \rt| \rt)(t) \rt\|_0
   \\
\le & \lt\|    \mathfrak{J}_0 \rt\|_{L^\iy} \lt(\lt\| {\pl_t^2  v'  }  \rt\|_ {0}  +\lt\| \sa \pl_t^2  v''  \rt\|_{0} \rt) +\lt\|    \mathfrak{J}_1 \rt\|_{0} \lt(\lt\| {\sa^{1/2} \pl_t   v'(0) }  \rt\|_ {L^\iy}  +\lt\| \sa^{3/2} \pl_t  v''(0)  \rt\|_{L^\iy} \rt)
 \\
&+ \int_0^t   \lt( \lt\|   {\pl_t^2 v'}     \rt\|_{0}+ \lt\|\sa  \pl_t^2   v''  \rt\|_{0}\rt)  ds \lt\|  \sa^{1/2} \mathfrak{J}_1 \rt\|_{L^\iy} +\lt\|     \mathfrak{J}_2 \rt\|_{0} \lt(\lt\| {    v'  }  \rt\|_ {L^\iy}  +\lt\| \sa   v''   \rt\|_{L^\iy} \rt)  \\
\le & M_0 + CtP\lt(\sup_{[0,t]}E\rt),
\end{split}\eee
where we have used \ef{norm}, \ef{f0}-\ef{lb0} and $\|\cdot\|_{L^\iy}\le C \|\cdot\|_1$. Similarly, one has
\bee\label{}\begin{split}
  \sum_{\aa=0 }^2     \lt\|   \sa^{3/2} \mathfrak{L}_{2-\aa}   (t)   \pl_t^\aa v'(t)
  \rt\|_0^2 \le & \sum_{\aa=0 }^2    \lt\|  \sa^{1/2} \lt( \sa \mathfrak{L}_{2-\aa}\rt) \pl_t^\aa v'  \rt\|_0^2
   \le     M_0 + CtP\lt(\sup_{[0,t]}E\rt),
 \end{split}\eee
and
\bee\label{}\begin{split}
   & \sum_{\aa=0}^2  \lt(   \lt\|  \sa^{1/2} \mathfrak{J}_{2-\aa} (t)  {\pl_t^\aa v} (t)\rt\|_0^2 +  \lt\|  \sa^{3/2}  \mathfrak{L}_{2-\aa}  (t) \pl_t^\aa v  (t)\rt\|_0^2 \rt)\\
   \le & C \lt(\lt\|     \mathfrak{J}_{2 } \rt\|_0^2 + \lt\| \sa  \mathfrak{L}_{2 } \rt\|_0^2 \rt)  \lt\|  v \rt\|_{L^\iy}^2 +  C \lt(\lt\|   \mathfrak{J}_{1 } \rt\|_0^2 + \lt\|  \sa \mathfrak{L}_{1} \rt\|_0^2 \rt)  \lt\| \pl_t v \rt\|_{L^\iy}^2 \\
    &+ C \lt(\lt\|  \mathfrak{J}_{0 } \rt\|_{L^\iy}^2 + \lt\|  \sa  \mathfrak{L}_{ 0 } \rt\|_{L^\iy}^2 \rt)  \lt\| \pl_t^2 v \rt\|_{0}^2\\
   \le &  M_0 + CtP\lt(\sup_{[0,t]}E\rt),
 \end{split}\eee
where we have used \ef{tnorm}, \ef{tegn}, \ef{f0}, $\ef{lb0}_2$ and $\|\cdot\|_{L^\iy}\le C\|\cdot\|_1$.
Hence, it holds that
\be\label{bt34}\begin{split}
  \lt\| \chi \sa^{1/2} \frac{1}{\sa} \left[ \sa^2  (I_{11}+I_{12}) \right]' (t) \rt\|_0^2  \le M_0 + CtP\lt(\sup_{[0,t]}E\rt),
 \end{split}\ee
Similarly, one can also obtain easily that
 \be\label{bt35}\begin{split}
    \lt\|\chi \sa^{1/2}(\sa/x) (3I_{21}+I_{22})(t)\rt\|_0 ^2   \le  M_0+ C t P\lt(\sup_{[0,t]} E\rt).
    \end{split}\ee
Finally, one has
\be\label{bt36}\begin{split}
   \lt\|  \chi \sa^{1/2}  \phi   { x^2  } \pl_t^{4}\lt(\frac{x^2}{r^2}\rt)(t)
 \rt\|_0^2
  \le  C \lt\|     x  \mathcal{I}_3(t)
 \rt\|_0^2
 \le   M_0+ C t P\lt(\sup_{[0,t]} E\rt).
    \end{split}\ee
Here \ef{ik} and \ef{pjk} were used. Now, it follows from \ef{eie}, \ef{bt31}-\ef{bt36}, by applying \ef{ik3''} with $\ba=\chi\sa^{1/2}$, that
\bee\label{}\begin{split}
\sup_{[0,t]}\lt\|\chi \sa^{1/2} H_0 \rt\|_0^2 \le M_0+ C t P\lt(\sup_{[0,t]} E\rt).
\end{split}\eee
Thanks to \ef{h-B} and \ef{tnorm}, one can then get
 \bee\label{}\begin{split}
\sup_{[0,t]}\lt\|\chi \sa^{1/2} B \rt\|_0^2 \le 2\sup_{[0,t]}\lt\|\chi\sa^{1/2} H_0 \rt\|_0^2
  + C(\da)\sup_{[0,t]}\lt\|\pl_t^3 v\rt\|_0^2   \le M_0+ C t P\lt(\sup_{[0,t]} E\rt).
\end{split}\eee
It then follows from $\ef{bht}_1$ and \ef{tnorm} that
\bee\label{}\begin{split}
\sup_{[0,t]}\lt(  \lt\| \chi \sa^{3/2} \pl_t^3 v''\rt\|_0^2 + \lt\| \chi \sa^{1/2} \sa' \pl_t^3 v' \rt\|_0^2 \rt)\le & \sup_{[0,t]}\lt(
4 \lt\| \chi \sa^{1/2} B \rt\|_0^2 + C \lt\|  \sa  \pl_t^3 v'  \rt\|_0^2\rt)\\
\le & M_0+ C t P\lt(\sup_{[0,t]} E\rt).
 \end{split}\eee
This, together with \ef{tnorm} and the Sobolev embedding \ef{weightedsobolev}, yields
\be\label{bde3}\begin{split}
\sup_{[0,t]}\lt( \lt\| \chi \sa^{3/2} \pl_t^3 v''\rt\|_0^2 + \lt\| \chi \sa^{1/2}   \pl_t^3 v' \rt\|_0^2 +  \lt\| \chi     \pl_t^3 v \rt\|_{1/2}^2 \rt)
 \le  M_0+ C t P\lt(\sup_{[0,t]} E\rt),
 \end{split}\ee
because of
\bee\label{}\begin{split}
\sup_{[0,t]} \lt\| \chi   \sa^{1/2}  \pl_t^3 v ' \rt\|_0^2 \le & \sup_{[0,t]}\lt(C   \lt\| \chi  \sa^{3/2} \pl_t^3 v'  \rt\|_0^2 + C \lt\| \chi \sa^{1/2} \sa' \pl_t^3 v ' \rt\|_0^2\rt)\\
 \le & \sup_{[0,t]}\lt(C   \lt\| \chi  \sa  \pl_t^3 v'  \rt\|_0^2 + C \lt\| \chi \sa^{1/2} \sa' \pl_t^3 v ' \rt\|_0^2\rt)
 \le  M_0+ C t P\lt(\sup_{[0,t]} E\rt)
 \end{split}\eee
and
\bee\label{}\begin{split}
 \sup_{[0,t]}\lt\| \chi     \pl_t^3 v  \rt\|_{1/2}^2 \le &\sup_{[0,t]}\lt( C   \lt\| \sa^{1/2} \pl_t^3 v \rt\|_0^2 + C \lt\| \chi \sa^{1/2}  \pl_t^3 v ' \rt\|_0^2
 \rt)
 \le  M_0+ C t P\lt(\sup_{[0,t]} E\rt).
 \end{split}\eee

\subsubsection{Boundary estimates for $ \pl_t v$}
Consider equation \ef{eie} with $k=1$. Our goal is to bound $\lt\|\chi \sa^{1/2} B'\rt\|_0$. It follows from \ef{tnorm} and \ef{bde3} that
 \be\label{bt11}\begin{split}
&\lt\|\chi \sa^{1/2}\lt( \frac{1}{2}x   \pl_t^{3} v   - 4 \lt(\frac{\sa}{x}\rt)' \pl_t v \rt)'(t)\rt\|_0^2 \\
  \le & C  \lt\|   \pl_t^{3} v (t) \rt\|_0^2 + \lt\| \chi \sa^{1/2} \pl_t^{3} v'  (t)\rt\|_0^2+  C(\da)\lt\|
 \pl_t v  (t)\rt\|_1^2
  \le   M_0+ C t P\lt(\sup_{[0,t]} E\rt).
\end{split}\ee
For $\mathfrak{H}_1$ and $\mathfrak{H}_2$, it follows from \ef{norm}, \ef{egn}, \ef{tegn},  $\ef{lb1}_1$  and \ef{ho3} that
\be\label{bt12}\begin{split}
  \lt\| \chi  \sa^{1/2} \mathfrak{H}_1' (t)\rt\|_0^2 \le & C  \lt\{ \lt\|     \frac{ x^4 } {r^4 r' }-1 \rt\|_{L^\iy}^2 +
    \lt\| \frac{x^3  } {r^3 r'^2}-1\rt\|_{L^\iy}^2
 +\lt\|\frac{x^2  } {r^2 r'^3}-1\rt\|_{L^\iy}^2\rt\}\\
 &\times \lt(\lt\| {\pl_t   v}  \rt\|_1^2  + \lt\|    \pl_t  v' \rt\|_0^2+\lt\|  \sa^{1/2}  \pl_t  v'' \rt\|_0^2\rt) \\
 &+ C \lt\|\mathcal{R}_0\rt\|_0^2 \lt\{\lt\| \sa^{1/2} \pl_t  v' \rt\|_{L^\iy}^2    + \lt\| {\pl_t v}    \rt\|_{L^\iy}^2\rt\}
 \le   C t P\lt(\sup_{[0,t]} E\rt)
    \end{split}\ee
and
\be\label{bt13}\begin{split}
 &\lt\|\chi \sa^{1/2} \mathfrak{H}_2' (t)\rt\|_0^2 \\
 \le &
  C \lt\|\mathcal{R}_0\rt\|_0^2 \lt\{\lt\|    \pl_t  v \rt\|_{L^\iy}^2+ \lt\|  \sa^{1/2} \pl_t  v' \rt\|_{L^\iy}^2     + \lt\|  \sa^{3/2} \pl_t  v''  \rt\|_{L^\iy}^2  \rt\} \\
   &+C \lt(\lt\|\mathcal{R}_0\rt\|_0^2 \lt\|\sa\mathcal{R}_0\rt\|_{L^\iy}^2
 + \lt\|\sa \mathcal{R}_1\rt\|_0^2 \rt) \lt\{\lt\|   \pl_t  v \rt\|_{L^\iy}^2 + \lt\| \sa^{1/2}  \pl_t  v' \rt\|_{L^\iy}^2   \rt\}
  \\
  & +C\lt\{
    \lt\| \frac{x^3  } {r^3 r'^2}-1\rt\|_{L^\iy}^2 +\lt\|\frac{x^2  } {r^2 r'^3}-1\rt\|_{L^\iy}^2
 \rt\}  \lt\{ \lt\|  \sa^{1/2}  \pl_t  v'' \rt\|_0^2  + \lt\|   \sa^{3/2}\pl_t  v''' \rt\|_0^2
    +\lt\|      {\pl_t  v} \rt\|_1^2   \rt\} \\
\le  &  C t P\lt(\sup_{[0,t]} E\rt).
    \end{split}\ee
For the term involving $I_{11}$ and $I_{12}$, it follows from \ef{14}, \ef{314} and \ef{114} that
\bee\label{}\begin{split}
 & \lt\| \chi \sa^{1/2}  \lt\{\frac{1}{\sa} \left[  \sa^2  (I_{11}+I_{12}) \right]' \rt\}' \rt\|_0 \\
   \le & C  \lt\| \chi \sa^{3/2} \mathcal{Q} \lt(\lt|v/x\rt| + \lt| v'\rt|\rt)\rt\|_0 + C\lt\| \chi \sa^{1/2} \mathfrak{L}_0 \lt(|v/x|+\lt|v'\rt| + |\sa(v/x)'|+\lt|\sa v''\rt|\rt)\rt\|_0 \\
    & + C\lt\|\chi \sa^{1/2}\mathfrak{J}_0 \lt(|v/x|+|v'|+|(v/x)'|+\lt|v''\rt| + |\sa(v/x)''| + \lt|\sa v'''\rt|\rt)\rt\|_0 \\
    \le & C \lt\| \chi \sa^{3/2} \mathcal{Q} \lt(\lt|v\rt| + \lt| v'\rt|\rt)\rt\|_0 + C\lt\| \chi \sa^{1/2} \mathfrak{L}_0 \lt(|v|+\lt|v'\rt| +\lt|\sa v''\rt|\rt)\rt\|_0 \\
    & + C\lt\|\chi \sa^{1/2}\mathfrak{J}_0 \lt(|v|+|v'| + \lt|v''\rt|   + \lt|\sa v'''\rt|\rt)\rt\|_0 .
    \end{split}\eee
Note that one can derive from \ef{f0} and $\ef{lb1}_{1,2}$ that
  \bee\label{}\begin{split}
 \lt\|\chi   \mathfrak{L}_{0}(t) \rt\|_{0}^2 \le & C\lt(\lt\|  v(t)\rt\|_{2}   + \lt\|      \mathcal{R}_0(t)\rt\|_{0}  \lt\|\mathfrak{J}_{0}(t)\rt\|_{L^\iy}  \rt)^2
\le     M_0 + CtP\lt(\sup_{[0,t]}E\rt),
\end{split}\eee
\bee\label{}\begin{split}
   \lt\| \chi \sa^{} \mathcal{Q} (t)\rt\|_{0}^2
  \le &  C \lt\|v(t)\rt\|_2^2 + C \lt\|\sa v'''(t)\rt\|_{0}^2
  + C \lt\|  \mathcal{R}_0 (t) \rt\|_0^2\lt\|  \sa \mathfrak{L}_0(t)\rt\|_{L^\iy }^2
   + C \lt(  \lt\|  \sa \mathcal{R}_1(t)\rt\|_{0} \rt. \\
    &\lt.+  \lt\| \sa\mathcal{R}_0 (t) \rt\|_{L^\iy}\lt\|  \mathcal{R}_0 (t) \rt\|_0\rt)^2\lt\|\mathfrak{J}_0 (t) \rt\|_{L^\iy }^2
  \le  M_0+ C t P\lt(\sup_{[0,t]} E\rt);
    \end{split}\eee
which implies, due to \ef{f0} and $\ef{lb0}_2$, that
\bee\label{}\begin{split}
 &\lt\|\chi \sa \mathcal{Q}(t) \rt\|_0^2 + \lt\|\chi \mathfrak{L}_0(t)\rt\|_0^2 + \lt\| \chi \mathfrak{J}_0(t)\rt\|_{L^\iy}^2 + \|\sa v(t)\|_3^2+ \|v(t)\|_2^2
\le   M_0+ C t P\lt(\sup_{[0,t]} E\rt).
        \end{split}\eee
So, we  obtain
\be\label{bt14}\begin{split}
 & \lt\| \chi \sa^{1/2}  \lt\{\frac{1}{\sa} \left[  \sa^2  (I_{11}+I_{12}) \right]' \rt\}' (t)\rt\|_0^2 \\
  \le  &  C   \lt\|\chi \sa  \mathcal{Q}\rt\|_0^2 \|v\|_2^2  + C\lt\|\chi \mathfrak{L}_0 \rt\|_0 \lt(\|v\|_2^2 +\|\sa v\|_3^2 \rt)  + C\lt\|\chi \mathfrak{J}_0\rt\|_{L^\iy}^2 \lt(\|v\|_2^2 +\|\sa v\|_3^2 \rt) \\
    \le &  M_0+ C t P\lt(\sup_{[0,t]} E\rt),
    \end{split}\ee
where we have used the fact that $\|\cdot\|_{L^\iy}\le C \|\cdot\|_1$. Similarly, one can show that
 \be\label{bt15}\begin{split}
    \lt\|\chi \sa^{1/2} \lt[ (\sa/x)  (3I_{21}+I_{22})\rt]'(t)\rt\|_0^2    \le  M_0+ C t P\lt(\sup_{[0,t]} E\rt).
    \end{split}\ee
It follows from  \ef{tnorm}, $\ef{3lk}_1$ and  \ef{ho3} that
\be\label{bt16}\begin{split}
   & \lt\| \chi \sa^{1/2} \lt[   \phi   { x^2  } \pl_t^2 \lt(\frac{x^2}{r^2}\rt)\rt]'(t)
 \rt\|_0^2 \\
  \le &  C \lt(\lt\|   \pl_t    v\rt\|_1 + \lt\|  v \rt\|_{L^\iy} \|v \|_0 + \|v\|_{L^\iy} \|v'\|_0
 \rt)^2 + C \|\mathcal{R}_0\|_0^2\lt( \|\pl_t v\|_{L^\iy} + \|v\|_{L^\iy}^2   \rt)^2 \\
 \le & C \lt(\lt\|   \pl_t    v\rt\|_1 +   \|v \|_2^2   \rt)^2 + C \|\mathcal{R}_0\|_0^2\lt( \|\pl_t v\|_1 + \|v\|_1^2   \rt)^2
 \le    M_0+ C t P\lt(\sup_{[0,t]} E\rt).
    \end{split}\ee
It yields from \ef{eie}, \ef{ik3'''} and \ef{bt11}-\ef{bt16} that
 \bee\label{}\begin{split}
\sup_{[0,t]}\lt\|\sa^{1/2}\chi H_0' \rt\|_0^2 \le M_0+ C t P\lt(\sup_{[0,t]} E\rt),
\end{split}\eee
which implies
 \bee\label{}\begin{split}
\sup_{[0,t]}\lt\|\sa^{1/2}\chi B' \rt\|_0^2 \le   \sup_{[0,t]}\lt(2\lt\|\sa^{1/2}\chi H_0' \rt\|_0^2 + C(\da)\lt\| \pl_t v  \rt\|_1^2  \rt) \le  M_0+ C t P\lt(\sup_{[0,t]} E\rt),
\end{split}\eee
due to \ef{h-B} and \ef{tnorm}.
We can then obtain, using $\ef{bht}_3$ and \ef{tnorm}, that
\bee\label{}\begin{split}
&\sup_{[0,t]} \lt(\lt\| \chi \sa^{3/2} \pl_t  v'''\rt\|_0^2 + \lt\| \chi \sa^{1/2} \sa' \pl_t  v'' \rt\|_0^2 \rt)\\
\le&\sup_{[0,t]} \lt(
4 \lt\| \chi \sa^{1/2} \lt( B' - 2\sa'' \pl_t v'\rt) \rt\|_0^2 + C \lt\|  \sa  \pl_t  v''  \rt\|_0^2 \rt)\\
 \le  & C\sup_{[0,t]}  \lt(\lt\| \chi \sa^{1/2}  B' \rt\|_0^2  +   \lt\| \pl_t v'   \rt\|_0^2 +   \lt\|  \sa  \pl_t v'' \rt\|_0^2 \rt)\le M_0+ C t P\lt(\sup_{[0,t]} E\rt).
 \end{split}\eee
This,  together with \ef{tnorm} and the Sobolev embedding \ef{wsv},  yields
\be\label{bde1}\begin{split}
\sup_{[0,t]}\lt( \lt\| \chi \sa^{3/2} \pl_t  v'''\rt\|_0^2 + \lt\| \chi \sa^{1/2}   \pl_t  v'' \rt\|_0^2 + \lt\| \chi   \pl_t  v  \rt\|_{3/2}^2\rt)
 \le  M_0+ C t P\lt(\sup_{[0,t]} E\rt).
 \end{split}\ee

\section{Existence for the case  $\gamma=2$}

Summing over inequalities \ef{et4}, \ef{hk3}, \ef{iry2}, \ef{hk1}, \ef{iry0}, \ef{bde2}, \ef{bde0}, \ef{bde3} and \ef{bde1}, we find that
\bee\label{}\begin{split}
\sup_{[0,t]}E
 \le  M_0+ C t P\lt(\sup_{ [0,t]} E \rt), \ \ t\in [0,T];
 \end{split}\eee
which implies that for small $T$,
\be\label{final}\begin{split}
\sup_{t\in[0,T]}E(t)
 \le  2 M_0.
 \end{split}\ee
With this $\mu$-independent estimate, one can use the standard compactness argument \cite{10} to show the existence of the solutions to the problem \ef{419'} for some time $T$.

\section{Case $1<\ga< 2$}
In this section, we use similar arguments to those used to deal with the case for $\ga=2$ to handle the case for general $\ga$.  {It should be noted that the value of $\gamma$ determines the rate of degeneracy near the vacuum boundary, since $\rho_0$ appears as the coefficient in front of $\partial_t v$ in (3.7) and the physical vacuum condition indicates that $\rho_0(x)\sim (1-x)^{\frac{1}{\gamma-1}}$ as
$x\to 1$. Thus the smaller value of $\gamma$ is, the more degenerate equation (3.7) is near the vacuum boundary.  Although the rate of degeneracy near the origin is the same no matter what $\ga$ is,
we need higher order derivatives in the energy functional to control the $H^2$-norm of $v$ (and thus the $C^1$-norm of $v$) for smaller $\gamma$, since we have to match the norms in the intermediate region.}

{We first define the higher-order energy functional for $1<\ga<2$.} Set
$$\nu:=(2-\ga)/(2\ga-2)>0, \ \  l:= 3 + 2 \lceil 1/2 +\nu \rceil,$$
where $\lceil\cdot\rceil$ is the ceiling function defined for any real number $q\ge 0$ as
$$\lceil q\rceil :=\min\{m: \ \  m\ge q, \  m   \ {\rm is~an~integer}\}.$$
Define
\be\label{ggnorm}\begin{split}
 \widetilde{{E}}(v, t):= & \lt\| \sa (\sa/x)^\nu \pl_t^l v'(\cdot,t) \rt\|_0^2 +    \lt\|   (\sa/x)^{1+\nu}  \pl_t^l v (\cdot,t)\rt\|_0^2 \\
      & + \sum_{j=1}^{ \frac{l+1}{2} } \lt\{ \lt\|  \sa^{3/2+\nu} \pl_t^{l-2j +1 } \pl_x^{ j +1 }  v (\cdot,t)\rt\|_{0}^2
      + \sum_{i=0}^j \lt\|  \sa^{1/2+\nu} \pl_t^{l-2j +1 } \pl_x^{ i }  v (\cdot,t)\rt\|_{0}^2 \rt\}  \\
      & + \sum_{j=1}^{ \frac{l-1}{2} } \lt\{ \lt\|  \sa^{2+\nu} \pl_t^{l-2j  } \pl_x^{ j +2 }  v (\cdot,t)\rt\|_{0}^2
      +\sum_{i=-1}^j \lt\|  \sa^{1+\nu} \pl_t^{l-2j  } \pl_x^{ i +1 }  v (\cdot,t)\rt\|_{0}^2 \rt\}\\ &+\sum_{j=1}^{ \frac{l+1}{2}} \lt\{ \lt\|\za \sa \pl_t^{l-2j +1 } v (\cdot,t)\rt\|_{j+1}^2
   +\lt\|\za  \pl_t^{l-2j +1 } v (\cdot,t)\rt\|_{j }^2
   +\lt\|\za \frac {\pl_t^{l-2j +1 } v}{x} (\cdot,t)\rt\|_{j-1}^2  \rt\}\\
 &+ \sum_{j=1}^{\frac{l-1}{2}} \lt\{ \lt\|\za \sa \pl_t^{l-2j  } v (\cdot,t)\rt\|_{j+2}^2
   +\lt\|\za  \pl_t^{l-2j  } v (\cdot,t)\rt\|_{j +1 }^2
   +\lt\|\za \frac{ \pl_t^{l-2j  } v }{x}  (\cdot,t)\rt\|_{j}^2 \rt\}
   ,
\end{split}\ee
where, as before,
\bee
 \za=1 \  \ {\rm on} \ \ [0,\da_\ga], \ \ \za=0 \ \ {\rm on} \ \  [2\da_\ga, 1], \ \  |\za'|\le s_0/\da_\ga.
\eee
Here $\da_\ga$ is a given constant depending on $\rho_0$ and $\ga$ which will be determined in \ef{daga} later.
It follows from the Hardy type embedding for the weighted Sobolev spaces \ef{wsv} that
\be\label{}\begin{split}
 \|v\|_2^2 \le  \|v \|_{\frac{l+1}{2}-\lt(\frac{1}{2}+\nu\rt)}^2  \le  C  \sum_{i=0}^ {\frac{l+1}{2}}  \lt\|  \sa^{1/2+\nu}   \pl_x^{ i }  v \rt\|_{0}^2 \le  C\widetilde{{E}},
\end{split}\ee
which indicates that the high-order energy functional $\widetilde{{E}}$ is  suitable  for the study of the physical vacuum problem \ef{419'} when $\ga\in(1,2)$. {In fact, the norm chosen in \ef{ggnorm}
is in the same spirit of but slightly different from that in \ef{norm} for $\ga=2$. Since the energy estimate gives
the bound of
$$\lt\| \sqrt{x \sa } (\sa/x)^\nu  \pl_t^{l+1} v  (t)\rt\|_0  +    \lt\| \sa (\sa/x)^\nu \pl_t^l v' (t)\rt\|_0 +    \lt\|  (\sa/x)^{1+\nu}  \pl_t^l v (t)\rt\|_0 ,$$
from which we can derive the bound of $\|\pl_t^l v\|_0$ for $\ga=2$. But for $\ga\in(1,2)$, we cannot improve the spatial regularity as that for $\ga=2$ due to $\nu>0$ (or equivalently, the higher degeneracy of the equation). So, the norm chosen for $\pl_t^{l-2i}v$ ($i=1,2,\cdots$) is based on $\|\pl_t^l v\|_0$ for $\ga=2$ and on $\lt\| \sa (\sa/x)^\nu \pl_t^l v' (t)\rt\|_0$ for $\ga\in(1,2)$. This is the difference between \ef{norm} and \ef{ggnorm}.}

For $\mu>0$, we use the following parabolic approximation to $\ef{419'}_1$:
\be\label{93}\begin{split}
 x \sa \pl_t v  + \left[ \sa^2 \lt( \frac{x}{r}\rt)^{2\ga-2} \lt(\frac{1}{r'}\rt)^\ga    \right]' - 2 \frac{\sa^2}{x}\lt( \frac{x}{r}\rt)^{2\ga-1} \lt(\frac{1}{r'} \rt)^{\ga-1} + \phi \sa x^2 \lt(\frac{x^2}{r^2}\rt) \\
 + \frac{  2-\ga }{\ga-1}
\sa x \lt(\frac{\sa}{x}\rt)' \lt( \frac{x}{r}\rt)^{2\ga-2} \lt(\frac{1}{r'}\rt)^\ga   = \frac{\ga\mu}{x} \lt[(x\sa)^2 \lt(\frac{\sa}{x}\rt)^{2\nu}\lt(\frac{v}{x}\rt)'\rt]',
\end{split}\ee
{which is the general form of $\ef{pe1-3}_1$ for $\ga=2$. This approximation matches the energy estimates and elliptic estimates in the sense that one can derive the uniform estimates with respect to $\mu$.
The existence and uniqueness of the solution to the approximate parabolic problem with the same initial and boundary data as in \ef{419'} can be checked easily as before. To reduce the length of this paper, we will only derive the a priori estimates that guarantees the existence of the solution to problem \ef{419'}.}

\subsection{Energy estimates}
{As for $\ga=2$, taking the $(k+1)-$th time derivative of \ef{93} yields}
\be\label{gg1}\begin{split}
&x \sa \pl_t^{k+2} v   -\left\{ \sa^2 \lt[(2\ga-2)  \lt( \frac{x}{r}\rt)^{2\ga-1} \lt(\frac{1}{r'}\rt)^\ga \frac{\pl_t^k v}{x}   +\ga   \lt( \frac{x}{r}\rt)^{2\ga-2} \lt(\frac{1}{r'}\rt)^{\ga+1} \pl_t^k v'\rt]\right\}' \\
& + 2 \frac{\sa^2}{x}\lt[  (2\ga-1) \lt( \frac{x}{r}\rt)^{2\ga } \lt(\frac{1}{r'}\rt)^{\ga-1} \frac{\pl_t^k v}{x}   +(\ga-1)   \lt( \frac{x}{r}\rt)^{2\ga-1} \lt(\frac{1}{r'}\rt)^{\ga} \pl_t^k v'  \rt]\\
& - 2\nu
\sa x \lt(\frac{\sa}{x}\rt)' \left[  (2\ga-2) \lt( \frac{x}{r}\rt)^{2\ga-1} \lt(\frac{1}{r'}\rt)^\ga \frac{\pl_t^k v}{x}   +\ga   \lt( \frac{x}{r}\rt)^{2\ga-2} \lt(\frac{1}{r'}\rt)^{\ga+1} \pl_t^k v'\right]   \\
= & \lt\{  \sa^2\lt[ (2\ga-2)  W_{11}   +\ga  W_{12}\right]\rt\}'- 2 \frac{\sa^2}{x}\lt[  (2\ga-1) W_{21}   +(\ga-1) W_{22}  \rt]
\\&+2\nu
\sa x \lt(\frac{\sa}{x}\rt)'\left[  (2\ga-2)  W_{11}   +\ga   W_{12}\right]
- \phi \sa x^2 \pl_t^{k+1}\lt(\frac{x^2}{r^2}\rt),
\end{split}\ee
where
\bee\label{}\begin{split}
 &W_{11}=   \pl_t^{k } \lt( \lt( \frac{x}{r}\rt)^{2\ga-1} \lt(\frac{1}{r'}\rt)^\ga  \frac{v}{x}\rt)-  \lt( \frac{x}{r}\rt)^{2\ga-1} \lt(\frac{1}{r'}\rt)^\ga \frac{\pl_t^{k } v } {x}    ,\\
 &W_{12}=   \pl_t^{k } \lt( \lt( \frac{x}{r}\rt)^{2\ga-2} \lt(\frac{1}{r'}\rt)^{\ga+1} v'\rt)-  \lt( \frac{x}{r}\rt)^{2\ga-2} \lt(\frac{1}{r'}\rt)^{\ga+1}\pl_t^{k } v'  ,\\
 &W_{21}=  \pl_t^{k } \lt( \lt( \frac{x}{r}\rt)^{2\ga } \lt(\frac{1}{r'}\rt)^{\ga-1}  \frac{v}{x}\rt)-  \lt( \frac{x}{r}\rt)^{2\ga } \lt(\frac{1}{r'}\rt)^{\ga-1}  \frac{\pl_t^{k } v } {x}  ,\\
&W_{22}=  \pl_t^{k }\left( \lt( \frac{x}{r}\rt)^{2\ga-1} \lt(\frac{1}{r'}\rt)^{\ga} v'\right)- \lt( \frac{x}{r}\rt)^{2\ga-1} \lt(\frac{1}{r'}\rt)^{\ga}\pl_t^{k } v'  .
 \end{split}\eee
{Comparing it with \ef{e1-11} for $\ga=2$, we have to deal with an additional term, the last term on the left-hand side of \ef{gg1}, which does not appear in \ef{e1-11}. To do so, we introduce a weight $(\sa/x)^{2\nu}$ (or equivalently, $\rho_0^{2-\ga}$), which is $1$ for $\ga=2$.} Multiply \ef{gg1} with $k=l$ by $(\sa/x)^{2\nu} \pl_t^{l+1} v$ and integrate the resulting equation with respect to time and space to get
\be\label{geg}\begin{split}
 &\lt\| \sqrt{x \sa } (\sa/x)^\nu  \pl_t^{l+1} v  (t)\rt\|_0^2 +    \lt\| \sa (\sa/x)^\nu \pl_t^l v' (t)\rt\|_0^2 +    \lt\|  (\sa/x)^{1+\nu}  \pl_t^l v (t)\rt\|_0^2 \\
\le &  \widetilde{M}_0+ C t P\lt(\sup_{[0,t]} \widetilde{E}\rt),
\end{split}\ee
provided that $t$ is small. Here $\widetilde{M}_0=P(\widetilde{E}(0,v))$ is determined by the initial density $\rho_0$.
{It should be noted that \ef{geg} is the energy estimate parallelling to \ef{et4}  for $\ga=2$.}

{Based on this energy estimate, we can derive the higher-order spatial derivative of $\pl_t^{l-1}v$ and $\pl_t^{l-2}v$ associated with weights, respectively. Inductively, the weighted spatial derivative of $\pl_t^{l-{2i}+1}v$ and $\pl_t^{l-2i}v$ ($i=2,3,\cdots$) can then be achieved. Next, we use elliptic estimates to
obtain the other norms in the higher-order energy functional. This is done by the interior and boundary estimates}.

\subsection{Elliptic estimates -- interior part}
{To obtain the interior estimates, the key is to choose a suitable cut-off function to separate the whole region into interior and boundary regions such that  the energy norms can be matched in the intermediate regions. For this purpose,} note that
 \bee\label{}\begin{split}
 &\frac{1}{\sa}\left\{ \sa^2 \lt[(2\ga-2)    \frac{\pl_t^k v}{x}   +\ga   \pl_t^k v'\rt]\right\}'
- 2 \frac{\sa }{x}\lt[  (2\ga-1)  \frac{\pl_t^k v}{x}   +(\ga-1)    \pl_t^k v'  \rt] + 2 \ga \nu x \lt(\frac{\sa}{x}\rt)' \pl_t^k v  '     \\
 & =   \ga   H_0   + (6\ga-4  ) (\sa/x)' \pl_t^k v  + 2 \ga \nu x \lt(\frac{\sa}{x}\rt)' \pl_t^k v  ' \\
 & =   \ga \lt[ H_0 + 2\nu x(\sa/x)' \pl_t^k v' -  2\nu (\sa/x)' \pl_t^k v \rt] + (6\ga-4 + 2\nu\ga) (\sa/x)' \pl_t^k v\\
 &= \ga ({x\sa})^{-1}  \lt[(x\sa)^2 \lt( {\sa}/{x}\rt)^{2\nu}\lt( {\pl_t^k v}/{x}\rt)'\rt]' + (6\ga-4 + 2\nu\ga) (\sa/x)' \pl_t^k v,
\end{split}\eee
where $H_0$ is defined in \ef{H0}. Then equation \ef{gg1} reads
 \be\label{gi1}\begin{split}
 &\ga \lt[ H_0 + 2\nu x \lt(\frac{\sa}{x}\rt)' \pl_t^k v' -  2\nu \lt(\frac{\sa}{x}\rt)' \pl_t^k v \rt] \\
 = &x   \pl_t^{k+2} v - (6\ga-4+2\nu \ga) \lt(\frac{\sa}{x}\rt)'\pl_t^k v- \frac{1}{\sa} \lt\{  \sa^2\lt[ (2\ga-2)  W_{11}   +\ga  W_{12}\right]\rt\}'\\
&+ 2 \frac{\sa }{x}\lt[  (2\ga-1) W_{21}   +(\ga-1) W_{22}  \rt]
-2\nu
  x \lt(\frac{\sa}{x}\rt)'\left[  (2\ga-2)  W_{11}   +\ga   W_{12}\right]
+  \phi  x^4 \pl_t^{k+1} \lt(\frac{1}{r^2}\rt)
\\
  &- \frac{1}{\sa}\left\{ \sa^2 \lt[(2\ga-2)  \lt[ \lt( \frac{x}{r}\rt)^{2\ga-1} \lt(\frac{1}{r'}\rt)^\ga -1\rt] \frac{\pl_t^k v}{x}   +\ga   \lt[\lt( \frac{x}{r}\rt)^{2\ga-2} \lt(\frac{1}{r'}\rt)^{\ga+1} -1\rt] \pl_t^k v'\rt]\right\}' \\
& + 2 \frac{\sa }{x}\lt[  (2\ga-1) \lt[ \lt( \frac{x}{r}\rt)^{2\ga } \lt(\frac{1}{r'}\rt)^{\ga-1}-1\rt] \frac{\pl_t^k v}{x}   +(\ga-1)  \lt[ \lt( \frac{x}{r}\rt)^{2\ga-1} \lt(\frac{1}{r'}\rt)^{\ga} -1\rt] \pl_t^k v'  \rt]\\
& - 2\nu
  x \lt(\frac{\sa}{x}\rt)' \left[  (2\ga-2) \lt( \frac{x}{r}\rt)^{2\ga-1} \lt(\frac{1}{r'}\rt)^\ga \frac{\pl_t^k v}{x}   +\ga   \lt[\lt( \frac{x}{r}\rt)^{2\ga-2} \lt(\frac{1}{r'}\rt)^{\ga+1}-1\rt] \pl_t^k v'\right].\\
\end{split}\ee
In the interior region, one can see easily that the main part of the left-hand side of \ef{gi1} is $H_0$. So, we analyze $H_0$ to determine the length of the interior region, $\da_\ga$. Taking the $i$-th ($i\ge 2$) spatial derivative of $H_0$ ($i=0,1$ has been treated in the case of $\ga=2$) leads to
 \be\label{g54}\begin{split}
 H_0^{(i)}-{{H}}_{0 i} = \sa   f^{(i+2)}  + (i+2) \sa'   f^{(i+1)}   -2\sa'  \lt(
 \frac{ f }{x} \rt)^{(i)}
 =: {\widetilde{H}}_{0 i}   , \  \ {\rm where} \ \   f=\pl_t^k v
 \end{split}\ee
 and
 \bee
{{H}}_{0 i} = \sum_{\aa=2}^{i} C_\aa^i    \sa^{(\aa)}   f^{(i+2-\aa)}  + 2\sum_{\aa=1}^{i} C_\aa^i  \sa^{(\aa+1)}   f^{(i+1-\aa)}  - 2\sum_{\aa=1}^{i} C_\aa^i   \sa^{(\aa+1)}      \lt(
 \frac{ f }{x} \rt)^{(i -\aa)}
 \eee
 is the lower-order term. Here $\mathfrak{g}^{(i)}$ denotes $\pl_x^i \mathfrak{g}(x,t)$ for any function $\mathfrak{g}(x,t)$.
Note that
\bee
f^{(j)}=\lt(x \frac{f}{x}  \rt)^{(j)}= x \lt(\frac{f}{x}  \rt)^{(j)} +j\lt(\frac{f}{x}  \rt)^{(j-1)}, \ \  \ \  j=1,2,\cdots.
\eee
Then ${\widetilde{H}}_{0 i}$ $(i=2, 3, \cdots)$ can be rewritten as
\bee\label{}\begin{split}
 {\widetilde{H}}_{0 i}
   = & \sa   x g'' +(i+2) \lt(\sa x\rt)'  g'
    +   i(i+3) \sa'   g , \ \  {\rm where}
\ \  g=\lt(\frac{f}{x}  \rt)^{(i )}=\lt(\frac{\pl_t^{k } v}{x}  \rt)^{(i )};
 \end{split}\eee
or equivalently,
\bee\label{}\begin{split}
  \widetilde{H}_{0i} -(i+2)  \lt(\sa'x-\sa\rt)  g'
    = &  \sa   x g'' +2 (i+2)   \sa  g'
    +   i(i+3)   \sa'   g .
 \end{split}\eee
Therefore, we   obtain that
\bee\label{}\begin{split}
    &\lt\|\za {\widetilde{H}}_{0 i}-(i+2)\za \lt(\sa'x-\sa\rt)  g' \rt\|_0^2 \\
     =
  & \lt\| \za \sa   x g'' \rt\|_0^2+4 (i+2)^2 \lt\| \za \sa  g'\rt\|_0^2
    +   i^2(i+3)^2 \lt\| \za \sa'   g \rt\|_0^2   + 4(i+2) \int \za \sa   x g'' \za \sa  g' dx \\
     &+ 2i(i+3)  \int \za \sa   x g''    \za \sa'   g  dx + 4 i(i+2) (i+3) \int   \za \sa  g'   \za \sa'   g   dx,
 \end{split}\eee
and
\bee\label{}\begin{split}
  & \lt\| \za \sa   x g'' \rt\|_0^2+2\lt[ (i+1)^2+1\rt] \lt\| \za \sa  g'\rt\|_0^2
    +  i (i+3)\lt(  i^2+i-2 \rt) \lt\| \za \sa'   g \rt\|_0^2 \\
= &\lt\|\za \widetilde{H}_{ 0i } -(i+2)\za \lt(\sa'x-\sa\rt)  g' \rt\|_0^2\\
  &    + 4(i+2)\lt[ \int   \za \za '  x  \lt|\sa g'\rt|^2 dx +  \int \za^2 \sa (\sa'x -\sa ) \lt| g'\rt|^2 dx   \rt] \\
     &- 2i(i+3) \lt[ \int \za^2  \sa  \lt(\sa' x -\sa \rt)     g g''  dx - 2 \int \za \za'  \sa^2    g  g' dx  \rt]\\
     & + 2 i(i+2) (i+3) \lt[2 \int \za \za'    \sa \sa'  g^2   dx +\int \za^2     \sa \sa''  g^2   dx\rt]\\
\le & 2\lt\|\za \widetilde{H}_{ 0i }   \rt\|_0^2 + C(i,m_0,m_1)\da \lt[  \lt\| \za \sa   x g'' \rt\|_0^2 +   \lt\|  {\za}  \sa  g'\rt\|_0^2 +   \lt\|  {\za}     g \rt\|_0^2 \rt]\\
&+   C(i,m_0,s_0) \int_\da^{2\da} \lt[(\sa g')^2 + g^2   \rt] dx.
 \end{split}\eee
So, there exist constants $\bar{\da}_{i}=\bar{\da}_{i}(i,m_0,m_1)$ ($i=2,3,\cdots$) such that for   $\da\le  \min\{\da_0/2, \bar{\da}_{i}\}$,
\bee\label{}\begin{split}
  & \frac{1}{2}\lt\| \za \sa   x g'' \rt\|_0^2+ \lt[ (i+1)^2+1\rt] \lt\| \za \sa  g'\rt\|_0^2
    +   \frac{1}{2}  i (i+3)\lt(  i^2+i-2 \rt)  m_0^2 \lt\| \za   g \rt\|_0^2 \\
\le & 2\lt\|\za \widetilde{H}_{ 0i }   \rt\|_0^2  + C(i,m_0,s_0) \int_\da^{2\da} \lt[(\sa g')^2 + g^2   \rt] dx,
 \end{split}\eee
where one has used the fact $\sa'(x)\ge m_0$ on $[0,\da_0]$. Consequently,
 \bee\label{}\begin{split}
  &  \lt\| \za \sa   x g'' \rt\|_0^2+  \lt\| \za \sa  g'\rt\|_0^2
    +   \lt\| \za     g \rt\|_0^2
\le    C (i,m_0 )   \lt\|\za \widetilde{H}_{ 0i }   \rt\|_0^2     +   C(i,m_0,s_0) \int_\da^{2\da} \lt[(\sa g')^2 + g^2   \rt] dx  .
 \end{split}\eee
 It then follows from \ef{g54} that, for each $i\ge 2$,
\be\label{ghaha}\begin{split}
     &\lt\| \za \sa   \lt( {\pl_t^{k } v}   \rt)^{(i+2 )} \rt\|_0^2+  \lt\| \za \lt( {\pl_t^{k } v}   \rt)^{(i+1 )}\rt\|_0^2
    +   \lt\| \za    \lt( \frac{\pl_t^{k } v}{x}   \rt)^{(i )} \rt\|_0^2 \\
\le & C(i,m_0 ) \lt[\lt\|\za  {H}_{ 0  }^{(i)}   \rt\|_0^2 +  \lt\|\za  {H}_{ 0i }   \rt\|_0^2 \rt]     + C(i,m_0,s_0) \int_\da^{2\da} \lt[ \lt|\lt( {\pl_t^{k } v}   \rt)^{(i+1 )} \rt|^2+ \lt|\lt( \frac{\pl_t^{k } v}{x}   \rt)^{(i )}\rt|^2    \rt] dx .
 \end{split}\ee
Choose
\be\label{daga}
\da_\ga=\min\{\da_0/2,\da_1,\da_2,\da_3, \bar{\da}_2, \cdots, \bar{\da}_{ \frac{l-1}{2} } \}.
 \ee
 (Thus $\da_\ga$ depends on the initial density $\rho_0(x)$ and $\ga$.)
With this $\da_\ga$, we can derive from \ef{ghaha}, \ef{gi1} and \ef{geg} that
\be\label{gie}\begin{split}
& \sum_{j=1}^{ \frac{l+1}{2}} \lt\{ \lt\|\za \sa \pl_t^{l-2j +1 } v \rt\|_{j+1}^2
   +\lt\|\za  \pl_t^{l-2j +1 } v \rt\|_{j }^2
   +\lt\|\za \frac {\pl_t^{l-2j +1 } v}{x} \rt\|_{j-1}^2  \rt\}(t) \\
   &+
 \sum_{j=1}^{\frac{l-1}{2}} \lt\{ \lt\|\za \sa \pl_t^{l-2j  } v \rt\|_{j+2}^2
   +\lt\|\za  \pl_t^{l-2j  } v \rt\|_{j +1 }^2
   +\lt\|\za \frac{ \pl_t^{l-2j  } v }{x} \rt\|_{j}^2  \rt\}(t)
\le    \widetilde{M}_0+ C t P\lt(\sup_{[0,t]} \widetilde{E}\rt)  .
 \end{split}\ee
{This completes the interior estimates. Next, we show the boundary estimates using the same argument as that in Section 7.2.}

\subsection{Elliptic estimates -- boundary part}
As before, we can introduce a cut-off function as
$$
 \chi=1 \  \ {\rm on} \ \ [\da_\ga,1], \ \ \chi=0 \ \ {\rm on} \ \  [0, \da_\ga/2], \ \  |\chi'|\le s_0/\da_\ga,
$$
for some constant $s_0$, where $\da_\ga$ is given by \ef{daga}.
Note that in the boundary region, $x\in [\da_\ga/2, 1]$,  the main part of the left-hand side of \ef{gi1} is
$$
B_\ga:=\sa \pl_t^k v''+ (2+2\nu)\sa'\pl_t^k v'.
$$
Taking the $i$-th $(i\ge 0)$ spatial derivative of $B_\ga$ yields
$$
  B^{(i)}_\ga-B_{\ga i }
 =  {\sa   \pl_t^k v ^{(i+2)}  + (i+2+2\nu) \sa'   \pl_t^k v^{(i+1)}} =: \widetilde{B}_{\ga i},
$$
where  $${B}_{\ga i}=\sum_{\aa=2}^{i} C_\aa^i   \sa^{(\aa)}   \pl_t^k v^{(i+2-\aa)}  + 2(1+\nu)\sum_{\aa=1}^{i} C_\aa^i \sa^{(\aa+1)}   \pl_t^k v^{(i+1-\aa)}$$
denotes the lower-order term.
Since for any function $h=h(x,t)$ and integer $i\ge 0$, it holds that
\bee\label{}\begin{split}
 & \lt\| \chi \sa^{3/2+\nu}   h'\rt\|_0^2 + \lt\| \chi \sa^{1/2+\nu} \sa' h \rt\|_0^2 \le
4 \lt\| \chi \sa^{1/2+\nu} \lt(\sa   h' + (i+2+2\nu) \sa' h \rt)  \rt\|_0^2 + C \lt\|  \sa^{1+\nu}  h  \rt\|_0^2,\\
 &\lt\| \chi \sa^{2+\nu}   h'\rt\|_0^2 +  \lt\| \chi \sa^{1+\nu} \sa' h \rt\|_0^2
 \le     4 \lt\| \chi \sa^{1+\nu} \lt(\sa   h' + (i+3+2\nu) \sa' h \rt)  \rt\|_0^2 + C \lt\|  \sa^{3/2+\nu}  h  \rt\|_0^2 ;
 \end{split}\eee
then we have
\be\label{gbe}\begin{split}
 & \sum_{j=1}^{ \frac{l+1}{2} } \lt\{ \lt\|\chi \sa^{3/2+\nu} \pl_t^{l-2j +1 } \pl_x^{ j +1 }  v \rt\|_{0}^2
      +\lt\|\chi \sa^{1/2+\nu} \pl_t^{l-2j +1 } \pl_x^{ j }  v \rt\|_{0}^2 \rt\}\le    \widetilde{M}_0+ C t P\lt(\sup_{[0,t]} \widetilde{E}\rt) \\
      &\sum_{j=1}^{ \frac{l-1}{2} } \lt\{ \lt\|\chi \sa^{2+\nu} \pl_t^{l-2j  } \pl_x^{ j +2 }  v \rt\|_{0}^2
      +\lt\|\chi \sa^{1+\nu} \pl_t^{l-2j  } \pl_x^{ j +1 }  v \rt\|_{0}^2 \rt\}\le    \widetilde{M}_0+ C t P\lt(\sup_{[0,t]} \widetilde{E}\rt).
\end{split}\ee
{This yields the desired is elliptic estimates on the boundary.}

\subsection{Existence for case $1<\ga<2$}
It follows from \ef{ggnorm}, \ef{geg}, \ef{gie} and \ef{gbe} that
\bee\label{}\begin{split}
\widetilde{E}(t)
 \le  \widetilde{M}_0+ C t P\lt(\sup_{s\in[0,t]} \widetilde{E}(s)\rt), \ \ t\in [0,T];
 \end{split}\eee
which implies that for small $T$,
\bee\label{}\begin{split}
\sup_{t\in[0,T]}\widetilde{E}(t)
 \le  2 \widetilde{M}_0.
 \end{split}\eee
{With this a priori estimates,  one can then obtain the local existence of smooth solutions in  the functional space for which $\sup_{t\in[0,T]}\widetilde{{E}}(v, t)<\infty$
 provided that $\widetilde{{E}}(v, 0)<\infty $ ($\widetilde{{E}}(v, 0)$ is determined by the initial data and their spatial derivatives via the equation), by using the parabolic approximation in \ef{93} in a similar way as before.}

\section{Case $\ga >2$}
{In this section, we deal with the case when $\ga>2$, which is easier than
the case when $1<\ga<2$ because the rate of degeneracy of  equation $\ef{419'}_1$ near vacuum states
is lower and less derivatives are needed to control the $H^2$-norm of $v$.}
Set  $$\nu=(2-\ga)(2\ga-2)\in (-1/2,0).$$ The higher-order energy norm is chosen as follows:
 \be\label{less2}\begin{split}
 \widehat{E}(v, t)= & \lt\| \sa (\sa/x)^\nu \pl_t^4 v' (\cdot,t)\rt\|_0^2 +    \lt\|   (\sa/x)^{1+\nu}  \pl_t^4 v (\cdot,t)\rt\|_0^2 \\
      & + \sum_{j=1}^{  2 } \lt\{ \lt\|   \sa^{3/2+\nu} \pl_t^{5-2j  } \pl_x^{ j +1 }  v (\cdot,t)\rt\|_{0}^2
      + \sum_{i=0}^j \lt\|   \sa^{1/2+\nu} \pl_t^{5-2j  } \pl_x^{ i }  v (\cdot,t)\rt\|_{0}^2 \rt\}  \\
      & + \sum_{j=1}^{ 2 } \lt\{ \lt\|  \sa^{2+ \nu} \pl_t^{4-2j  } \pl_x^{ j +2 }  v (\cdot,t)\rt\|_{0}^2
      +\sum_{i=-1}^j \lt\|   \sa^{ 1+ \nu} \pl_t^{4-2j  } \pl_x^{ i +1  }  v (\cdot,t) \rt\|_{0}^2 \rt\}\\
      &+\sum_{j=1}^{ 2} \lt\{ \lt\|\za \sa \pl_t^{5-2j   } v (\cdot,t)\rt\|_{j+1}^2
   +\lt\|\za  \pl_t^{5-2j   } v (\cdot,t)\rt\|_{j }^2
   +\lt\|\za \frac {\pl_t^{5-2j } v}{x} (\cdot,t)\rt\|_{j-1}^2  \rt\}\\
 &+ \sum_{j=1}^{2} \lt\{ \lt\|\za \sa \pl_t^{4-2j  } v (\cdot,t)\rt\|_{j+2}^2
   +\lt\|\za  \pl_t^{4-2j  } v (\cdot,t)\rt\|_{j+1   }^2
   +\lt\|\za \frac{ \pl_t^{4-2j  } v }{x} (\cdot,t)\rt\|_{j }^2  \rt\}.
\end{split}\ee
Here $\za$ is defined in \ef{delta}. It follows from Sobolev embedding \ef{wsv} that
 \bee\label{}\begin{split}
\|v\|_2^2 \le \|v\|_{2-\nu}^2 \le  \sum_{i=0}^3   \| \sa^{1+\nu} \pl_x^i v\|_0  \le C\widehat{E}.
\end{split}\eee
As before, we can show
\bee\label{}\begin{split}
 \widehat{E}(t)
 \le  P(\widehat{E}(0)) + C t P\lt(\sup_{s\in[0,t]} \widehat{E}(s)\rt), \ \ t\in [0,T];
 \end{split}\eee
which implies that for small $T$,
\bee\label{}\begin{split}
\sup_{t\in[0,T]}\widehat{E}(t)
 \le  2  P(\widehat{E}(0)) .
 \end{split}\eee
{With the above estimates,  one can then obtain the local existence of smooth solutions in the  functional space $\sup_{t\in[0,T]}\widehat{E}(t)<\infty$.}

\section{Uniqueness of spherically symmetric motions for the three-dimensional compressible Euler equations}
 For the free-boundary problem of the compressible Euler equations without self-gravitation, we can prove that the uniqueness theorem is true for all values of $\gamma>1$ in a natural functional space for the spherically symmetric motion. (Indeed, a similar argument can be extended to the general three-dimensional motion.) In this case,  problem \ef{419'} becomes
\be\label{419''}\begin{split}
& \rho_0 \lt( \frac{x}{r}\rt)^2\pl_t v   + \pl_x \left[    \lt(\frac{x^2}{r^2}\frac{\rho_0}{\pl_x r}\rt)^\ga    \right]  =0 & {\rm in} & \ \ I \times (0,T],\\
& v(0, t)=0 & {\rm on} & \ \ \{x=0\}\times (0,T],\\
& v(x, 0)= u_0(x)  & {\rm on} & \ \ I \times \{t=0\},
\end{split}
\ee
where the initial density $\rho_0$ satisfies \ef{156}. For problem \ef{419''}, we have the following result:

  \begin{thm}\label{unique2} {\rm(uniqueness for Euler equations)}
 Suppose $\gamma>1$.  Let  $v_1$ and $v_2$ be two  solutions to the problem \eqref{419''} on $[0, T]$ for $T >0$ with
 $$r_i(x, t)=x+\int_0^t v_i(x, s)ds, \  \ i=1, 2.$$
If there exist some positive constants $w_1$,  $w_2$ and $w_3$ such that
 \be\label{a1'}w_1\le r_i'(x, t)\le w_2 \ \  {\rm and} \ \  |v_i'(x, t)|\le w_3, \ \ \  (x, t)\in [0, 1]\times [0, T] , \ \  i=1, 2,\ee
  then
\be\label{v1v2}v_1(x, t)=v_2(x, t), \ \  (x, t)\in [0, 1]\times [0, T]\ee
 provided that $v_1(x, 0)=v_2(x, 0)$ for $x\in [0,1]$.
 \end{thm}

The solution to the spherically symmetric problem of Euler equations in Eulerian coordinates can be obtained from the solution to \eqref{419''}. Denote this solution by
$(\rho, u)(r, t)$ ($0\le r\le R(t)$, $0\le t\le T$). For $({\bf x}, t)\in \mathbb{R}^3\times [0, T]$ with $|{\bf x}|<R(t)$, we set
   \begin{equation}\label{3.1'}\rho({\bf x}, t) = \rho(|{\bf x}|, t), \  \ {\bf u}({\bf x}, t) =  u(|{\bf x}|, t) {\bf x} /|{\bf x}|. \end{equation}
Then $(\rho, {\bf u}, R(t))$ is a solution of the following free boundary problem:
\be\label{2.1euler} \lt\{ \begin{split}
& \pl_t \rho  + {\rm div}(\rho {\bf u}) = 0 ,  & 0<|{\bf x}|< R(t), \ \ t\in [0, T], \\
 &\pl_t (\rho {\bf u})  + {\rm div}(\rho {\bf u}\otimes {\bf u})+\nabla_{\bf x} (\rho^{\gamma}) =0 ,  & 0<|{\bf x}|< R(t), \ \ t\in [0, T],\\
 &\rho>0,  &   0 \le |{\bf x}|< R(t), \ \ t\in [0, T],\\
 & \rho =0,  &  |{\bf x}|= R(t), \  \ t\in [0, T], \\
 &u({\bf 0}, t)=0, & t\in [0, T],\\
 & \mathcal{V}(\partial B_{R(t)})={\bf u}|_{\partial B_{R(t)}}\cdot {\bf n}, & t\in [0, T], \\
&(\rho,{\bf u})({\bf x}, 0)=(\rho_0, {\bf u}_0)(|{\bf x}|), \ & |{\bf x}|\le R_0,
 \end{split} \rt. \ee
 where $R_0>0$ is a constant, $B_{R(t)}=\{{\bf x}\in \mathbb{R}^3: |{\bf x}|<R(t)\}$, $\mathcal{V}(\partial B_{R(t)})$ and ${\bf n}$ represent, respectively, the normal velocity of $\partial B_{R(t)}$ and exterior unit normal vector to $\partial B_{R(t)}$.

As a direct consequence  of Theorem \ref{unique2}, we have
 \begin{coro}\label{corollary} Let $\gamma>1$. The solutions $(\rho, {\bf u}, R(t))$ of the form \ef{3.1'} to the free boundary problem \ef{2.1euler} are unique provided they satisfy the following regularity conditions:
 $$R(t)\in C^1\lt([0, T]\rt) \ \ {\rm and} \ \
 (\rho, {\bf u})\in C^1\cap W^{1, \infty}\lt(\{({\bf x}, t)\in \mathbb{R}^3\times [0, T]: 0<|{\bf x}|\le R(t)\}\rt).$$
 \end{coro}

\noindent{\bf Proof of Theorem \ref{unique2}}. We first present the proof for the case of $\gamma=2$ for simplicity. When $\gamma=2$,  equation  $\ef{419''}_1$ reduces to
$$
x \sa \pl_t^2 r  + \left[\sa^2 \frac{x^2 } {r^2   r'^2}  \right]' -2   \frac{\sa^2}{x} \frac{x^3}{r^3 r'}  =0  \ \ {\rm in} \ \  I\times (0, T].
$$
Set
$$\theta(x,t)=r_2(x,t)-r_1(x,t),$$
then
\be\label{1}\begin{split}
x \sa \pl_t^2 \ta  - \left[\sa^2 \lt( \frac{x^2 } {r_1^2   {r_1'}^{2}}  - \frac{x^2 }{r_2^2   {r_2'}^{2}}\rt) \right]' + 2   \frac{\sa^2}{x} \lt( \frac{x^3}{r_1^3 r'_1} -\frac{x^3}{r_2^3 r'_2}\rt) =0  \ \ {\rm in} \ \  I\times (0, T].
\end{split}\ee

Multiplying \eqref{1} by $\pl_t \ta$ and integrating the resulting equation with respect to $x$, we have
\bee\label{2}\begin{split}
 \frac{1}{2}\frac{d}{dt}\int x \sa \lt(\pl_t \ta \rt)^2 dx
 =  & -\int  \sa^2 \lt( \frac{x^2 } {r_1^2   {r_1'}^{2}}  - \frac{x^2 }{r_2^2   {r_2'}^{2}}\rt)  \lt(\pl_t \ta'\rt) dx
  \\
  &- 2     \int \frac{\sa^2}{x} \lt( \frac{x^3}{r_1^3 r'_1} -\frac{x^3}{r_2^3 r'_2}\rt) \lt(\pl_t \ta \rt)dx.
\end{split}\eee
Note that
\bee
   \frac{x^2 } {r_1^2   {r_1'}^{2}}  - \frac{x^2 }{r_2^2   {r_2'}^{2}}= \mathcal{A}_1 \ta' +  \mathcal{A}_2 (\ta/x)  \ \ {\rm and} \ \
   \frac{x^3}{r_1^3 r'_1} -\frac{x^3}{r_2^3 r'_2} = \mathcal{A}_3 \ta'  +\mathcal{A}_4 (\ta /x),
 \eee
where
\bee\begin{split}
 & \mathcal{A}_1 = \lt(\frac{x}{r_1}\rt)^2 \lt(\frac{1}{r_1'}+ \frac{1}{r_2'}\rt) \frac{1}{r_1' r_2'}, &
  \mathcal{A}_2 = \lt(\frac{1}{r_2'}\rt)^2 \lt(\frac{x}{r_1}+ \frac{x}{r_2}\rt) \frac{x}{r_1 }\frac{x}{r_2 },\\
  &\mathcal{A}_3 = \lt(\frac{x}{r_1}\rt)^3  \frac{1}{r_1' r_2'}, &
  \mathcal{A}_4=  \frac{1}{r_2'}  \lt[\lt(\frac{x}{r_1} \rt)^2 + \frac{x}{r_1} \frac{x}{r_2} +\lt( \frac{x}{r_2}\rt)^2 \rt]\frac{x}{r_1 }\frac{x}{r_2 }.
\end{split}\eee
Since $r_i(0, t)=0$ and $v_i(0, t)=0$ ($i=1,\ 2$) for $t\in [0, T]$, the bounds in \ef{a1'}  give the following bounds:
$$ w_1\le r_i (x, t)/x\le w_2 \ \ {\rm and } \ \ |v_i(x,t)/x|\le w_3,\ \ (x, t)\in [0, 1]\times [0, T], \ \ i=1, 2. $$
Then, using the integration by parts and the Cauchy inequality, we can get that
\be\label{2}\begin{split}
 &\frac{1}{2}\frac{d}{dt}\int  \lt\{ x \sa \lt(\pl_t \ta \rt)^2
 + \sa^2 \lt[ \mathcal{A}_1 (\ta')^2 +  2\mathcal{A}_2 (\ta/x) \ta' + 2\mathcal{A}_4 (\ta/x)^2 \rt] \rt\} dx\\
 = & \int  \sa^2 \lt[ \frac{1}{2}(\pl_t \mathcal{A}_1) (\ta')^2  +  (\pl_t \mathcal{A}_2) (\ta/x) \ta' +   (\pl_t \mathcal{A}_4)(\ta/x)^2    \rt]dx \\
 & + \int \sa^2 (\mathcal{A}_2-2 \mathcal{A}_3) (\pl_t \ta/x) \ta' dx \\
 \le & C(w_1,w_2,w_3) \int  \sa^2 \lt[   (\ta')^2  + (\ta/x)^2    \rt]dx   + 2 w_3 \int \sa^2 \lt|\mathcal{A}_2-2 \mathcal{A}_3\rt| \lt|\ta'\rt| dx ,
\end{split}\ee
where $C(w_1,w_2,w_3)$ is a positive constant depending on $w_1,w_2,w_3$; because
\bee\label{}\begin{split}
 \lt|\pl_t \mathcal{A}_1\rt| + \lt|\pl_t \mathcal{A}_2\rt|+\lt|\pl_t \mathcal{A}_4\rt| \le C(w_1,w_2) \lt(|v_1/x|+|v_2/x| + |v_1'| + |v_2'|\rt) \le C(w_1,w_2, w_3)
\end{split}\eee
and
\bee\label{}\begin{split}
 \lt|\pl_t \ta/x\rt| = \lt|(v_2-v_1)/x \rt| \le 2 w_3.
\end{split}\eee

To estimate \ef{2}, we need the following a priori assumption: there exists a small positive constant $\ea_0$ such that
\be\label{aa1} |\ta'(x,t)|+ |(\ta/x)(x,t)|\le \ea_0  \ \ {\rm for } \ \ {\rm all} \ \  (x, t)\in (0, 1)\times [0, T].\ee
Thus, a simple calculation  yields that
\bee\begin{split}
 & \mathcal{A}_1 \ge (2- C(w_1,w_2) \ea_0) \lt(\frac{x}{r_1}\rt)^2 \lt( \frac{1}{r_1'} \rt)^3 \ge \frac{7}{4}\lt(\frac{x}{r_1}\rt)^2 \lt( \frac{1}{r_1'} \rt)^3 ,\\
  &
  \mathcal{A}_2 \le (2+ C(w_1,w_2)  \ea_0)   \lt(\frac{x}{r_1}\rt)^3 \lt( \frac{1}{r_1'} \rt)^2 \le \frac{9}{4}\lt(\frac{x}{r_1}\rt)^3 \lt( \frac{1}{r_1'} \rt)^2 ,\\
  &
  \mathcal{A}_4 \ge (3- C(w_1,w_2) \ea_0)  \lt(\frac{x}{r_1}\rt)^4   \frac{1}{r_1'} \ge \frac{11}{4} \lt(\frac{x}{r_1}\rt)^4   \frac{1}{r_1'} ;
\end{split}\eee
which implies that for $(x, t)\in (0, 1)\times [0, T]$,
 \be\label{3}\begin{split}
 \mathcal{A}_1 (\ta')^2 +  2\mathcal{A}_2 (\ta/x) \ta' + 2\mathcal{A}_4 (\ta/x)^2 \ge &\frac{1}{4} \lt(\frac{x}{r_1}\rt)^2 \lt( \frac{1}{r_1'} \rt)^3 (\ta')^2 +  \frac{17}{8} \lt(\frac{x}{r_1}\rt)^4   \frac{1}{r_1'} \lt(\frac{\ta}{x}\rt)^2\\
 \ge & k_1 (\ta')^2 + k_2 (\ta/x)^2 .
\end{split}\ee
Here $k_1$ and $k_2$ are positive constants depending on $w_1$ and $w_2$. We use the cancelation of the leading terms to estimate of $\int \sa^2 \lt|\mathcal{A}_2-2 \mathcal{A}_3\rt| \lt|\ta'\rt| dx$.  Note  that
\bee\begin{split}
  &
  \mathcal{A}_2 =2   \lt(\frac{x}{r_1}\rt)^3 \lt( \frac{1}{r_1'} \rt)^2+ C(w_1,w_2) \lt(\lt|\frac{\ta}{x}\rt|+|\ta'|\rt) , \\
  &\mathcal{A}_3 =   \lt(\frac{x}{r_1}\rt)^3 \lt( \frac{1}{r_1'} \rt)^2+ C(w_1,w_2) |\ta'|.
\end{split}\eee
It then follows from the Cauchy's inequality that
 \be\label{4}\begin{split}
 \int \sa^2 \lt|\mathcal{A}_2-2 \mathcal{A}_3\rt| \lt|\ta'\rt| dx  \le
  C(w_1,w_2)  \int  \sa^2 \lt[   (\ta')^2  + (\ta/x)^2    \rt]dx .
\end{split}\ee
In view of \ef{2}, \ef{3} and \ef{4}, we see that
\bee\label{5}\begin{split}
 &\frac{1}{2} \int  \lt[ x \sa \lt(\pl_t \ta \rt)^2
 +    k_1 ( \sa \ta')^2 +  k_2 (\sa \ta/x)^2 \rt]  dx (t) \\
 \le & \frac{1}{2} \int  \lt\{ x \sa \lt(\pl_t \ta \rt)^2
 + \sa^2 \lt[ \mathcal{A}_1 (\ta')^2 +  2\mathcal{A}_2 (\ta/x) \ta' + 2\mathcal{A}_4 (\ta/x)^2 \rt] \rt\} dx (t=0)\\
 & +C(w_1,w_2,w_3) \int_0^t \int   \lt[   (\sa \ta')^2  + (\sa \ta/x)^2    \rt]dx \\
 \le & C(w_1,w_2,w_3)\int_0^t \int   \lt[   (\sa \ta')^2  + (\sa \ta/x)^2    \rt]dx,
\end{split}\eee
 provided that $v_1(x,0)=v_2(x,0)$. So, it gives from  Grownwall's inequality that for $t\in[0, T]$,
 \bee\label{5}\begin{split}
 &  \int  \lt[ x \sa \lt(\pl_t \ta \rt)^2
 +    k_1 ( \sa \ta')^2 +  k_2 (\sa \ta/x)^2 \rt](x, t)  dx \\
  \le & \exp\lt\{C(w_1,w_2,w_3)T\rt\} \int  \lt[ x \sa \lt(\pl_t \ta \rt)^2
 +    k_1 ( \sa \ta')^2 +  k_2 (\sa \ta/x)^2 \rt] (x, 0) dx
  =   0  ,
\end{split}\eee
 if $v_1(x,0)=v_2(x,0)$; which implies directly that
 $$v_2-v_1=\pl_t \ta = \ta' = \ta /x =0, \ \ (x,t)\in (0,1)\times [0,T] ,$$
because of $\sa(x)>0$ for all $x\in(0,1)$. This verifies the a priori assumption \ef{aa1} and completes the proof of Theorem \ref{unique2} when $\ga=2$.

When $\ga\neq 2$,  equation  $\ef{419''}_1$ reduces to
\bee\label{}\begin{split}
 x \sa \pl_t^2 r  + \left[ \sa^2 \lt( \frac{x}{r}\rt)^{2\ga-2} \lt(\frac{1}{r'}\rt)^\ga    \right]' - &2 \frac{\sa^2}{x}\lt( \frac{x}{r}\rt)^{2\ga-1} \lt(\frac{1}{r'} \rt)^\ga
 \\&+ \frac{  2-\ga }{\ga-1}
\sa x \lt(\frac{\sa}{x}\rt)' \lt( \frac{x}{r}\rt)^{2\ga-2} \lt(\frac{1}{r'}\rt)^\ga     =0,
\end{split}\eee
which implies that
\bee\label{}\begin{split}
x \sa \pl_t^2 \ta  -& \left[ \sa^2 \lt( \frac{x}{r_1}\rt)^{2\ga-2} \lt(\frac{1}{r'_1}\rt)^\ga  -\sa^2  \lt( \frac{x}{r_2}\rt)^{2\ga-2} \lt(\frac{1}{r'_2}\rt)^\ga \right]' \\
&+ 2 \frac{\sa^2}{x}\lt[\lt( \frac{x}{r_1}\rt)^{2\ga-1} \lt(\frac{1}{r'_1} \rt)^\ga-\lt( \frac{x}{r_2}\rt)^{2\ga-1} \lt(\frac{1}{r'_2} \rt)^\ga\rt]\\
&-\frac{  2-\ga }{\ga-1}
\sa x \lt(\frac{\sa}{x}\rt)' \lt[\lt( \frac{x}{r_1}\rt)^{2\ga-2} \lt(\frac{1}{r'_1}\rt)^\ga-\lt( \frac{x}{r_2}\rt)^{2\ga-2} \lt(\frac{1}{r'_2}\rt)^\ga\rt]=0.
\end{split}\eee
Set
$$\nu:=(2-\ga)/(2\ga-2).$$
Multiply the preceding equation with $(\sa/x)^{2\nu} \pl_t \ta$ and integrate the product with respect to time and space. Then using the same argument as to the proof of $\ga=2$, we can show that \ef{v1v2} is true for $\ga\neq2$.

This finishes the proof of Theorem \ref{unique2}.

\vskip 1cm

 \centerline{\bf Acknowledgements}. 

Xin's research was partially supported by  the Zheng Ge Ru Foundation, and Hong
Kong RGC Earmarked Research Grants CUHK-4041/11P, CUHK-4048/13P, a Focus Area Grant from
The Chinese University of Hong Kong, and a grant from Croucher Foundation. Zeng's research was partially supported by NSFC grant \#11301293/A010801.

\vskip 1cm

\section*{Appendix}
In this appendix, we verify \ef{3lk}, \ef{hzz},  \ef{ho3} and \ef{lb0}.

 \vskip 0.5cm
\noindent{\bf Verification of  \ef{3lk}}. For $\mathcal{R}_0$, it follows from \ef{norm} and \ef{egn} that
\bee\label{}\begin{split}
  \lt\| \mathcal{R}_0(t)\rt\|_0
  \le  \int_0^t \lt(  \lt\|\lt(\frac{v}{x}\rt)'\rt\|_0+ \lt\|v''\rt\|_0 \rt)ds
  \le C t \sup_{[0,t]} \sqrt{E} ,
    \end{split}\eee
\bee\label{rrr}\begin{split}
  \lt\|\sa \mathcal{R}_0(t)\rt\|_{L^\iy } \le &
 \lt\|\int_0^t  \sa \lt(\frac{v}{x}\rt)' ds\rt\|_{L^\iy }+ \lt\|\int_0^t \sa v'' ds\rt\|_{L^\iy }\\
  \le &  C\lt\|\int_0^t  x \lt(\frac{v}{x}\rt)' ds\rt\|_{L^\iy }+ \lt\|\int_0^t \sa v'' ds\rt\|_{L^\iy }\\
  \le & C\int_0^t \lt(\lt\|v'- \frac{v}{x} \rt\|_{L^\iy}+ C \lt\|\sa v''\rt\|_{L^\iy}\rt)ds \le  Ct  \sup_{[0,t]} \sqrt{E}.
    \end{split}\eee
Next, we will show $\ef{3lk}_2$. It follows from $\ef{3jk}_1$, \ef{egn} and $\ef{3egn}_1$ that for  $p\in(1,\iy)$,
\bee\label{}\begin{split}
  \lt\| \mathfrak{L}_{0} (t) \rt\|_{0} \le    \lt\|  v''\rt\|_{0} + \lt\| (v/x)'\rt\|_{0}+  \lt\| \mathcal{R}_0\rt\|_0
\lt\|\mathfrak{J}_{0}\rt\|_{L^\iy}
 \le     C \sqrt{E(t)} + Ct \sup_{[0,t]} {E},
 \end{split}\eee
\bee\label{}\begin{split}
  \lt\|\sa \mathfrak{L}_{0} (t)\rt\|_{L^\iy } \le     \lt\|\sa v''\rt\|_{L^\iy} + C\lt\| v'- \frac{v}{x}\rt\|_{L^\iy}+  \lt\| \sa \mathcal{R}_0\rt\|_{L^\iy }
\lt\|\mathfrak{J}_{0}\rt\|_{L^\iy}
 \le     C \sqrt{E(t)} + Ct \sup_{[0,t]} {E} ,
 \end{split}\eee
\bee\label{}\begin{split}
  \lt\|  \sa \mathfrak{L}_{1} (t)\rt\|_{L^p } \le &   \lt\| \sa \pl_t v''\rt\|_{L^p} + C \lt\| \pl_t v ' - \frac{ \pl_t v}{x}   \rt\|_{L^p}+   \lt\|\sa \mathcal{R}_0\rt\|_{L^\iy }
\lt\|  \mathfrak{J}_{1}\rt\|_{L^p}
\\& + \lt\|   \sa  \mathfrak{L}_{0}\rt\|_{L^\iy } \lt\|  \mathfrak{J}_{0}\rt\|_{L^\iy}
 \le     C P \lt(\sqrt{E(t)}\rt) + Ct P \lt(\sup_{[0,t]} \sqrt{E}\rt),
 \end{split}\eee
\bee\label{}\begin{split}
  \lt\|\za \sa \mathfrak{L}_{1}(t) \rt\|_{L^\iy} \le &   \lt\|\za \sa \pl_t v''\rt\|_{L^\iy} + C\lt\|\za  \lt(\pl_t v ' - \frac{ \pl_t v}{x}  \rt)\rt\|_{L^\iy}+   \lt\|\sa \mathcal{R}_0\rt\|_{L^\iy }
\lt\|\za \mathfrak{J}_{1}\rt\|_{L^\iy} \\
& + \lt\|  \sa \mathfrak{L}_{0}\rt\|_{L^\iy } \lt\|  \mathfrak{J}_{0}\rt\|_{L^\iy}
\le      C P \lt(\sqrt{E(t)}\rt) + Ct P \lt(\sup_{[0,t]} \sqrt{E}\rt),
 \end{split}\eee
\bee\label{}\begin{split}
  \lt\| \sa \mathfrak{L}_{2} (t) \rt\|_{0} \le &  \lt\| \sa \pl_t^2 v''\rt\|_{0} + C\lt\|   \pl_t^2 v ' - \frac{ \pl_t^2 v}{x}  \rt\|_{0}+  \lt\| \sa \mathcal{R}_0\rt\|_{L^\iy }
\lt\|  \mathfrak{J}_{2}\rt\|_{0}
+ \lt\| \sa  \mathfrak{L}_{0}\rt\|_{L^\iy } \lt\|  \mathfrak{J}_{1}\rt\|_{0}\\
& +
\lt\|   \sa \mathfrak{L}_{1}\rt\|_{0} \lt\|  \mathfrak{J}_{0}\rt\|_{L^\iy}
\le  C P \lt(\sqrt{E(t)}\rt) + Ct P \lt(\sup_{[0,t]} \sqrt{E}\rt).
 \end{split}\eee
 We now turn to the proof of $\ef{3lk}_3$. It follows from \ef{tnorm}, \ef{tegn}, $\ef{3jk}_2$ and $\ef{3egn}_{2,3}$ that
\bee\label{}\begin{split}
  \lt\| \za \mathfrak{L}_{0} (t) \rt\|_{0}^2 \le   \lt( \lt\| \za v''\rt\|_{0} + \lt\| \za (v/x)'\rt\|_{0}+  \lt\| \mathcal{R}_0\rt\|_0
\lt\|\mathfrak{J}_{0}\rt\|_{L^\iy}\rt)^2
 \le M_0+ C t P\lt(\sup_{[0,t]} E\rt)    ,
 \end{split}\eee
\bee\label{}\begin{split}
   \lt\|  \sa \mathfrak{L}_{0} (t) \rt\|_{L^p }^2  \le &   \lt( \lt\| \sa v''\rt\|_{L^p}   + C\lt\|   v'-    {v}/{x}\rt\|_{L^p} +   \lt\| \sa \mathcal{R}_0\rt\|_{L^\iy }
\lt\|  \mathfrak{J}_{0}\rt\|_{L^p}  \rt)^2 \\
 \le  &     M_0+ C t P\lt(\sup_{[0,t]} E\rt),
 \end{split}\eee
\bee\label{}\begin{split}
  \lt\|\za \sa \mathfrak{L}_{0} (t)\rt\|_{L^\iy }^2  \le &   \lt( \lt\|\za \sa v''\rt\|_{L^\iy}   + C\lt\| \za v'- \za {v}/{x}\rt\|_{L^\iy}   +
\lt\| \sa \mathcal{R}_0\rt\|_{L^\iy} \lt\|\za \mathfrak{J}_{0}\rt\|_{L^\iy}  \rt)^2\\
 \le  &  M_0+ C t P\lt(\sup_{[0,t]} E\rt),
 \end{split}\eee
\bee\label{}\begin{split}
  \lt\| \sa \mathfrak{L}_{1} (t)\rt\|_{0}^2 \le &   \lt( \lt\| \sa \pl_t v''\rt\|_{0} + C \lt\| \pl_t v ' -  { (\pl_t v)}/{x}   \rt\|_{0} +
\lt\| \sa \mathcal{R}_0\rt\|_{L^\iy }\lt\|  \mathfrak{J}_{1}\rt\|_{0}\rt. \\
&\lt.
  + \lt\|  \sa \mathfrak{L}_{0}\rt\|_{L^4} \lt\|  \mathfrak{J}_{0}\rt\|_{L^4}  \rt)^2
 \le      M_0+ C t P\lt(\sup_{[0,t]} E\rt),
 \end{split}\eee
\bee\label{}\begin{split}
  \lt\|\za \sa \mathfrak{L}_{2} (t)\rt\|_{0}^2 \le &   \lt( \lt\| \za \sa \pl_t^2 v''\rt\|_{0}  + C\lt\| \za\lt(  \pl_t^2 v ' -  ({ \pl_t^2 v})/{x}\rt)  \rt\|_{0} + \lt\| \sa \mathcal{R}_0\rt\|_{L^\iy }
\lt\| \za \mathfrak{J}_{2}\rt\|_{0}
\rt.\\
&  \lt.
+ \lt\| \za  \sa \mathfrak{L}_{0}\rt\|_{L^\iy }  \lt\|  \mathfrak{J}_{1}\rt\|_{0} +
\lt\|  \sa \mathfrak{L}_{1}\rt\|_0 \lt\| \za \mathfrak{J}_{0}\rt\|_{L^\iy}^2\rt)^2
\le  M_0+ C t P\lt(\sup_{[0,t]} E\rt).
 \end{split}\eee

\vskip 0.2cm
\noindent {\bf Verification  of \ef{hzz}.} Note that
\bee\label{}\begin{split}
 \lt\| \chi \lt(\sa   h' + i \sa' h \rt)  \rt\|_0^2 =  \lt\| \chi \sa   h'\rt\|_0^2 + i^2 \lt\| \chi \sa' h \rt\|_0^2 +2i \int \chi^2 \sa \sa' h h' dx
 \end{split}\eee
and
\bee\label{}\begin{split}
   & 2i \int \chi^2 \sa \sa' h h' dx = - i\int \lt(\chi^2 \sa \sa'\rt)' h^2 dx \\ \ge & - i \lt\| \chi \sa'   h \rt\|_0^2 - C(i) \lt\| \chi \sa^{1/2}  h  \rt\|_0^2 - C(i,\da)\int_{\da/2}^\da   \chi \sa  h^2 dx \\
   \ge  & - i \lt\| \chi \sa'   h \rt\|_0^2 - C (i,\da) \lt\|   \sa^{1/2}  h  \rt\|_0^2 .
 \end{split}\eee
 Then we have for $i\ge 2$ that
\bee\label{}\begin{split}
  \lt\| \chi \sa   h'\rt\|_0^2 +  \lt\| \chi \sa' h \rt\|_0^2 \le \lt\| \chi \sa   h'\rt\|_0^2 + i(i-1) \lt\| \chi \sa' h \rt\|_0^2
 \le     \lt\| \chi \lt(\sa   h' + i \sa' h \rt)  \rt\|_0^2 + C \lt\|  \sa^{1/2}  h  \rt\|_0^2 .
 \end{split}\eee
This is $\ef{hzz}_1$. Next, we will show $\ef{hzz}_2$.
Note that
\bee\label{}\begin{split}
 \lt\| \sa^{1/2}\chi \lt(\sa   h' + i \sa' h \rt)  \rt\|_0^2 =  \lt\| \chi \sa^{3/2}   h'\rt\|_0^2 + i^2 \lt\| \chi \sa^{1/2} \sa' h \rt\|_0^2 +2i \int \chi^2 \sa^2 \sa' h h' dx
 \end{split}\eee
and
\bee\label{}\begin{split}
   & 2i \int \chi^2 \sa^2 \sa' h h' dx = - i\int \lt(\chi^2 \sa^2 \sa'\rt)' h^2 dx \\ \ge & - 2 i \lt\| \chi \sa^{1/2}\sa'   h \rt\|_0^2 - C(i) \lt\| \chi \sa   h  \rt\|_0^2 - C(i,\da)\int_{\da/2}^\da   \chi \sa^2  h^2 dx \\
   \ge  & - 2 i \lt\| \chi \sa^{1/2}\sa'   h \rt\|_0^2 - C (i,\da) \lt\|   \sa   h  \rt\|_0^2 .
 \end{split}\eee
Then, one has  for $i\ge 2$
\bee\label{}\begin{split}
 \lt\| \chi \sa^{3/2}   h'\rt\|_0^2\le & \lt\| \chi \sa^{3/2}   h'\rt\|_0^2 + i(i-2) \lt\| \chi \sa^{1/2} \sa' h \rt\|_0^2
 \le     \lt\| \chi \sa^{1/2} \lt(\sa   h' + i \sa' h \rt)  \rt\|_0^2 + C \lt\|  \sa  h  \rt\|_0^2 .
 \end{split}\eee
Since the estimate on $\lt\| \chi \sa^{1/2} \sa' h \rt\|_0$ is missed, one has to
use Minkowski's inequality to find it again. That is,
\bee\label{}\begin{split}
\lt\| \chi \sa^{1/2} \sa' h \rt\|_0^2 \le i^2  \lt\| \chi \sa^{1/2} \sa' h \rt\|_0^2 \le  & \lt(\lt\| \chi \sa^{1/2} \lt(\sa   h' + i \sa' h \rt)  \rt\|_0 + \lt\| \chi \sa^{3/2}   h'\rt\|_0\rt)^2\\
 \le & 2\lt(\lt\| \chi \sa^{1/2} \lt(\sa   h' + i \sa' h \rt)  \rt\|_0^2 + \lt\| \chi \sa^{3/2}   h'\rt\|_0^2\rt) \\
 \le & 3 \lt\| \chi \sa^{1/2} \lt(\sa   h' + i \sa' h \rt)  \rt\|_0^2 + C \lt\|  \sa  h  \rt\|_0^2,
 \end{split}\eee
provided that $i\ge 1$. Therefore, we obtain
\bee\label{}\begin{split}
 \lt\| \chi \sa^{3/2}   h'\rt\|_0^2 + \lt\| \chi \sa^{1/2} \sa' h \rt\|_0^2 \le
4 \lt\| \chi \sa^{1/2} \lt(\sa   h' + i \sa' h \rt)  \rt\|_0^2 + C \lt\|  \sa  h  \rt\|_0^2.
 \end{split}\eee

\vskip 1cm
\noindent {\bf Verification of   \ef{ho3}. }
 In view of \ef{norm}, one has
\bee\label{}\begin{split}
  \lt\| \lt( \sa^{3/2}  \pl_t  v'',    \ \sa^{3/2} \pl_t^3 v'
\rt)(\cdot,t)\rt\|_{1}^2
 \le  C E(t),
\end{split}\eee
which implies
\bee\label{}\begin{split}
 & \lt\| \lt( \sa^{3/2}  \pl_t v'',    \ \sa^{3/2} \pl_t^3 v'
\rt)(\cdot,t)\rt\|_{L^\iy}^2  \le C{E(t)}.
\end{split}\eee
Using the embedding $W^{1,4/3}(\mathbb{R}) \subset W^{3/4,2}(\mathbb{R})$, one has
\bee\label{}\begin{split}
  \lt\|  \sa^{1/2} \pl_t^3 v(t) \rt\|_{L^\iy} \le  &  C\lt\|\sa^{1/2} \pl_t^3 v \rt\|_{3/4}=C\lt\|\sa^{1/2} \pl_t^3 v \rt\|_{W^{3/4,2}}\le C \lt\|\sa^{1/2} \pl_t^3 v \rt\|_{W^{1,4/3}}\\
  \le & C \lt\|\sa^{1/2} \pl_t^3 v \rt\|_{L^{4/3}} +C \lt\|\lt(\sa^{1/2} \pl_t^3 v \rt)'\rt\|_{L^{4/3}} \le C \sqrt{E(t)},
\end{split}\eee
since
\bee\label{}\begin{split}
 \lt\|\sa^{1/2} \pl_t^3 v (t)\rt\|_{L^{4/3}}\le  \lt\|\sa^{1/2} \pl_t^3 v \rt\|_{L^{2}}\|1\|_{L^4} \le C \lt\|  \pl_t^3 v \rt\|_{0} \le C \sqrt{E(t)}
\end{split}\eee
and
\bee\label{}\begin{split}
 \lt\|\lt(\sa^{1/2} \pl_t^3 v \rt)' (t)  \rt\|_{L^{4/3}}\le  & \lt\| \sa^{1/2} \pl_t^3 v '  \rt\|_{L^{4/3}} + C\lt\| \sa^{-1/2} \pl_t^3 v   \rt\|_{L^{4/3}}\\
  \le & \lt\| \sa^{1/2} \pl_t^3 v '  \rt\|_{0} + C\lt\| \sa^{-1/2} \rt\|_{L^{5/3}}\lt\| \pl_t^3 v   \rt\|_{L^{20/3}} \\
  \le & \lt\| \sa^{1/2} \pl_t^3 v '  \rt\|_{0} + C\lt\| \pl_t^3 v   \rt\|_{1/2}  \le C \sqrt{E(t)}.
\end{split}\eee
Here we have used the H$\ddot{o}$lder inequality and the fact $\|\cdot\|_{L^p} \le C \|\cdot\|_{1/2}$ for any $p\in(1,\iy)$.
Similarly,
$$ \lt\|  \sa^{1/2} \pl_t  v' \rt\|_{L^\iy}   \le C \sqrt{E(t)}.
$$

\vskip 1cm
\noindent {\bf Verification of   \ef{lb0}. }
One can obtain $\ef{lb0}_1$ by using \ef{egn}, \ef{ho3} and $\ef{lb1}_{1,2}$, since
\bee\label{}\begin{split}
 & \lt\|   \sa^{1/2} \mathfrak{J}_1 (t)\rt\|_{L^\iy}\le C \|\pl_t v/x\|_{L^\iy} + \lt\| \sa^{1/2} \pl_t v' \rt\|_{L^\iy}  +   \lt\|   \mathfrak{J}_0 \rt\|_{L^\iy}^2  \le
P\lt(\sqrt{E (t)}\rt),
\end{split}\eee
\bee\label{}\begin{split}
 \lt\| \sa^{3/2} \mathfrak{L}_{1}(t) \rt\|_{L^\iy} \le & C\lt(\lt\|\sa^{3/2} \pl_t v''\rt\|_{L^\iy} +  \lt\|  \sa^{1/2} \lt( \pl_t v ' - \pl_t v/x\rt) \rt\|_{L^\iy} +\lt\|  \pl_t v  \rt\|_{L^\iy}\rt. \\
  &\lt. + \lt\|  \sa  \mathcal{R}_0\rt\|_{L^\iy }  \lt\| \sa^{1/2}\mathfrak{J}_{1}\rt\|_{L^\iy} + \lt\|   \sa  {\mathfrak{L}}_0\rt\|_{L^\iy}  \lt\|\mathfrak{J}_{0}\rt\|_{L^\iy} \rt) \\
\le  & C P\lt(\sqrt{E(t)}\rt) +C t P\lt( \sup_{[0,t]} \sqrt{E}\rt).
\end{split}\eee
For $\ef{lb0}_2$, it follows from \ef{f0}, \ef{tnorm}, \ef{tjk} and $\|\cdot\|_{L^\iy}\le \|\cdot\|_1$ that
\bee\label{}\begin{split}
\lt\|  \mathfrak{J}_0(t)\rt\|_{L^\iy}^2\le C \lt(\lt\|   {  v} /x \rt\|_{L^\iy}^2 + \lt\|      v'\rt\|_{L^\iy}^2\rt)\le C \lt(\lt\|   {  v} /x \rt\|_{1}^2 + \lt\|      v\rt\|_{2}^2\rt)
\le  M_0+ C t P\lt(\sup_{[0,t]} E\rt),
\end{split}\eee
\bee\label{}\begin{split}
\lt\|  \mathfrak{J}_1(t)\rt\|_{0}^2
\le  M_0+ C t P\lt(\sup_{[0,t]} E\rt),
\end{split}\eee
\bee\label{}\begin{split}
\lt\|   \mathfrak{J}_2(t)\rt\|_{0}^2\le & C \lt(\lt\|  \pl_t^2 {  v} /x\rt\|_{0}^2 + \lt\|     \pl_t^2    v'\rt\|_{0}^2 +
\lt\|  \mathfrak{J}_0\rt\|_{L^\iy}^2\lt\|  \mathfrak{J}_1\rt\|_{0}^2\rt)
\le   M_0+ C t P\lt(\sup_{[0,t]} E\rt),
\end{split}\eee
and
\bee\label{}\begin{split}
 \lt\| \sa \mathfrak{L}_{0} (t)\rt\|_{L^\iy}^2 \le & C\lt(\lt\|\sa v''\rt\|_{L^\iy}^2 +  \lt\|  v '-v/x \rt\|_{L^\iy}^2 + \lt\| \sa  \mathcal{R}_0\rt\|_{L^\iy}^2 \lt\|\mathfrak{J}_{0}\rt\|_{L^\iy }^2 \rt)\\
\le & C\lt(\lt\|\sa v \rt\|_{3}^2 +  \lt\|  v  \rt\|_{2}^2+ \lt\|  v /x  \rt\|_{1}^2 + \lt\|  \sa  \mathcal{R}_0\rt\|_{L^\iy}^2 \lt\|\mathfrak{J}_{0}\rt\|_{L^\iy }^2 \rt)
 \le   M_0+ C t P\lt(\sup_{[0,t]} E\rt),
\end{split}\eee
\bee\label{}\begin{split}
 \lt\| \sa  \mathfrak{L}_{1}(t) \rt\|_{0}^2
\le   M_0 + CtP\lt(\sup_{[0,t]}E\rt),
\end{split}\eee
\bee\label{}\begin{split}
 \lt\|  \sa   \mathfrak{L}_{2} (t)\rt\|_{0}^2 \le & C\lt(\lt\|\sa  \pl_t^2 v''\rt\|_{0}  +  \lt\|   \pl_t^2 v' - \pl_t^2 v/x  \rt\|_{0}    + \lt\|  \sa  \mathcal{R}_0\rt\|_{L^\iy }  \lt\| \mathfrak{J}_{2}\rt\|_{0}  + \lt\|  \sa  \mathfrak{L}_0\rt\|_{L^\iy}  \lt\|\mathfrak{J}_{1}\rt\|_{0}\rt.\\
 &\lt. +
 \lt\|  \sa  \mathfrak{L}_1\rt\|_{0}  \lt\|\mathfrak{J}_{0}\rt\|_{L^\iy } \rt)^2
\le   M_0 + CtP\lt(\sup_{[0,t]}E\rt),
\end{split}\eee
where we have used  $\ef{3lk}_{1,3}$.

\noindent Tao Luo\\
Dept. of Mathematics and Statistics,\\
Georgetown University,\\
Washington DC, USA.\\
tl48@georgetown.edu\\
\noindent Zhouping Xin\\
Institute of Mathematical Sciences,\\
Chinese University of Hong Kong,\\
Hong Kong, China. \\
zpxin@ims.cuhk.edu.hk\\
\noindent Huihui Zeng\\
Mathematical Sciences Center,\\
Tsinghua University, \\
Beijing, China. \\
hhzeng@mail.tsinghua.edu.cn


\begin{thebibliography}{10}

\bibitem{AB}G. Auchmuty and R. Beals , Variational solutions
of some nonlinear free boundary problems, Arch. Rat.  Mech.
Anal.43, 255-271 (1971).

\bibitem{1} Ambrose, D., Masmoudi, N.: The zero surface tension limit of three-dimensional water
waves. Indiana Univ. Math. J. 58, 479-521 (2009)

\bibitem{6'} S. Chandrasekhar, Introduction to the Stellar
Structure, {\it University of Chicago Press} (1939).

\bibitem{CL00}
Christodoulou, D., Lindblad, H.: On the motion of the free surface of a liquid.
\newblock Comm. Pure Appl. Math. \textbf{53}(12), 1536--1602 (2000).





\bibitem{7} Coutand, D., Lindblad, H., Shkoller, S.: A priori estimates for the free-boundary
3-D compressible Euler equations in physical vacuum. Commun. Math. Phys. 296,
559-587 (2010)

\bibitem{9} Coutand, D., Shkoller, S.: Well-posedness of the free-surface incompressible Euler
equations with or without surface tension. J. Am. Math. Soc. 20, 829-930 (2007)

\bibitem{10} Coutand, D., Shkoller, S.: Well-posedness in smooth function spaces for the moving-
boundary 1-D compressible Euler equations in physical vacuum. Commun. Pure
Appl. Math. 64, 328-366 (2011)

\bibitem{10'}Coutand, Daniel; Shkoller, Steve; Well-Posedness in Smooth Function Spaces for the Moving-Boundary Three-Dimensional Compressible Euler Equations in Physical Vacuum. Arch. Ration. Mech. Anal. 206 (2012), no. 2, 515-616.
\bibitem{cox}  Cox, J.P., Giuli, R.T.: Principles of stellar structure, I.,II. New York: Gordon and Breach, 1968.

\bibitem{dafermos} C. M. Dafermos, Hyperbolic conservation laws in continuum physics,
Springer-Verlag, Berlin-New York (2005),


\bibitem{DLYY} Deng, Y., Liu, T.P., Yang, T., Yao, Z.: Solutions of Euler-Poisson equations for gaseous stars. Arch. Rat. Mech. Anal. 164(3), 261-285 (2002).

\bibitem{diperna} DiPerna, Ronald J.,  Uniqueness of solutions to hyperbolic conservation laws. Indiana Univ. Math. J. 28 (1979), no. 1, 137-188.

\bibitem{zhenlei} Gu, Xumin; Lei, Zhen Well-posedness of 1-D compressible Euler-Poisson equations with physical vacuum. J. Differential Equations 252 (2012), no. 3, 2160-2188.


\bibitem{15} Jang, J.: Local well-posedness of dynamics of viscous gaseous stars. Arch. Rational
Mech. Anal. 195, 797-863 (2010)
\bibitem{16} Jang, J., Masmoudi, N.:Well-posedness for compressible Euler with physical vacuum
singularity. Commun. Pure Appl. Math. 62, 1327-1385 (2009)
\bibitem{16'} Jang, J., Masmoudi, N.: Well-posedness of compressible Euler equations in a physical vacuum, arXiv:1005.4441.
\bibitem{17'} Jang, J. Nonlinear Instability Theory of Lane-Emden stars, arXiv:1211.2463.
\bibitem{17} Kreiss, H.O.: Initial boundary value problems for hyperbolic systems. Commun. Pure
Appl. Math. 23, 277-296 (1970)
\bibitem{18} Kufner, A.: Weighted Sobolev Spaces. Wiley-Interscience, New York, 1985
\bibitem{18'} A. Kufner, L. Maligranda, L.-E. Persson: The Hardy inequality. Vydavatelsky
Servis, Plzen, 2007. About its history and some related results.
\bibitem{19} Lannes, D.:Well-posedness of the water-waves equations. J. Am. Math. Soc. 18, 605-
654 (2005)
\bibitem{liebyau} Lieb, E.H., Yau, H.T.: The Chandrasekhar theory of stellar collapse as the limit of quantum mechanics. Commun. Math. Phys. 112(1), 147-174 (1987)
\bibitem{20} Lin, L.W.: On the vacuum state for the equations of isentropic gas dynamics. J. Math.
Anal. Appl. 121, 406-425 (1987)
\bibitem{21}Lindblad, H.: Well-posedness for the motion of an incompressible liquid with free
surface boundary. Ann. Math. 162, 109-194 (2005)
\bibitem{22} Lindblad, H.:Well posedness for the motion of a compressible liquid with free surface
boundary. Commun. Math. Phys. 260, 319-392 (2005)
\bibitem{23} Liu, T.-P.: Compressible flow with damping and vacuum. Jpn. J. Appl.Math. 13, 25-32
(1996)
\bibitem{24} Liu, T.-P., Yang, T.: Compressible Euler equations with vacuum. J. Differ. Equ. 140,
223-237 (1997)
\bibitem{25} Liu, T.-P., Yang, T.: Compressible flow with vacuum and physical singularity. Methods
Appl. Anal. 7, 495-310 (2000)
\bibitem{26} Liu, T.-P., Smoller, J.: On the vacuum state for isentropic gas dynamics equations.
Adv. Math. 1, 345-359 (1980).
\bibitem{luosmoller1} Luo, T.; Smoller, J.: Nonlinear dynamical stability of Newtonian rotating and non-rotating white dwarfs and rotating supermassive stars. Comm. Math. Phys. 284 (2008), no. 2, 425-457.
\bibitem{luosmoller2} Luo, T.; Smoller, J.: Existence and non-linear stability of rotating star solutions of the compressible Euler-Poisson equations. Arch. Ration. Mech. Anal. 191 (2009), no. 3, 447-496.
\bibitem{29} Makino, T.: On a local existence theorem for the evolution equation of gaseous stars.
Patterns and Waves. Stud. Math. Appl., Vol. 18. North-Holland, Amsterdam, 459-479,
1986
\bibitem{30} Matusu-Necasova, S., Okada, M., Makino, T.: Free boundary problem for the equation
of spherically symmetric motion of viscous gas III. Jpn. J. Indust. Appl. Math. 14,
199-213 (1997)


\bibitem{33}Okada, M., Makino, T.: Free boundary problem for the equation of spherically symmetric
motion of viscous gas. Jpn. J. Indust. Appl. Math. 10, 219-335 (1993).
\bibitem{rein} Rein, G.: Non-linear stability of gaseous stars. Arch. Rat. Mech. Anal. 168(2), 115-130 (2003)
\bibitem{34} Shatah, J., Zeng, C.: Geometry and a priori estimates for free boundary problems of
the Euler equation. Commun. Pure Appl. Math. 61, 698-744 (2008).


\bibitem{35} Trakhinin, Y.: Local existence for the free boundary problem for the non-relativistic
and relativistic compressible Euler equations with a vacuum boundary condition.
Commun. Pure Appl. Math. 62, 1551-1594 (2009)


\bibitem{36} Wu, S.:Well-posedness in Sobolev spaces of the fullwaterwave problem in 2-D. Invent.
Math. 130, 39-72 (1997)
\bibitem{37} Wu, S.: Well-posedness in Sobolev spaces of the full water wave problem in 3-D.
J. Am. Math. Soc. 12, 445-495 (1999)
\bibitem{38}Xu, C.-J., Yang, T.: Local existence with physical vacuum boundary condition to Euler
equations with damping. J. Differ. Equ. 210, 217-231 (2005)
\bibitem{39} Yang, T.: Singular behavior of vacuum states for compressible fluids. J. Comput. Appl.
Math. 190, 211-231 (2006)

\bibitem{40} Zhang, P., Zhang, Z.: On the free boundary problem of three-dimensional incompressible
Euler equations. Commun. Pure Appl. Math. 61, 877-940 (2008)



\end{thebibliography}
\end{document}